\documentclass[11pt, reqno]{amsart}
\usepackage{amsmath, amssymb, amscd, amsthm, amsxtra}
\usepackage{graphicx}
\usepackage{subfigure}
\usepackage{ragged2e}
\usepackage{color}
\usepackage{comment}
\usepackage{mathtools}
\usepackage{tikz}
\usetikzlibrary{arrows.meta}

\usepackage{todonotes}

\usepackage[pdftex]{hyperref}

\usepackage[format=plain,labelfont=bf,up,width=.9\textwidth]{caption}

\theoremstyle{plain}
\newtheorem{thm}{Theorem}[section]
\newtheorem{lemma}[thm]{Lemma}
\newtheorem{prop}[thm]{Proposition}
\newtheorem{remark}[thm]{Remark}
\newtheorem{cor}{Corollary}[thm]

\theoremstyle{definition}
\newtheorem{defn}[thm]{Definition}

\theoremstyle{plain}
\newtheorem{thmx}{\bf Theorem}

\newcommand{\cal}{\mathcal}

\newcommand{\N}{\mathbb{N}}
\newcommand{\C}{\mathbb{C}}
\newcommand{\Z}{\mathbb{Z}}
\newcommand{\D}{\mathbb{D}}

\newcommand{\s}{\mathbb{S}}

\newcommand{\ds}{\displaystyle}
\newcommand{\HOV}{\rm HOV}

\newcommand{\Crit}{\mathcal{C}} 
\newcommand{\Cs}{\mathcal{C}^{s}} 
\newcommand{\Cu}{\mathcal{C}^{u}} 

\newcommand\restr[2]{{
  \left.\kern-\nulldelimiterspace 
  #1 
  \vphantom{\big|} 
  \right|_{#2} 
  }}

\def\proof{\noindent {\bf Proof.\ }}

\def\qed{\hfill $\square$\\ }

\def\He{{H\'enon }}

\makeatletter
\@namedef{subjclassname@2020}{%
  \textup{2020} Mathematics Subject Classification}
\makeatother

\begin{document}

\title[Critical locus for H\'enon maps in an HOV region]{Critical locus for H\'enon maps\\ in an HOV region}

\author[Tanya Firsova]{Tanya Firsova}
\address{Kansas State University, Manhattan, KS, United States}
\email{tanyaf@math.ksu.edu}

\author{Remus Radu}
\address{Institute of Mathematics of the Romanian Academy, Bucharest, Romania}
\email{rradu@imar.ro}

\author{Raluca Tanase}
\address{Institute of Mathematics of the Romanian Academy, Bucharest, Romania}
\email{rtanase@imar.ro}


\begin{abstract} We prove that the characterization of the critical locus for complex H\'enon maps that are small perturbations of quadratic polynomials with disconnected Julia sets given by Firsova holds in a much larger HOV-like region from the complex horseshoe locus. The techniques of this paper are non-perturbative. 
\end{abstract}

\maketitle

\section{Introduction}\label{sec:Intro}

In dimension one, there is a strong connection between the orbits of the critical points of a polynomial $p$ and the topology of its Julia set. The Julia set is connected if and only if the critical points do not escape to infinity under forward iterations by $p$. When the degree of the polynomial is two, there is the classical dichotomy: the Julia set is either connected or a Cantor set (see e.g. \cite{M1}, \cite{H}).

This dichotomy does not have a proper counterpart in several dimensions, in the study of polynomial automorphisms of $\C^2$. Friedland and Milnor \cite{FM} have classified the polynomial automorphisms of $\C^2$ and shown that the ones with non-trivial dynamics can be reduced to compositions of H\'enon maps with simpler functions. A \He map is defined by $H(x,y)=(p(x)-ay,x)$, where $p$ is a polynomial of degree $d\geq 2$ and $a\neq 0$ is a complex parameter. H\'enon maps have a seemingly easy formula and their dynamics bares some similarity to the dynamics of the polynomial $p$ when the Jacobian $a$ is small. However, their dynamics can be quite complicated in general. As an automorphism of $\C^2$, the H\'enon map has no critical points in the usual sense. Critical loci, sets of tangencies between dynamically defined foliations/laminations often serve as a good analog of the critical points. Several notions of critical points have thus emerged, depending on their location in the dynamical space of the \He map. 

The study of the dynamics of \He maps leads to the definition of the sets $K^{\pm}$ of points which remain bounded in forward/backward time. Let $U^{\pm}=\C^2 - K^{\pm}$ denote the set of points which escape to infinity  under forward, and respectively backward iterations of the \He map. $U^{\pm}$ are called the escaping sets. The sets $J^{\pm}=\partial K^{\pm}=\partial U^{\pm}$, and $J=J^{+}\cap J^{-}$ are called the Julia sets of the \He map. The sets $U^{+}$ and $U^{-}$ are naturally foliated by Riemann surfaces isomorphic to the complex plane, given by the rates at which points escape to infinity (see e.g. \cite{HOV}). The critical locus $\mathcal{C}$ is the set of tangencies between the foliations of $U^+$ and $U^-$. By \cite{BS5}, $\mathcal{C}$  is a closed nonempty analytic subvariety of $U^+\cap U^-$, which may have singularities, and is invariant under the \He map. When $J^+$ or $J^-$ are laminated by copies of the stable, respectively unstable manifolds of points in $J$, one can define other critical loci: $\mathcal{C}^s$, the set of tangencies between the foliation of $U^-$ and the lamination of $J^+$, and respectively  $\mathcal{C}^u$, the set of tangencies between the foliation of $U^+$ and the lamination of $J^-$. Bedford and Smillie \cite{BS6} showed that important dynamical and topological properties such as Lyapunov exponents, or the connectivity of the Julia set $J$ are related to the stable/unstable critical loci $\mathcal{C}^{s}/\mathcal{C}^{u}$. For a given H\'enon map, either of the sets $\mathcal{C}^{s}$ or $\mathcal{C}^{u}$ may be empty, however $\mathcal{C}$ is never empty. Therefore, understanding the geometric and topological properties of the critical locus $\mathcal{C}$, as well as its relation with the other critical loci, can be very relevant for the dynamics of the H\'enon map.

Homoclinic/heteroclinic tangencies between the stable and unstable manifolds of some saddle periodic point(s) are probably the most known types of ``critical points'' as they are the basis of many coexisting phenomena in the H\'enon family, which differentiate the dynamics of the H\'enon map from  one-dimensional dynamics (see e.g. Palis, Takens \cite{PT},  Newhouse \cite{N}). When the \He map is hyperbolic, the laminations of $J^{+}$ and $J^{-}$ are  transverse to each other, thus there are no ``critical points'' in $J$. 

In some cases, there is also a notion of critical points in the interior of $K^{+}$. For moderately dissipative H\'enon maps where $int(K^{+})$ contains the basin of attraction of an attracting or semi-parabolic cycle, Dujardin and Lyubich \cite{DL} show the existence of ``critical points'' in these basins (as tangencies between the foliations of the basins and the unstable manifolds in $J^{-}$). In our setting we have no such critical points as the interior of $K^{+}$ is empty.

So far, an explicit model for the critical locus $\mathcal{C}$ was given only for special classes of \He maps with small Jacobian. In \cite{LR}, Lyubich and Robertson described the critical locus for \He maps that are small perturbations of degree $d\geq 2$ polynomials with simple critical points and with connected Julia set. We state their result for degree two:
\begin{thm}[Lyubich, Robertson  \cite{LR}]\label{thm: LyubichRobertson}
\justifying  Let $H$ be a quadratic \He map which is a small perturbation of a hyperbolic polynomial
$p(x)=x^2+c$ with connected Julia set.
There exists a unique primary component $\mathcal{C}_0$ of the {\it
critical locus} asymptotic to the $x$-axis. $\mathcal{C}_0$ is biholomorphic to $\C-\overline{\D}$, and it is everywhere transverse to the foliations of $U^+$ and $U^-$. All other connected components of the critical locus are forward or backward iterates of $\mathcal{C}_0$ under the \He map $H$.
\end{thm}

The same model holds true for quadratic \He maps with a semi-parabolic periodic point that are small perturbations of quadratic polynomial with a parabolic fixed point \cite{T}. By \cite{RT1} and \cite{RT2}, these \He maps and some nearby hyperbolic perturbations have connected Julia set $J$. In these cases, the critical locus was used in an essential way in \cite{T} to describe the Julia set $J^{+}$ in terms of discrete groups acting on $\s^{1}\times \C$.

In \cite{F}, the first author gave an explicit description of the critical locus for \He maps that are small perturbations of quadratic polynomials with disconnected Julia set. In this paper we prove that the same model holds true throughout a large subset of the complex horseshoe region. To our knowledge, these are the first results on the structure of the critical locus in a non-perturbative setting.

Unlike the quadratic polynomials, \He maps with disconnected Julia sets can have a wild behavior, and a general understanding of their dynamics is currently out of reach. The region
\begin{equation}\label{eq:HOV}
\HOV = \left\{ (c,a)\in \C^{2}: |c|>2(1+|a|)^{2} \right\}
\end{equation}
is a subset of the complex horseshoe region and was introduced by Hubbard and Oberste-Vorth to study a part of the parameter space of \He maps with disconnected Julia set (see \cite{O}, \cite{MNTU}, \cite{BS}). The boundary of the HOV region intersects the Mandelbrot set in the $c$-plane at the tip $c=-2$. 
When $(c,a)\in\HOV$, the \He map 
\begin{equation}\label{eq:Hca}
H_{c,a}(x,y)=(x^2+c-ay,x),\ a\neq0
\end{equation}
 is hyperbolic on $J$ and conjugate to the standard Smale horseshoe map. In this context, $J^{\pm}=K^{\pm}$, $J$ is homeomorphic to a Cantor set, and the \He map on $J$ is topologically conjugate to the full shift on two symbols (see \cite{MNTU}, \cite{FS} in perturbative setting). 
 
Consider now a real parameter $\beta$ and define
\begin{equation}\label{eq:HOVb}
\HOV_{\beta} = \left\{ (c,a)\in \C^{2}: |c|>\beta(1+|a|)^{2} \right\}
\end{equation}
a subset of the HOV region. 
While it is clear that the optimal lower bound on $\beta$ for a set defined by the algebraic conditions \eqref{eq:HOVb} to be part of the horseshoe region is $\beta=2$, higher values of $\beta$ were also considered in the literature. 
 The region HOV$_{\beta}$ with $\beta=(5+2\sqrt{5})/4$ was analyzed by Oberste-Vorth in \cite{O}, where increasing $\beta$ had certain advantages, for example being able to work with standard families of vertical/horizontal cones when showing hyperbolicity on $J$.

In the theorems that follow, we give a comprehensive description of the critical locus and its boundary for \He maps in the HOV$_{\beta}$, where $\beta=18.75$. We expect the same model outlined in Theorems \ref{thm:CL}, \ref{thm:SM} and \ref{thm:Acc} to hold throughout the HOV region. However, we choose to work in the $\HOV_{\beta}$ region for technical reasons, to relax some of the estimates needed to trap the critical locus.  

\begin{thmx}\label{thm:CL}  For all $(c,a)\in \HOV_{\beta}$ the critical locus of $H_{c,a}$ is a smooth, irreducible complex analytic subvariety of $U^{+}\cap U^{-}$, of pure dimension one.
\end{thmx}

We now use the dynamics of the polynomial to define the main building block of the construction: {\it the truncated spheres}. We refer to Milnor \cite{M1} and Hubbard \cite{H} for a treatment of one-dimensional complex dynamics. In Section \ref{sec:criticalModel} we show how they naturally arise in the dynamics of the H\'enon map.

\begin{figure}[htb]
\begin{center}
\includegraphics[scale=0.165]{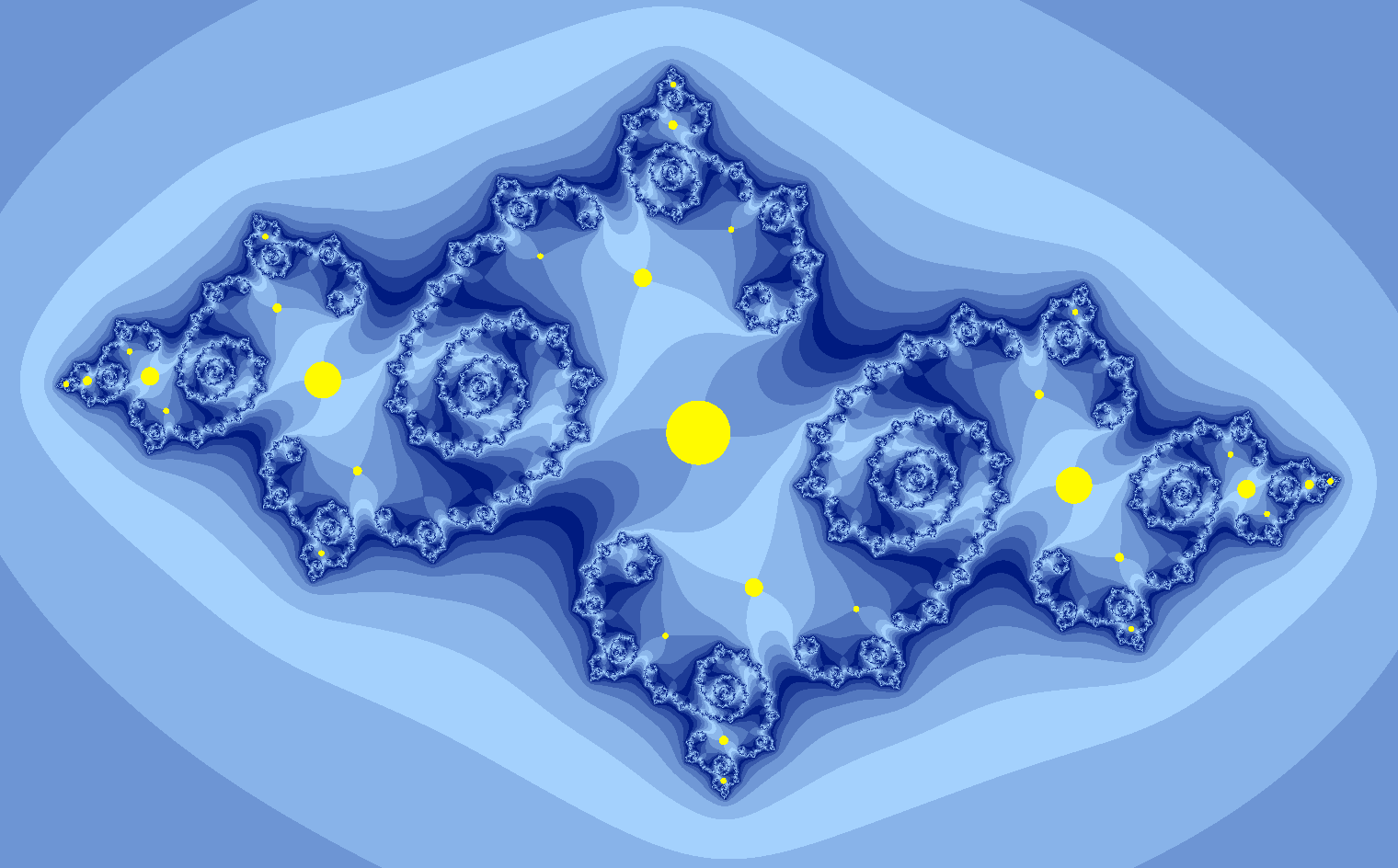}
\end{center}
\caption{A Cantor Julia set of a quadratic polynomial together with the removed disks $W_{\alpha_{n}}$.}
\label{fig:Cantor}
\end{figure}

Consider a quadratic polynomial $p(x)=x^{2}+c$ such that $c$ belongs to the complement of the Mandelbrot set in the plane. Let $\Sigma^+=\{0,1\}^{\N}$. The Julia set $J_p$ is a Cantor set, and the polynomial $p$ on $J_p$ is conjugate to the shift map $\sigma:\Sigma^+\rightarrow\Sigma^+$, $\sigma(\alpha_0\alpha_1\ldots)=\alpha_1\alpha_2\ldots$. Each point of the Julia set $J_{p}$ of $p$ is uniquely characterized by a one-sided infinite sequence $\alpha$ of $0$'s and $1$'s, so we identify $J_p$ with $X=\{x_{\alpha}, \alpha\in \Sigma^{+}\}$. We consider an open topological disk $D$ bounded by an equipotential of $J_p$ 
such that $p:p^{-1}(D)\rightarrow D$ is a degree two covering map. Denote by $D_0$ and $D_1$ the two components of $p^{-1}(D)$.  $D_0$ and $D_1$ are topological disks, compactly contained in $D$. We remove from $D$ a small disk $W:=W_{\emptyset}$ around the critical point $0$. We then remove two more disks $W_0$ from $D_0$ and $W_1$  from $D_1$, centered at the preimages of $0$ (see Figure \ref{fig:Cantor}). We proceed inductively. For each finite sequence $\alpha_n$ of $0$'s and $1$'s, of length $n>0$, we remove  from the corresponding preimage $D_{\alpha_n}$ of $D$ under $p^{-n}$ a small disk $W_{\alpha_n}$, centered at the corresponding preimage $p^{-n}(0)$ and compactly contained in $\C\setminus J_p$.
The sequence $\overline{D_{\alpha_n}}$ is a nested sequence of compact sets with diameters converging to $0$, therefore it converges to a point $x_{\alpha}\in J_p$ as $\alpha_n$ converges to $\alpha\in \Sigma^+$. As a consequence, the sequence 
of removed disks $W_{\alpha_n}$ also converges to $x_{\alpha}$.  

\begin{figure}[htb]
\begin{center}
\includegraphics[scale=0.4175]{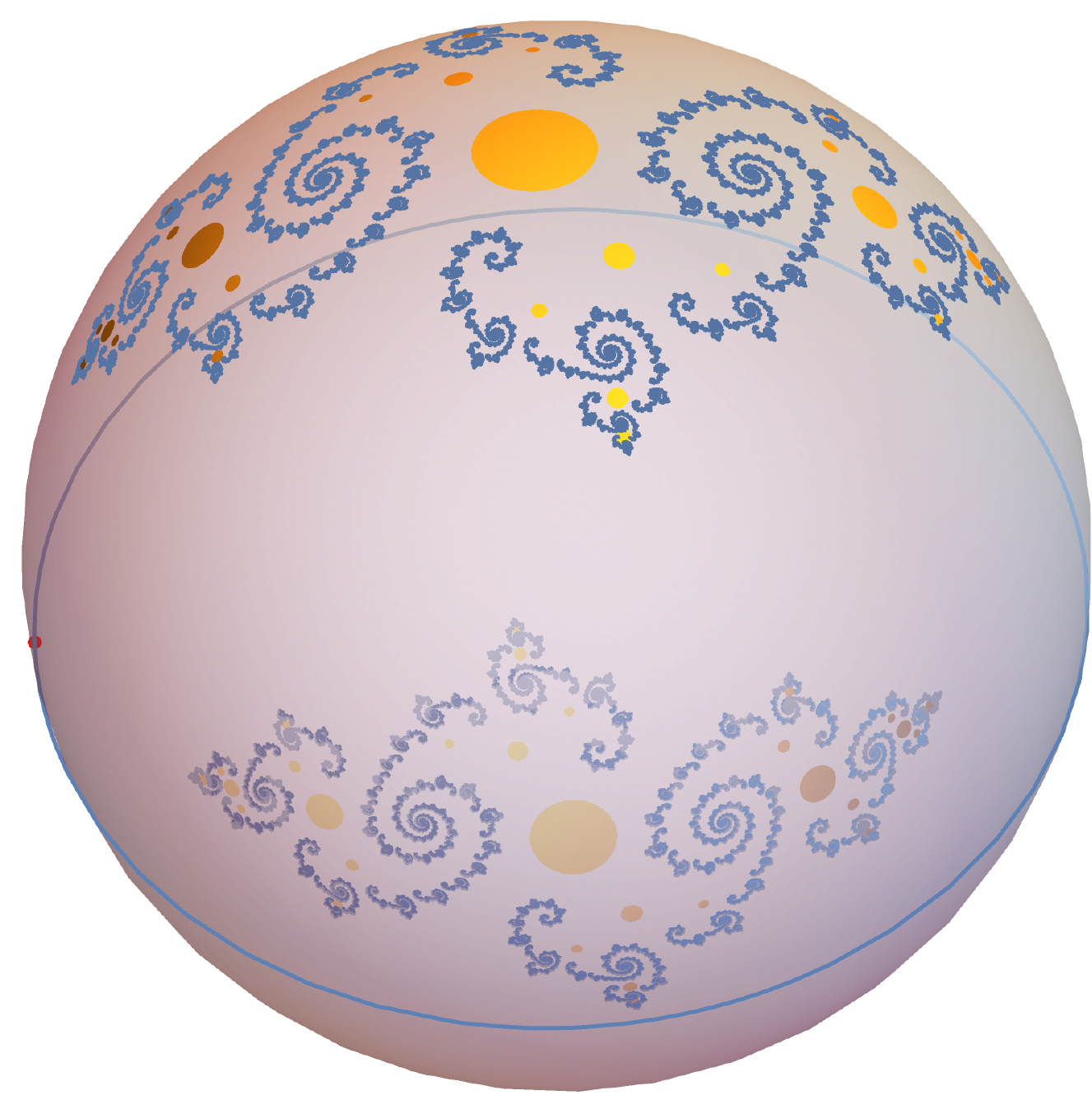}
\end{center}
\caption{A truncated sphere.}
\label{fig:truncated}
\end{figure}

Let $S^{+}$ denote the disk $\overline{D}$ from which we  remove the disks $W_{\alpha_n}$ and the Julia set $X$. We glue $S^{+}$ to itself along the  boundary of $D$ and remove an extra point from the equator $\partial D$. The resulting set is a {\it truncated sphere}, which we denote $\mathcal{S}$ (see Figure \ref{fig:truncated}). The set $S^+$ lies in the upper hemisphere. To distinguish between the two hemispheres, we denote the lower hemisphere by $S^{-}$, the removed disks from $D$ by $U_{\alpha_n}$, where $\alpha_{n}\in \{0,1\}^{n}$, $n\geq 0$, and the removed Cantor set by $\Omega$. Notice that the truncated spheres are all homeomorphic, regardless of the choice of the parameter $c$ from the outside of the Mandelbrot set. 

\begin{thmx}\label{thm:SM}  
For all $(c,a)\in \HOV_{\beta}$ the critical locus of $H_{c,a}$ is homeomorphic to the following Riemann surface:
consider countably many truncated spheres $(\mathcal{S}_{k})_{k\in \Z}$ such that for each $n,k\in\Z$, $n\geq 0$, and each $\alpha_{n}\in\{0,1\}^{n}$, the boundary of $W_{\alpha_{n}}$ on the sphere $\mathcal{S}_{k}$ is glued to the boundary of $U_{\alpha_{n}}$ on $\mathcal{S}_{k+n+1}$ by a handle $\Upsilon_{\alpha_{n}}^{k}$. The \He map $H_{c,a}$ acts on the critical locus by sending the sphere $\mathcal{S}_{k}$ to $\mathcal{S}_{k+1}$, and the handle $\Upsilon_{\alpha_{n}}^{k}$ to $\Upsilon_{\alpha_{n}}^{k+1}$. 
\end{thmx}

The model of the critical locus given by Theorem \ref{thm:SM} was conjectured by John Hubbard. 

\begin{thmx}\label{thm:Acc} In the $\HOV_{\beta}$ region, the accessible boundary of the critical locus $\mathcal{C}$ is $\Cs\cup\Cu$. The boundary of the critical locus $\mathcal{C}$ is $J^{+}\cup J^{-}$.
\end{thmx}

By analyzing the cases in \cite{LR}, \cite{F}, \cite{T} and the $\HOV$-like regions in this paper, where we understand the critical locus, a relation between the topology of the Julia set $J$ and the critical locus has emerged: The Julia set $J$ is connected iff the critical locus $\mathcal{C}$ is disconnected. 
We conjecture this to be true for all dissipative (partially) hyperbolic quadratic complex H\'enon maps.

\medskip
\noindent \textit{Acknowledgements.} TF was partly supported by the Simons Collaboration Grant 965460.
RR was partly supported by grant PN-III-P4-ID-PCE-2020-2693 from the Romanian Ministry of Education and Research, CNCS-UEFISCDI. RT was partly supported by grant PN-III-P1-1.1-TE-2019-2275 from the Romanian Ministry of Education and Research, CNCS-UEFISCDI. The authors participated in a program hosted by the MSRI, Berkeley during the Spring 2022 semester (supported by the NSF under Grant No. DMS-1928930) and are grateful to MSRI for the hospitality.

\section{Dynamical filtration}\label{sec:filtration}

When $a\neq 0$, the \He map $H_{c,a}$ from Equation \eqref{eq:Hca} is invertible with polynomial inverse $H_{c,a}^{-1}(x,y)=(y, (y^2+c-x)/a)$. For simplicity, we denote $H_{c,a}$ by $H$. 
Choose any constant $R$ such that 
\begin{equation}\label{eq:right}
R>R_{1}:=\frac{1+|a|+\sqrt{(1+|a|)^{2}+4|c|}}{2}.
\end{equation}

\begin{figure}[htb]
\begin{center}
\includegraphics[scale=1.1]{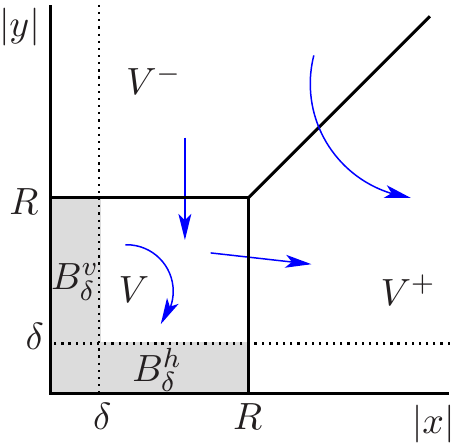}
\end{center}
\caption{Dynamical filtration of $\C^2$.}
\label{fig:filtration}
\end{figure}

As in \cite{HOV}, the dynamical spaces of the H\'enon map can be divided into three sets: 
\begin{eqnarray}\label{eq:VV+V-}
V &=& \left\{(x,y)\in \C^{2}: |x|\leq R, |y|\leq R\right\}\\ \nonumber
V^{-} &=& \left\{(x,y)\in \C^{2}: |y|>R, |x|\leq |y|\right\}\\ \nonumber
V^{+} &=& \left\{(x,y)\in \C^{2}: |x|>R, |x|\geq |y|\right\} 
\end{eqnarray}
One has $H(V^+)\subset V^+$ and $H^{-1}(V^-)\subset V^-$. In addition $H(V)\subset V\cup V^+$, and $H^{-1}(V)\subset V\cup V^-$. Any point in $V^{\mp}$ will be  eventually mapped in $V\cup V^{\pm}$ under forward (respectively backward) iterates. Therefore, the escaping set $U^+$ can be
described as the union of backward iterates of $V^+$. Similarly, the set $U^-$ is the union of the forward iterates of $V^-$ under the \He map. The sets $J^+$ and $K^+$ are contained in $V\cup V^+$, whereas the sets $J^-$ and $K^-$ are contained in $V\cup V^-$. Lastly, the sets $J$ and $K$ are compact sets contained in the polydisk $V$. A partition of $\C^{2}$ like in Equation \eqref{eq:VV+V-} satisfying the properties above is called a Hubbard filtration of $\C^{2}$ (see Figure \ref{fig:filtration}).

The escaping sets $U^{+}$ and $U^{-}$ are naturally foliated by Riemann surfaces isomorphic to $\C$. These foliations, which we denote by ${\cal F}^+$ and ${\cal F}^-$, can be
characterized in terms of the B\"{o}ttcher coordinates $\varphi^{+}$ and $\varphi^{-}$, or in terms of the rate of escape functions $G^{+}$ and $G^{-}$. These functions were constructed in \cite{HOV} and we will summarize their properties. The B\"{o}ttcher coordinate $\varphi^{+}$ is a holomorphic map defined on $V^{+}$ with values in $\C-\D$, which semiconjugates the \He map $H$ to $z\mapsto z^2$ and behaves like the projection on the first coordinate as $(x,y)\mapsto \infty$ in $V^{+}$. Analogously, the B\"{o}ttcher coordinate $\varphi^{-}$ is a complex-valued holomorphic map defined on $V^{-}$  which semiconjugates the inverse \He map $H^{-1}$  to $z\mapsto z^2/a$
and behaves like the projection on the second coordinate as we approach infinity in $V^{-}$. 
Unlike the one-dimensional case where the Bottcher coordinates can be extended as holomorphic maps to the entire escaping sets, in two dimensions the topology of the escaping sets is more complicated and does not permit a holomorphic extension of $\varphi^{\pm}$ to $U^{\pm}$. The level sets
of $\varphi^{\pm}$ define the holomorphic foliation ${\cal F}^{\pm}$ in $V^{\pm}$ which can be propagated by
dynamics to all of $U^{\pm}$, using the level sets of the functions $(\varphi^{+})^{2^{k}}$ and $(\varphi^{-})^{2^{k}}$. These dyadic powers are well defined and holomorphic on $H^{-\circ k}(V^{+})$, respectively on $H^{\circ k}(V^{-})$ for all for $k\geq 0$, which can be seen using the inductively obtained relations 
\begin{eqnarray}\label{eq:dyadic}
&\left(\varphi^{+}\right)^{2^{k}} = \varphi^{+}\circ H^{\circ k}\ \ & \mbox{ on  }\  H^{-\circ k}(V^{+}), \nonumber\\
&\left(\varphi^{-}/a\right)^{2^{k}} = (\varphi^{-}/a)\circ H^{-\circ k} &\mbox{ on }\  H^{\circ k}(V^{-}). 
\end{eqnarray}
Naively, one would like to think of these foliations as a coordinate system in
$U^+\cap U^-$. However, for any H\'{e}non map there exists a
codimension one subvariety consisting of tangencies between ${\cal F}^+$
and ${\cal F}^-$ \cite{BS5}, which we refer to as the critical locus $\mathcal{C}$. 

Equation \eqref{eq:dyadic} implies that $\varphi^{+}$ is well defined on $ U^{+}$ up to a local choice of a dyadic root of unity. This makes $\log \varphi^{+}$ well defined on $U^{+}$ up to a local addition of an additive constant, therefore $d\log\varphi^{+}$ is a holomorphic $1$-form on $U^{+}$. In a similar fashion, by equation \eqref{eq:dyadic}, we get that $d\log\varphi^{-}$ is a holomorphic $1$-form on $U^{-}$.
The critical locus $\mathcal{C}$ can also be described as  the set of zeroes of the holomorphic function $w(x,y)$, where 
 \begin{equation}\label{eq:definingfunction}
 d\log \varphi^{+}\wedge d\log \varphi^{-}=w(x,y)dx\wedge dy.
 \end{equation}

Throughout the paper, when there is no danger of confusion, we will use the notations $H^{k}$ and $H^{-k}$ in place of $H^{\circ k}$ and $H^{-\circ k}$, to denote forward, respectively backward iterations of the H\'enon map.

Suppose that $(c,a)\in \HOV$. Since $|c|>2(1+|a|)^{2}$, one can choose a constant $R$ which satisfies simultaneously the lower bound $R_{1}$ from Equation \eqref{eq:right}, and the upper bound
\begin{equation}\label{eq:R}
\frac{|c|}{1+|a|}>R>R_{1}.
\end{equation}

In the $\HOV$ region, the H\'enon map expands horizontally and contracts vertically, and satisfies the 2-fold horseshoe condition with respect to the bidisk $V$: each of the sets $V\cap H^{-1}(V)$ and $V\cap H(V)$ consists of two connected components, biholomorphic to $V$. In fact, Equation \eqref{eq:R} implies that all iterates $H^{\pm n}$ satisfy the $2^{n}-$fold horseshoe condition with respect to $V$. The sequences $\left\{V\cap H^{-n}(V)\right\}_{n\geq 0}$ and $\left\{V\cap H^{n}(V)\right\}_{n\geq 0}$ are nested sequences with $2^{n}$ connected components of decreasing widths, respectively heights, and 
\begin{equation}\label{eq: horseshoe}
V\cap J^{+}=\bigcap\limits_{n=0}^{\infty} H^{-n}(V), \qquad  V\cap J^{-}=\bigcap\limits_{n=0}^{\infty} H^{n}(V).
\end{equation}
 The $2^{n }$ connected components of $V\cap H^{-n}(V)$ can be labeled inductively using finite strings with $n$ letters from the alphabet $\{0,1\}$. This induces a matching labelling on the $2^{n }$ connected components of $V\cap H^{n}(V)$, since $H^{n}\left( V\cap H^{-n}(V)\right)=V\cap H^{n}(V)$.

Let $\Sigma=\{0,1\}^{\Z}$. The Julia set $J$ is a Cantor set, and the map $H$ on $J$ is conjugate to the full shift on $2$ symbols $\sigma:\Sigma\rightarrow \Sigma$, $\sigma(\alpha)=\alpha'$, where $\alpha'_n=\alpha_{n+1}$, for all $n\in\Z$.

Throughout the paper, we set
\begin{equation}\label{eq:delta}
\delta := \frac{1}{2}(|a|+1).   
\end{equation}
Let $S^v_{\delta}$ and $S^h_{\delta}$ denote the vertical and horizontal strips of width $\delta$:
\begin{equation*}
S^v_{\delta}=\{(x,y)\in\C^2, |x|\leq \delta\}, \qquad  S^h_{\delta}=\{(x,y)\in\C^2, |y|\leq \delta\}.
\end{equation*}
Define three subsets of $V$, as in Figure \ref{fig:filtration}: 
\begin{equation}
B^{v}_{\delta} = V\cap S^v_{\delta}, \qquad  B^{h}_{\delta} = V\cap S^h_{\delta}, \qquad  B_{\delta}:=B^{v}_{\delta} \cap B^{h}_{\delta}
\end{equation}
Consider also the subsets $V^{\pm}_{\delta} \subset V^{\pm}$ given by
\begin{equation*}
V^{+}_{\delta} = V^{+}\setminus S^h_{\delta}, \qquad  V^{-}_{\delta} = V^{-}\setminus S^v_{\delta}.
\end{equation*}
One can immediately  verify that $H^{\pm}(V^{\pm}) \subset V^{\pm}_{\delta}$ since $R>\delta$. 

It is easy to see that $H(S^v_{\delta})=S^h_{\delta}$. Moreover, points from $S^h_{\delta} \setminus V$ are mapped by $H$ in $V^{+}\cap\{|y|>R\}$. Similarly,  points from $S^v_{\delta} \setminus V$ are mapped by $H^{-1}$ in $V^{-} \cap\{|x|>R\}$. 

Proposition \ref{prop:Bdelta} implies that $B_{\delta}^{h}$ is disjoint from the sets $H(V)\cap V$ and $H(V)\cap S_{\delta}^h$. Similarly, $B_{\delta}^{v}$ is disjoint from the sets $H^{-1}(V)\cap V$ and $H^{-1}(V)\cap S_{\delta}^v$.

\begin{prop}\label{prop:Bdelta} For all $(c,a)\in HOV_{\beta}$ there exists $R>0$ satisfying Equation \eqref{eq:R} such that  $H(B_{\delta}^{v})\subset V^{+}\cap S^h_{\delta}$ and $H^{-1}(B_{\delta}^{h})\subset V^{-}\cap S^v_{\delta}$.
\end{prop}
\proof
Let $|x|\leq\delta$ and $|y|\leq R$ and write $\delta=\alpha(|a|+1)$ as in \eqref{eq:delta}. We wish to prove that 
\[
|x^{2}+c-ay|\geq  |c|-|x|^{2}-|a||y|\geq |c|-\alpha^{2}(|a|+1)^{2}-|a|R>R,
\]
which is equivalent to 
\[
\frac{|c|}{|a|+1}-\alpha^{2}(|a|+1)>R.
\]
We can find a value of $R$ to verify this inequality and  \eqref{eq:R} provided that 
\[
\frac{|c|}{|a|+1}-\alpha^{2}(|a|+1)>\frac{1+|a|+\sqrt{(1+|a|)^{2}+4|c|}}{2},
\]
which is equivalent to 
\[
|c|^{2}-2(\alpha^{2}+1)(|a|+1)^{2}|c|+ \alpha^{2}(\alpha^{2}+1)(|a|+1)^{4}>0.
\]
The last inequality holds because $|c|>\beta(1+|a|)^{2}$ in the $\HOV_{\beta}$ region. Note that 
for $\alpha=1/2$ we have $\alpha^{2}+1+\sqrt{\alpha^{2}+1}=\beta=(5+2\sqrt{5})/4$ and the choice of $\alpha$ is optimal. This shows that there exists $R>0$ such that $H(B_{\delta}^{v})\subset V^{+}$. The computations for $B_{\delta}^{h}$ are similar.
\qed

\begin{defn}\label{def:vcones} Define the {\it horizontal cone} at a point $(x,y)\in \C^{2}$ as
\[
\mathcal{C}^{h}_{(x,y)}=\left\{(\xi,\eta)\in
T_{(x,y)}\C^{2},\ |\xi|\geq|\eta|\right\}
\]
and the {\it vertical cone} at a point $(x,y)\in \C^{2}$  as
\[
\mathcal{C}^{v}_{(x,y)}=\left\{(\xi,\eta)\in
T_{(x,y)}\C^{2},\ |\xi|\leq|\eta|\right\}.
\]
\end{defn}

\noindent We will use the vector norm $||(\xi,\eta)||:= \mbox{max} (|\xi|,|\eta|)$. We denote by $Int\ \mathcal{C}^{v/h}_{(x,y)}$ the topological interior of the cone $\mathcal{C}^{v/h}_{(x,y)}$ together with the vector $\{0\}$.

\begin{prop}[\textbf{Invariant families of cones}]\label{prop:cones}$\ $
 
\noindent If $(x,y)\in \C^2\setminus S^h_{\delta}$ then 
 $DH^{-1}_{(x,y)}(\mathcal{C}^{v}_{(x,y)})\subset Int\ \mathcal{C}^{v}_{H^{-1}(x,y)}$ and 
\[
\|DH^{-1}_{(x,y)}(w)\| > \frac{2|y|-1}{|a|}\|w\|, \ \ w\in \mathcal{C}^{v}_{(x,y)}.
\]
If $(x,y)\in \C^2\setminus S^v_{\delta}$ then $DH_{(x,y)}(\mathcal{C}^{h}_{(x,y)})\subset Int\ \mathcal{C}^{h}_{H(x,y)}$ and 
\[
\|DH_{(x,y)}(w)\| >(2|x|-|a|) \|w\|, \ \ w\in \mathcal{C}^{h}_{(x,y)}.
\]
\end{prop}
\proof Choose any point $(x,y)$ with $|y|> \delta$ and a vector $(\xi,\eta)$ with $|\eta|\geq |\xi|$. Let $(\xi',\eta')=DH^{-1}_{(x,y)}(\xi,\eta)=(\eta, -\xi/a +2y\eta/a)$. We can estimate 
\begin{eqnarray*}
|\eta'|\geq (2|y||\eta|-|\xi|)/|a|\geq k|\eta|=k|\xi'|,
\end{eqnarray*}
where $k=(2|y|-1)/|a|>1$. This inequality implies that $(\xi',\eta')\in Int\ \mathcal{C}^{v}_{H^{-1}(x,y)}$ and at the same time $\|DH^{-1}_{(x,y)}(\xi,\eta)\| > \|(\xi,\eta)\|(2|y|-1)/|a|$.

Choose any point $(x,y)$ with $|x|> \delta$ and a vector $(\xi,\eta)$ with $|\xi|\geq |\eta|$. Define $(\xi',\eta')=DH_{(x,y)}(\xi,\eta)=(2x\xi-a\eta, \xi)$. We can estimate 
\begin{eqnarray*}
|\xi'|\geq 2|x||\xi|-|a||\eta|\geq k|\xi|=k|\eta'|,
\end{eqnarray*}
where $k=(2|x|-|a|)>1$. Therefore $(\xi',\eta')\in Int\ \mathcal{C}^{h}_{H(x,y)}$ and 
\[\|DH_{(x,y)}(\xi,\eta)\| > \|(\xi,\eta)\|(2|x|-|a|),
\]
which concludes the proof.
\qed

\begin{defn} We will call a foliation $\mathcal{F}$ horizontal/vertical-like if for every point $x$ of any leaf $\mathcal{L}_x$ of the foliation, the tangent space to the leaf $T_x(\mathcal{L}_x)$ is contained in the horizontal/vertical cone  $\mathcal{C}^{h}_x/\mathcal{C}^{v}_x$.
\end{defn}

\begin{prop}\label{prop:foliation1} The foliation of $U^+$ in the set $U^+\setminus  \bigcup_{n\geq 0} H^{-n}(S_{\delta}^v) $ is vertical-like. The foliation of $U^-$ inside the set $U^-\setminus \bigcup_{n\geq 0} H^{n}(S_{\delta}^h) $ is horizontal-like.
\end{prop}
\proof 
The function $\varphi^+$ behaves like a projection on the first coordinate $\varphi^+(x,y)\simeq x$ when $(x,y)\rightarrow \infty$ in $V^+$. Therefore, for $|x|$ large enough, the level set $(\varphi^{+})^{-1}(x)$ is a vertical-like holomorphic disk in $V^+$. By Proposition \ref{prop:cones}, on $V^+_{\delta}$ we have an invariant family of vertical cones: if $x\in V^+_{\delta}$ then $DH^{-1}_x\left(\mathcal{C}^v_x\right)\subset Int\ \mathcal{C}^v_{H^{-1}(x)}$. Since $H(V^+)\subset V^+_{\delta}$, we can conclude that the foliation is vertical-like throughout the entire set $V^+$.

In fact, by Proposition \ref{prop:cones}, we can iterate the vertical-like foliation of $V^+$ backward and ensure that it remains vertical-like for as long as the backward iterates do not enter the vertical strip $S^v_{\delta}$. The escaping set $U^+$ is given by the relation $U^+=\bigcup_{n\geq 0} H^{-n}(V^+)$. Therefore the foliation of $U^+$ is vertical-like throughout the set $U^+\setminus \bigcup_{n\geq 0} H^{-n}(S_{\delta}^v)$.

Similarly, one can show that the foliation of $U^-$ is horizontal-like throughout the set $U^-\setminus \bigcup_{n\geq 0} H^{n}(S_{\delta}^h)$.
\qed

\begin{cor}\label{cor:foliation1} The foliation of $U^+$ in the set $V^+\cup V\setminus \bigcup_{n\geq 0} H^{-n}(B_{\delta}^v) $ is vertical-like. The foliation of $U^-$ inside the set $V^-\cup V \setminus \bigcup_{n\geq 0} H^{n}(B_{\delta}^h) $ is horizontal-like.
\end{cor}
\proof 
By the proof of Proposition \ref{cor:foliation1}, the foliation of $U^+$  is vertical-like on $V^+$. Also, the foliation of $U^+$ is vertical-like in the set $U^+\setminus \bigcup_{n\geq 0} H^{-n}(S_{\delta}^v)$. The set $S_{\delta}^v\setminus B_{\delta}^v$ is a subset of $V^-$, so all backward images of $S_{\delta}^v\setminus B_{\delta}^v$ under $H^{-1}$ are contained in $V^-$. Therefore $\bigcup_{n\geq 0} H^{-n}(S_{\delta}^v) \cap V\subset \bigcup_{n\geq 0} H^{-n}(B_{\delta}^v)$. Hence, the foliation of $U^+$ is vertical-like in the region $V\setminus \bigcup_{n\geq 0} H^{-n}(B_{\delta}^v)$.
\qed

We will end this chapter of preliminary observations with a lemma about intersections of the form $H^{n}(B^{h}_{r'})\cap H^{-m}(B^{v}_{r})$, where $0<r,r'<R$, $B^{v}_{r}=\{|x|\leq r, |y|\leq R\}$, $B^{h}_{r'}=\{|x|\leq R, |y|\leq r'\}$, and $m,n\geq0$, and show that they are biholomorphic to a standard polydisk. Denote by  $(h^{1}_{k}(x,y), h^{2}_{k}(x,y))=H^{k}(x,y)$, $k\neq 0$, the components of the forward or backward iterates of the H\'enon map.

\begin{lemma}[\textbf{Polydisks}]\label{lemma:polydisks} For every positive integer $n\geq1$ we have 
\begin{eqnarray}
H^{n}(B^{h}_{r'})\cap B^{v}_{r} &=& \{(x,y)\in\C^{2}, |x|<r \mbox{ and } |h^{2}_{-n}(x,y)|\leq r'\} \label{Eq:poly} \\
H^{-n}(B^{v}_{r})\cap B^{h}_{r'} &=& \{(x,y)\in\C^{2}, |y|<r' \mbox{ and } |h^{1}_{n}(x,y)|\leq r\}.
\end{eqnarray}
Each connected component of the set $H^{n}(B^{h}_{r'})\cap B^{v}_{r} $, respectively of the set $H^{-n}(B^{v}_{r})\cap B^{h}_{r'}$, is biholomorphic to the standard polydisk $\D_{r}\times \D_{r'}$ via the map $G(x,y)=(x, h^{2}_{-n}(x,y))$, and respectively via the map $F(x,y)=(h^{1}_{n}(x,y), y)$.
\end{lemma}
\proof It suffices to show only the first identity. Let $A$ and $B$ denote the sets from the left hand side and the right hand side of Equation \eqref{Eq:poly}.  It is easy to check that $A\subseteq B$, since the set $A$ is given by the conditions 
\[
	\{(x,y)\in\C^{2},\ |x|<r,\, |y|<R,\, |h^{1}_{-n}(x,y)|<R,\, |h^{2}_{-n}(x,y)|<r' \}.
\] 
To show that $B\subseteq A$, pick now any point $(x,y)\in B$. Since $|x|<r<R$, it follows that $(x,y)\in V\cup V^{-}$. By the properties of the dynamical filtration of $\C^{2}$, we must also have $H^{-n}(x,y)\in V\cup V^{-}$. By hypothesis, $|h^{2}_{-n}(x,y)|\leq r'<R$, which implies that $H^{-n}(x,y)\in V$, hence automatically $H^{-n}(x,y)\in B^{h}_{r'}$. Also, since $H^{-n}(x,y)\in V$ we get  $(x,y)=H^{n}(H^{-n}(x,y))\in V\cup V^{+}$. However, we already knew that $(x,y)\in V\cup V^{-}$, therefore $(x,y)\in V$. Combining this with the inequalities from the definition of the set $B$, we get that $(x,y)\in B^{v}_{r}$ and so $B\subseteq A$.

In the HOV region, all iterates $H^{n}$ satisfy the $2^{n}-$fold horseshoe condition with respect to the bidisk $V$, a condition which is equivalent to $\partial_{y}h^{2}_{-n}(x,y)\neq 0$ on $H^{n}(V)\cap V$ (see \cite{MNTU}, Chapter 7).  For example, when $n=1$, it is easy to see that 
\[
\partial_{y}h^{2}_{-1}(x,y)=2y/a\neq 0
\] 
on $H(V)\cap V$ since this set omits the $x$-axis. 
The function $G:A\rightarrow \D_{r}\times \D_{r'}$, $G(x,y)=(x,h^{2}_{-n}(x,y))$ is a covering map when $n\geq 1$, since the Jacobian of $G$ is equal to $\partial_{y}h^{2}_{-n}$ which never vanishes on the set $A\subset H^{n}(V)\cap V$. The Monodromy Theorem states that any connected cover above a simply connected domain is single-sheeted, hence each of the $2^{n}$ connected components of the set $A$ is mapped biholomorphically by $G$ onto the polydisk  $\D_{r}\times \D_{r'}$. 
\qed

\section{Finding a suitable HOV$_{\beta}$-region}\label{sec:DynamicalBox}

We will outline a strategy to find regions HOV$_{\beta}$ such that the polydisk $V$ can be replaced by a dynamical polydisk $\widetilde{V}$ whose boundaries consist of leaves of the foliation of $U^{+}$ and $U^{-}$, the set $B^{v}_{\delta}$ can be replaced by a set $B^{v}$ whose vertical boundary is foliated by the leaves of the foliation of $U^{+}$, and similarly  $B^{h}_{\delta}$ can be replaced by $B^{h}$ whose horizontal boundary is foliated by the leaves of the foliation of $U^{-}$.

We revisit Relation \eqref{eq:R} and note that it gives information about how $H^{-1}(V)\cap V$ and $H(V)\cap V$ fit inside $V$. Set 
\begin{equation}\label{eq:eta}
\eta_{1}=\sqrt{|c|+R(1+|a|)}\ \ \mbox{and}\ \  \eta_{2}=\sqrt{|c|-R(1+|a|)}.
\end{equation}

\begin{prop}
The distance between the set $H^{-1}(V)\cap V$ and the vertical boundary $\{|x|=R\}$ of $V$, respectively the distance between the set $H(V)\cap V$ and the horizontal boundary $\{|y|=R\}$ of $V$ is greater than or equal to $\eta_{1}$. 

\noindent The distance between $H^{-1}(V)\cap V$ and the vertical line $\{x=0\}$, respectively the distance between $H(V)\cap V$ and the horizontal line $\{y=0\}$ is greater than or equal to $\eta_{2}$. 
\end{prop}

\proof The proof follows by direct computations, but we include the details below, for completion.
Fix any complex number $y_{0}\in\C$ with $|y_{0}|\leq R$, and look at the set of points 
$(x_{0},y_{0}) \in V\cap H^{-1}(V)$, for which $(x_{0},y_{0})=H^{-1}(x,y)$ for some $(x,y)\in V$ with $|x|=R$.
Since $x_{0}=y$ and $y^{2}+c-x=ay_{0}$,  an upper bound on $|x_{0}|=|y|$ can be obtained by
\[
|y^{2}|-|c|-R\leq|y^{2}+c-x|\leq|a|R
\]
which implies $|y|\leq \sqrt{|c|+(1+|a|)R}$. Note that this is a non-trivial bound on $y$ because $\sqrt{|c|+(1+|a|)R}< R$, 
which follows from Relation \eqref{eq:right} because $R>R_{1}$, where $R_{1}$ is the largest root of the quadratic polynomial $R^{2}-(1+|a|)R-|c|$.
Under the same assumptions, a lower bound on $|y|$ can be deduced from the inequality
\[
|c|-R-|y|^{2}\leq|y^{2}+c-x|\leq|a|R
\]
which implies that $|y|\geq \sqrt{c-(1+|a|)R}$. Note that Relation \eqref{eq:R} ensures that the expression under the square root is positive. Hence
\[
\eta_{2}\leq |y|\leq \eta_{1}
\]
Fix any point $x_{0}\in\C$ with $|x_{0}|\leq R$, and look at the set of points 
$(x_{0},y_{0})\in V\cap H(V)$, for which $(x_{0},y_{0})=H(x,y)$ for some $(x,y)\in V$ with $|y|=R$.
Since $x=y_{0}$ and $x^{2}+c-ay=x_{0}$,  an upper bound on $|y_{0}|=|x|$ can be obtained from the inequality
\[
|x^{2}|-|c|-|a|R\leq|x^{2}+c-ay|\leq R.
\]
A lower bound on $|x|$ can be found from the inequality
\[
|c|-|a|R-|x|^{2}\leq|x^{2}+c-ay|\leq R.
\]
This leads to the same bounds as above: $\eta_{2}\leq |x|\leq \eta_{1}$.
\qed

\begin{prop}[\textbf{Technical estimates}]\label{prop:technical} There exist constants $1>\gamma_{1}>\gamma_{2}>0$ and $\beta>2$ such that in the parametric region HOV$_{\beta}$ there exists $R$ satisfying \eqref{eq:R} and the inequality 
 \[
 0<\gamma_{2}R\leq \eta_{2}<\eta_{1}\leq \gamma_{1}R<R.
 \]
 \end{prop}
 
\proof From $\eta_{1}=\sqrt{|c|+R(1+|a|)}\leq \gamma_{1}R$ and $\eta_{2}=\sqrt{|c|-R(1+|a|)}\leq \gamma_{2}R$ we get 
\[
\gamma_{2}^{2}R^{2}+R(1+|a|)-|c|\leq 0 \leq \gamma_{1}^{2}R^{2}-R(1+|a|)-|c|,
\]
which yields a lower and an upper bound on $R$:
\begin{equation}\label{eq:2}
\frac{1+|a|+\sqrt{(1+|a|)^{2} +4|c|\gamma_{1}^{2}}\, }{2\gamma_{1}^{2}}\leq R\leq 
\frac{-(1+|a|)+\sqrt{(1+|a|)^{2} +4|c|\gamma_{2}^{2}}\, }{2\gamma_{2}^{2}}.
\end{equation}
To simplify the computations, let $c_{0}=4|c|/(1+|a|)^{2}$.
The last inequality becomes:
\begin{equation}\label{eq:Rlowup}
\frac{1+|a|}{2}\left(\frac{1+\sqrt{1+c_{0}\gamma_{1}^{2}}}{\gamma_{1}^{2}}\right)\leq R\leq 
\frac{1+|a|}{2}\left(\frac{-1+\sqrt{1+c_{0}\gamma_{2}^{2}}}{\gamma_{2}^{2}}\right),
\end{equation}
which gives 
\begin{equation*}
\frac{1+\sqrt{1+c_{0}\gamma_{1}^{2}}}{\gamma_{1}^{2}}\leq 
\frac{-1+\sqrt{1+c_{0}\gamma_{2}^{2}}}{\gamma_{2}^{2}}
\end{equation*}
or equivalently 
\begin{equation}\label{eq:aux1}
\frac{1}{\gamma_{1}^{2}}+\frac{1}{\gamma_{2}^{2}}\leq 
-\frac{\sqrt{1+c_{0}\gamma_{1}^{2}}}{\gamma_{1}^{2}}+
\frac{\sqrt{1+c_{0}\gamma_{2}^{2}}}{\gamma_{2}^{2}}.
\end{equation}
Note that the function $x\mapsto \sqrt{1+c_{0}x^{2}}/x^{2}$ 
is decreasing and $\gamma_{1}>\gamma_{2}$, so the right hand side of \eqref{eq:aux1} is positive. After squaring both sides of \eqref{eq:aux1} and rearranging the terms we obtain 
\begin{equation}\label{eq:aux2}
2\sqrt{1+c_{0}\gamma_{1}^{2}}\sqrt{1+c_{0}\gamma_{2}^{2}}\leq c_{0}(\gamma_{1}^{2}+\gamma_{2}^{2})-2.
\end{equation}
We ask that the right hand side of inequality \eqref{eq:aux2} is positive, which means requiring that $ |c|>\frac{1}{2(\gamma_{1}^{2}+\gamma_{2}^{2})}(1+|a|)^{2}$.
After squaring inequality \eqref{eq:aux2} we obtain
\[
8c_{0}(\gamma_{1}^{2}+\gamma_{2}^{2})\leq c_{0}^{2}(\gamma_{1}^{2}-\gamma_{2}^{2})^{2},
\]
and further on
\begin{equation*}
|c_{0}|>\frac{8(\gamma_{1}^{2}+\gamma_{2}^{2})}{(\gamma_{1}^{2}-\gamma_{2}^{2})^{2}},
\end{equation*}
which places a new HOV-type constraint on $c$, of the form:
\begin{equation*}
|c|>\frac{2(\gamma_{1}^{2}+\gamma_{2}^{2})}{(\gamma_{1}^{2}-\gamma_{2}^{2})^{2}}(1+|a|)^{2}.
\end{equation*}
Clearly $\frac{2(\gamma_{1}^{2}+\gamma_{2}^{2})}{(\gamma_{1}^{2}-\gamma_{2}^{2})^{2}}>\frac{1}{2(\gamma_{1}^{2}+\gamma_{2}^{2})}$, so we can choose $\beta=\frac{2(\gamma_{1}^{2}+\gamma_{2}^{2})}{(\gamma_{1}^{2}-\gamma_{2}^{2})^{2}}$, and $\gamma_{1}>\gamma_{2}$.

We need to address the size of $R$ in the HOV$_{\beta}$ region.  Relation \eqref{eq:2} gives the lower bound
\begin{equation*}
 R \geq \frac{1+|a|}{2}\cdot \frac{1+\sqrt{1+c_{0}\gamma_{1}^{2}}}{\gamma_{1}^{2}}
 > 
\frac{1+|a|}{2}\cdot \frac{1+ \sqrt{1+\frac{8\gamma_{1}^{2}(\gamma_{1}^{2}+\gamma_{2}^{2})}{(\gamma_{1}^{2}-\gamma_{2}^{2})^{2}}}}{\gamma_{1}^{2}}.
\end{equation*}
Since $(\gamma_{1}^{2}-\gamma_{2}^{2})^{2}+8\gamma_{1}^{2}(\gamma_{1}^{2}+\gamma_{2}^{2})=(3\gamma_{1}^{2}+\gamma_{2}^{2})^{2}$, we find that
\begin{equation}\label{eq:Rlower}
R>\frac{2}{\gamma_{1}^{2}-\gamma_{2}^{2}}(1+|a|).
\end{equation}
Note also that
\[
\frac{|c|}{1+|a|}>\beta (1+|a|)>\frac{2}{\gamma_{1}^{2}-\gamma_{2}^{2}}(1+|a|),
\]
where the latter inequality is verified since $\gamma_{1}^{2}+\gamma_{2}^{2}>\gamma_{1}^{2}-\gamma_{2}^{2}$ holds true regardless of the values of $\gamma_{1},\gamma_{2}>0$. So a choice of $R$ which verifies Relation \eqref{eq:R} is  possible.
\qed

\begin{lemma}[\textbf{Projections of fibers}]\label{lemma:proj} Let $R_{2}, R_{3}$ be such that $\delta<R_{2}<R_{3}\leq R$.
\begin{enumerate}
\item
Consider the set
\[
W_{1}=\left\{ (x,y)\in V\setminus H^{-1}(V),\, R_{2}\leq |x|\leq R_{3},\ |y|\leq R\right\}.
\]
The foliation of $U^{+}$ in $W_{1}$ is vertical-like. If
\[
\frac{2R|a|}{2R_{2}-1}<R_{3}-R_{2},
\]
then there are no connected components of fibers of the foliation of $U^{+}$ in $V$ which pass through both vertical boundaries $|x|=R_{3}$ and $|x|=R_{2}$ of the set $W_{1}$. 
\item
Consider the set
\[
W_{2}=\left\{ (x,y)\in V\setminus H(V),\, R_{2}\leq |y|\leq R_{3},\ |x|\leq R\right\}.
\]
The foliation of $U^{-}$ in $W_{2}$ is vertical-like. If
\[
\frac{2R|a|}{2R_{2}-1}<R_{3}-R_{2},
\]
then there are no connected components of the fibers of the foliation of $U^{-}$ in $V$ which pass through both horizontal boundaries $|y|=R_{3}$ and $|y|=R_{2}$ of $W_{2}$. 
\end{enumerate}
\end{lemma} 
\proof
Since $W_{1}$ is disjoint from $H^{-1}(V)$, it is automatically disjoint from all the backward images $\bigcup_{n\geq 0}H^{-n}(B^{v}_{\delta})$, so the leaves of the foliation of $U^{+}$ inside $W_{1}$ are vertical-like. Similarly the set $W_{2}$ is disjoint from $H(V)$, so it is automatically disjoint from $\bigcup_{n\geq 0}H^{n}(B^{h}_{\delta})$, thus the leaves of the foliation of $U^{-}$ inside $W_{2}$ are horizontal-like. 

Let $L(y)=(\varphi(y),y),\ |y|<R$, be a connected component of a vertical-like leaf of the foliation of $U^{+}$ in $V$ which intersects the boundary $|x|=R_{3}$ at some point $(x_{0},y_{0})$. Let $\gamma(t)=\varphi(ty_{0}+(1-t)y_{1})$ be any path in $L$, $t\in[0,1]$, connecting $(x_{0},y_{0})$ to some other point $(x_{1},y_{1})$ in $L$.
 We will show that the length of $\gamma$ is strictly less than $R_{3}-R_{2}$, which in turn implies that the diameter of the horizontal projection of the leaf $L$ is strictly less that $R_{3}-R_{2}$, so $L$ cannot cross both boundaries $|x|=R_{2}$ and $|x|=R_{3}$.
Using Proposition \ref{prop:cones} we have
\begin{eqnarray}\label{eq:aux3}
\ell(\gamma) &=& \int_{0}^{1}|\varphi'(ty_{0}+(1-t)y_{1})||y_{1}-y_{0}|dt \leq |y_{1}-y_{0}|\int_{0}^{1}\frac{|a|}{2|\varphi(ty_{0}+(1-t)y_{1})|-1}dt\nonumber \\ 
&\leq&|y_{1}-y_{0}|\frac{|a|}{2R_{2}-1}\leq 2R\frac{|a|}{2R_{2}-1}<R_{3}-R_{2}
\end{eqnarray}

Let $L(x)=(x,\psi(x)),\ |x|<R$, be a horizontal-like leaf of the foliation of $U^{-}$ which intersects the boundary $|y|=R_{3}$ at some point $(x_{0},y_{0})$. Let $\gamma(t)=\psi(tx_{0}+(1-t)x_{1})$ be any path in $L$, $t\in[0,1]$, connecting $(x_{0},y_{0})$ to some other point $(x_{1},y_{1})$ in $L$.
 We will show that the length of $\gamma$ is strictly less than $R_{3}-R_{2}$, which in turn implies that the diameter of the vertical projection of the leaf $L$ is strictly less that $R_{3}-R_{2}$, so $L$ cannot cross both boundaries $|y|=R_{2}$ and $|y|=R_{3}$.
Using Proposition \ref{prop:cones} we have
\begin{eqnarray}\label{eq:aux4}
\ell(\gamma) &=& \int_{0}^{1}|\psi'(tx_{0}+(1-t)x_{1})||x_{1}-x_{0}|dt \nonumber\\
&\leq& |x_{1}-x_{0}|\int_{0}^{1}\frac{1}{2|\psi(tx_{0}+(1-t)x_{1})|-|a|}dt\nonumber\\
&\leq&|x_{1}-x_{0}|\frac{1}{2R_{2}-|a|}\leq \frac{2R}{2R_{2}-|a|}<R_{3}-R_{2},
\end{eqnarray}
 which concludes the proof. \qed

 We will actually apply Lemma \ref{lemma:proj} for more specific leaves of the foliation of $U^{+}$ and $U^{-}$, for which the following corollary is better suited:
 \begin{cor}\label{cor:proj}  Let $\delta<R_{2}<R_{3}\leq R$, and $W_{1}, W_{2}$ be as in Lemma \ref{lemma:proj}. 
 \begin{enumerate}
\item
If
\[
\frac{R|a|}{2R_{2}-1}<R_{3}-R_{2}
\]
then there are no fibers of the foliation of $U^{+}$ which pass through both $|x|=R_{2}$ and $\{|x|=R_{3}, y=0\}$, or through $|x|=R_{3}$ and $\{|x|=R_{2}, y=0\}$. 

\item If
\[
\frac{R|a|}{2R_{2}-1}<R_{3}-R_{2}
\]
then there are no fibers of the foliation of $U^{-}$ which pass through the sets $|y|=R_{2}$ and $\{|y|=R_{3}, x=0\}$, or through $\{|y|=R_{2}, x=0\}$  and $|y|=R_{3}$.
\end{enumerate}
  \end{cor}
\proof The only changes to make in the proof of Lemma \ref{lemma:proj} is in Equations \eqref{eq:aux3} and \eqref{eq:aux4}. We now have $y_{0}=0$, and respectively $x_{0}=0$, so $R$ can be used as an upper bound for $|y_{1}-y_{0}|$ and $|x_{1}-x_{0}|$, in place of $2R$.
 \qed

 Consider $\gamma_{1}$, $\gamma_{2}$, $\beta$ and $R$ as in Proposition \ref{prop:technical}. Let $V_{1}=\left\{|x|\leq \gamma_{1}R,\ |y|\leq \gamma_{1}R\right\}$.

 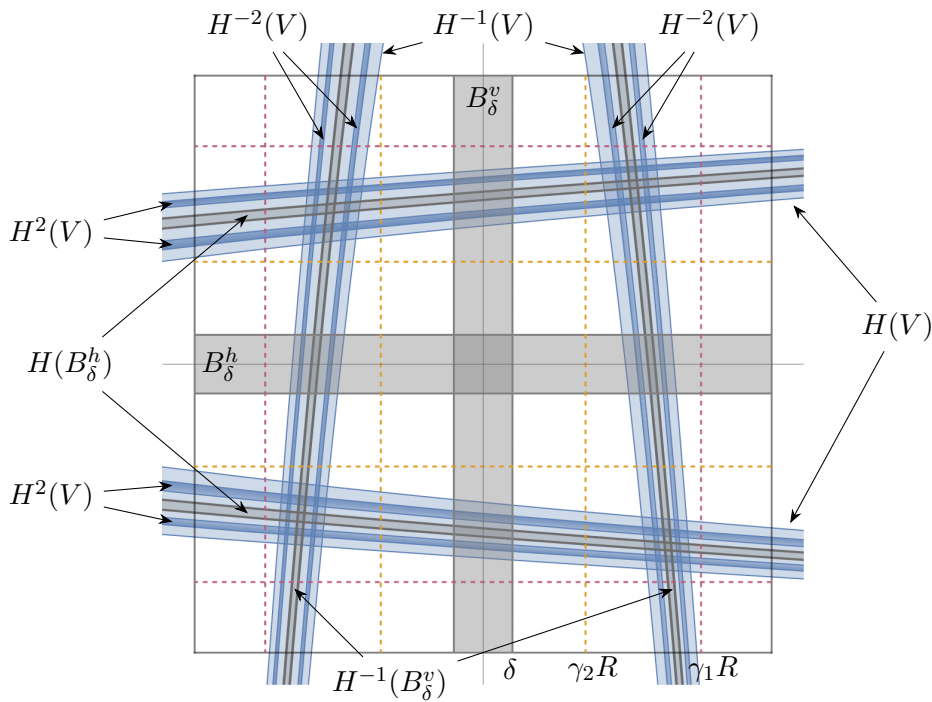
\begin{figure}[htb]
\begin{center}
\begin{tikzpicture}
\node[] (image) at (0,0) {\includegraphics[scale=.45]{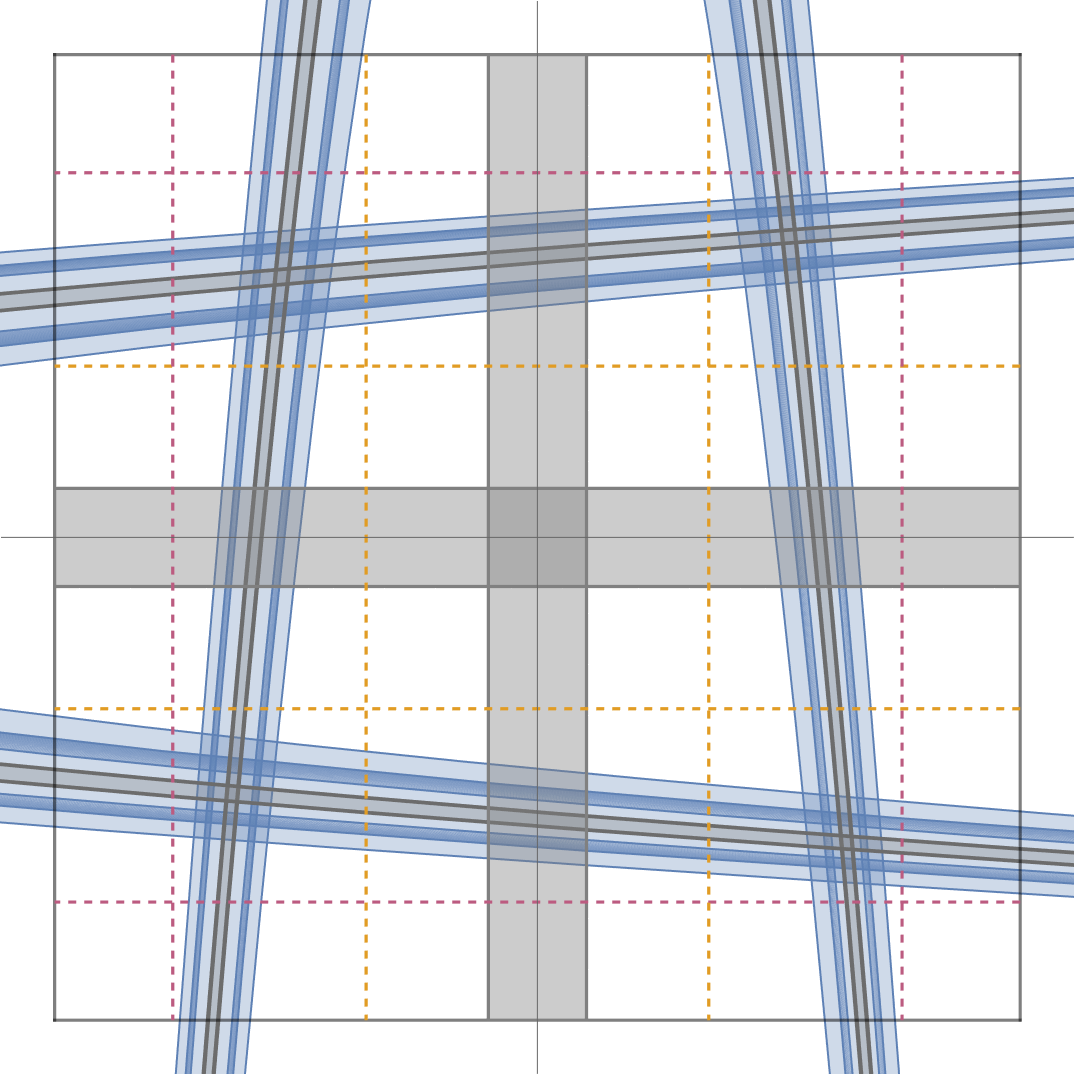}};
\node[] (a) at (1.45, -4.05) {$\gamma_{2}R$};
\node[] (b) at (3.05, -4.05) {$\gamma_{1}R$};
\node[] (c) at (0.35, -4.05) {$\delta$};
\node[] (d) at (-5.75, 1.75) {$H^{2}(V)$};
\node[] (dA) at (-4, 2.175) {};
\node[] (dB) at (-4, 1.575) {};
\node[] (e) at (5.5, 0.5) {$H(V)$};
\node[] (eA) at (4, 2.25) {};
\node[] (eB) at (4, -2.25) {};
\node[] (f) at (-5.75, -1.75) {$H^{2}(V)$};
\node[] (fA) at (-4, -2.125) {};
\node[] (fB) at (-4, -1.6) {};
\node[] (g) at (0, 4.5) {$H^{-1}(V)$};
\node[] (gA) at (-1.485, 4.1) {};
\node[] (gB) at (1.475, 4.1) {};
\node[] (h) at (-3.5, 0) {$B_{\delta}^{h}$};
\node[] (i) at (0, 3.5) {$B_{\delta}^{v}$};
\node[] (j) at (-5.5, 0) {$H(B_{\delta}^{h})$};
\node[] (jA) at (-3, 2.1) {};
\node[] (jB) at (-3, -2.1) {};
\node[] (k) at (-3, 4.5) {$H^{-2}(V)$};
\node[] (kA) at (-1.55, 3) {};
\node[] (kB) at (-2.05, 2.85) {};
\node[] (l) at (3, 4.5) {$H^{-2}(V)$};
\node[] (lA) at (1.575, 3) {};
\node[] (lB) at (2.03, 2.85) {};
\node[] (m) at (-1.25, -4.25) {$H^{-1}(B_{\delta}^{v})$};
\node[] (mA) at (-2.655, -2.85) {};
\node[] (mB) at (2.65, -2.95) {};
\draw[-{Stealth[scale=1.1]}] (e)--(eA);
\draw[-{Stealth[scale=1.1]}] (e)--(eB);
\draw[-{Stealth[scale=1.1]}] (d)--(dA);
\draw[-{Stealth[scale=1.1]}] (d)--(dB);
\draw[-{Stealth[scale=1.1]}] (f)--(fA);
\draw[-{Stealth[scale=1.1]}] (f)--(fB);
\draw[-{Stealth[scale=1.1]}] (g)--(gA);
\draw[-{Stealth[scale=1.1]}] (g)--(gB);
\draw[-{Stealth[scale=1.1]}] (j)--(jA);
\draw[-{Stealth[scale=1.1]}] (j)--(jB);
\draw[-{Stealth[scale=1.1]}] (k)--(kA);
\draw[-{Stealth[scale=1.1]}] (k)--(kB);
\draw[-{Stealth[scale=1.1]}] (l)--(lA);
\draw[-{Stealth[scale=1.1]}] (l)--(lB);
\draw[-{Stealth[scale=1.1]}] (m)--(mA);
\draw[-{Stealth[scale=1.1]}] (m)--(mB);
\end{tikzpicture}
\end{center}
\caption{A typical picture for the set $V$ and its forward/backward images (shaded in blue) in the HOV$_{\beta}$ region. The lines $\pm\gamma_{2}R$ are shown in orange, the lines $\pm\gamma_{1}R$ in red. The tubes $B_{\delta}^{v}$,  $H^{-1}(B_{\delta}^{v})$, $B_{\delta}^{h}$,  $H(B_{\delta}^{h})$ are shaded in grey.}
\label{fig:horseshoe}
\end{figure}

\begin{figure}[htb]
\begin{center}
\includegraphics[scale=1]{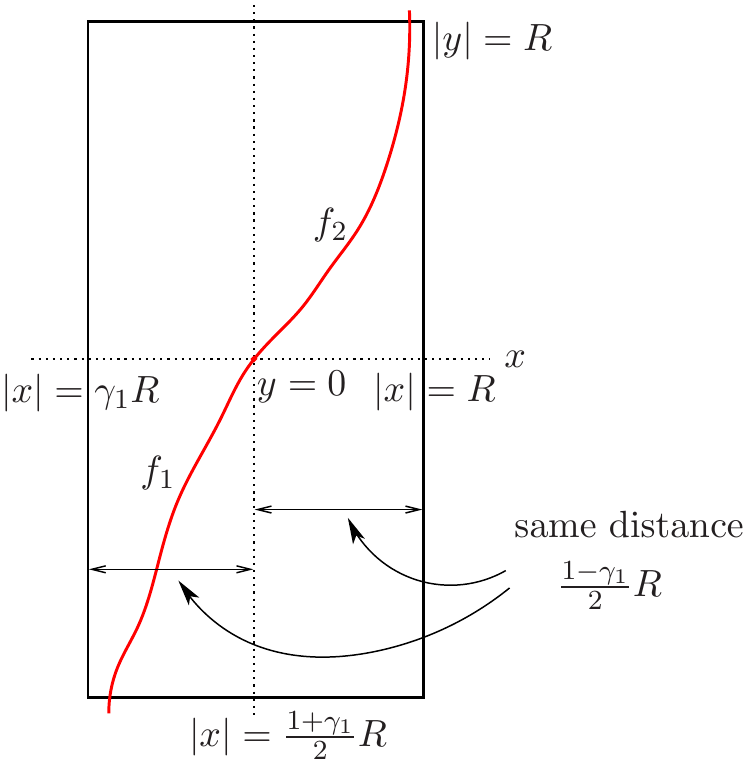}
\end{center}
\caption{A fiber that exits through the horizontal boundary.}
\label{fig:fibers1B}
\end{figure}

\begin{figure}[htb]
\begin{center}
\includegraphics[scale=1]{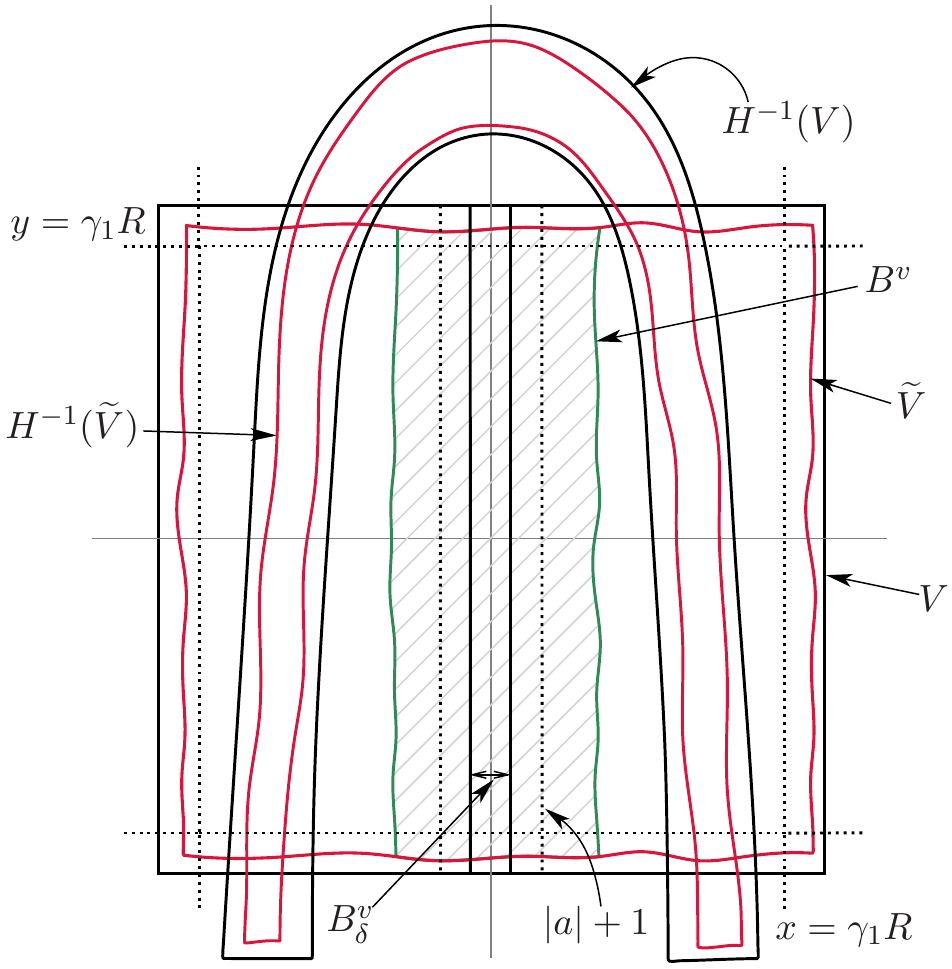}
\end{center}
\caption{A picture of the set $\widetilde{V}$ and the inner dynamical tube $B^{v}$. 
}
\label{fig:VandTube}
\end{figure}

 \begin{prop}[\textbf{Dynamical $V$}]\label{prop:dynV} 
 Suppose that $$\gamma_{2}<\gamma_{1}\leq \frac{1+\sqrt{1+3\gamma_{2}^{2}}}{3}.$$ There exists a set $\widetilde{V}$ with the same dynamical properties of the set $V$ from the Hubbard filtration of $\C^{2}$ such that $V_{1}\subset \widetilde{V}\subset V$, and moreover, the vertical-like boundary of $\widetilde{V}$ is foliated by leaves of the foliation of $U^{+}$, respectively the horizontal-like boundary of $\widetilde{V}$ is foliated by leaves of the foliation of $U^{-}$.  
 
 See Figure \ref{fig:VandTube} for a schematic drawing of the dynamical set $\widetilde{V}$.
 \end{prop}
 \proof We first discuss how to create the vertical boundary of the desired set $\widetilde{V}$. Let us choose the connected components of the leaves of the foliation of $U^{+}$ inside $V$ which intersect the circle $|x|=\frac{1+\gamma_{1}}{2}R$, $y=0$ (see Figure \ref{fig:fibers1B}). We will apply Corollary \ref{cor:proj}
 to show that these components intersect neither $|x|=\gamma_{1}R$, nor $|x|=R$. Indeed, we first restrict to the polydisk $$W_{1}=\left\{\gamma_{1}R\leq|x|\leq \frac{1+\gamma_{1}}{2}R,\ |y|\leq R\right\}$$ and verify that the conditions of Corollary \ref{cor:proj} are met. 
  
 We have
 \[
 \frac{R|a|}{2\gamma_{1}R-1}<R-\frac{1+\gamma_{1}}{2}R\, \Leftrightarrow\,  R>\frac{2|a|+(1-\gamma_{1})}{2\gamma_{1}(1-\gamma_{1})}
 \]
 By Inequality \eqref{eq:Rlower}, we also have $\ds R>\frac{2}{\gamma_{1}^{2}-\gamma_{2}^{2}}(1+|a|)$, so we need a relation between $\gamma_{1}$ and $\gamma_{2}$ so that 
 \[
 \frac{2}{\gamma_{1}^{2}-\gamma_{2}^{2}}(1+|a|)\geq \frac{2|a|+(1-\gamma_{1})}{2\gamma_{1}(1-\gamma_{1})},
 \]
 which is equivalent to 
 \[
 4\gamma_{1}(1-\gamma_{1})|a|+4\gamma_{1}(1-\gamma_{1})\geq 2(\gamma_{1}^{2}-\gamma_{2}^{2})|a|+(1-\gamma_{1})(\gamma_{1}^{2}-\gamma_{2}^{2}).
 \]
  Note that since $\gamma_{1}<1$, the free terms already satisfy the inequality
 \[
 4\gamma_{1}(1-\gamma_{1})>\gamma_{1}^{2}(1-\gamma_{1})>(1-\gamma_{1})(\gamma_{1}^{2}-\gamma_{2}^{2}).
 \]
 It is therefore convenient to ask that the coefficients of $|a|$ verify $4\gamma_{1}(1-\gamma_{1})\geq 2(\gamma_{1}^{2}-\gamma_{2}^{2})$. This reduces to imposing that $3\gamma_{1}^{2}-2\gamma_{1}-\gamma_{2}^{2}\leq 0$, which is satisfied for  $\gamma_{1}\leq \frac{1+\sqrt{1+3\gamma_{2}^{2}}}{3}$.
 
In conclusion, Corollary \ref{cor:proj} implies that a fiber which passes through $|x|=\frac{1+\gamma_{1}}{2}R$, $y=0$ does not cross $|x|=\gamma_{1}R$.
 
 We look at the second polydisk $$W_{1}'=\left\{\frac{1+\gamma_{1}}{2}R\leq|x|\leq R,\ |y|\leq R\right\}$$ and note that 
 \[
  \frac{R|a|}{2\frac{1+\gamma_{1}}{2}R-1}=  \frac{R|a|}{(1+\gamma_{1})R-1}<\frac{R|a|}{(2\gamma_{1})R-1}<
\frac{1-\gamma_{1}}{2}R.
 \]
The last inequality is exactly what we discussed above. By Corollary \ref{cor:proj} for the polydisk $W_{1}'$ we know that the connected components of the fibers of the foliation of $U^{+}$ inside $V$ which pass through $|x|=\frac{1+\gamma_{1}}{2}R$, $y=0$  do not cross $|x|=R$ either. Note that this implies that they will only exit the polydisk $V$ through the horizontal boundary $|y|=R$.

\begin{figure}[htb]
\begin{center}
\includegraphics[scale=1]{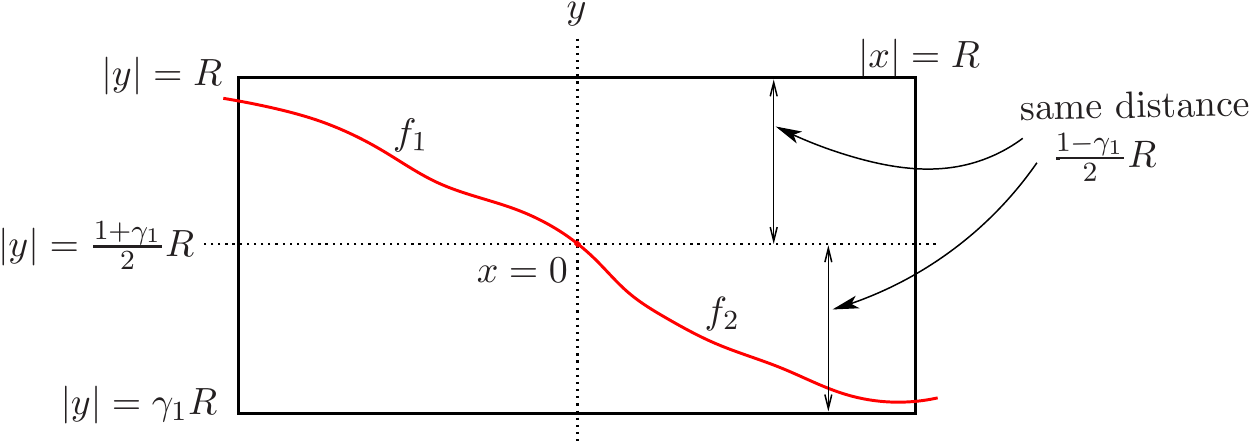}
\end{center}
\caption{A fiber that exits through the vertical boundary.}
\label{fig:fibers2B}
\end{figure}

We will now focus on constructing horizontal-like boundaries of the set $\widetilde{V}$. We  look at the connected components of the fibers of the foliation of $U^{-}$ inside $V$ which intersect the circle $|y|=\frac{1+\gamma_{1}}{2}R$, $x=0$ (see Figure \ref{fig:fibers2B}). In order to apply Corollary \ref{cor:proj}
 to show that they intersect neither $|y|=\gamma_{1}R$, nor $|y|=R$, we need to verify two conditions:
 \[
 \frac{R}{2\gamma_{1}R-|a|}<\frac{1-\gamma_{1}}{2}R\ \ \  \mbox{and}\ \ \  \frac{R}{2\frac{1+\gamma_{1}}{2}R-|a|}<\frac{1-\gamma_{1}}{2}R.
 \]
The second inequality is weaker than the first, so we just need to analyze the first one, which is equivalent to
\[
R>\frac{2+|a|(1-\gamma_{1})}{2\gamma_{1}(1-\gamma_{1})}.
\]
By \eqref{eq:Rlower}, we only need to argue that
\[ 
R>\frac{2}{\gamma_{1}^{2}-\gamma_{2}^{2}}(1+|a|)>\frac{2+|a|(1-\gamma_{1})}{2\gamma_{1}(1-\gamma_{1})},
\]  
which is equivalent to 
 \[
 4\gamma_{1}(1-\gamma_{1})|a|+4\gamma_{1}(1-\gamma_{1})\geq (1-\gamma_{1)}(\gamma_{1}^{2}-\gamma_{2}^{2})|a|+2(\gamma_{1}^{2}-\gamma_{2}^{2}).
 \]
Since $0<\gamma_{1}<1$, the coefficients of $|a|$ satisfy the inequality
 \[
 4\gamma_{1}(1-\gamma_{1})>\gamma_{1}^{2}(1-\gamma_{1})>(1-\gamma_{1})(\gamma_{1}^{2}-\gamma_{2}^{2}).
 \]
 
Clearly  $$4\gamma_{1}(1-\gamma_{1})\geq 2(\gamma_{1}^{2}-\gamma_{2}^{2}) \Leftrightarrow 3\gamma_{1}^{2}-2\gamma_{1}-\gamma_{2}^{2}\leq 0$$ is true because  $0<\gamma_{1}\leq \frac{1+\sqrt{1+3\gamma_{2}^{2}}}{3}$.
Thus we can create a horizontal-like boundary for $\widetilde{V}$ using the collection of fibers of the foliation of $U^{-}$ inside $V$ which pass through $|y|=\frac{1+\gamma_{1}}{2}R$, $x=0$.  Note also that these fibers are long, as they will only exit $V$ through the vertical boundary $|x|=R$.

To complete the argument, each connected component of a vertical-like leaf which intersects the circle $|x|=\frac{1+\gamma_{1}}{2}R$, $y=0$ has the form $\{(\varphi(z),z), |z|\leq R\}$. Each connected component of a horizontal-like leaf which intersects the circle $|y|=\frac{1+\gamma_{1}}{2}R$, $x=0$ has the form $\{(z,\psi(z)), |z|\leq R\}$. Since by the arguments above there exists $R'$ such that $|\varphi(z)|<R'<R$ and $|\psi(z)|<R'<R$ for all $z\in\D_{R}$, it follows by Rouch\'e's Theorem that each such pair of leaves intersect in exactly one point.

The vertical-like and horizontal-like boundaries described above determine a closed bounded set $\widetilde{V}$ with the same dynamical properties as $V$.
\qed

 \begin{remark} We point out that $\gamma_{2}< \frac{1+\sqrt{1+3\gamma_{2}^{2}}}{3}$, so the constraints on $\gamma_{1}$ and $\gamma_{2}$ from Proposition \ref{prop:dynV} are not vacuous. 
 \end{remark}

 \begin{prop}[\textbf{Constructing the inner tubes}]\label{prop:dyntubes} Suppose that $\gamma_{2}\geq \frac{1}{\sqrt{2}}\gamma_{1}$.
 \begin{enumerate}
 \item There exists a set $B^{v}$ such that
 \begin{enumerate}
 \item[a)] $\{(x,y)\in\widetilde{V},\ |x|\leq 1+|a|\}\subset B^{v}\subset \{(x,y)\in\widetilde{V},\ |x|\leq \gamma_{1}R\}$;
 \item[b)] the vertical boundary of $B^{v}$ consists of connected components of leaves of the foliation of $U^{+}$ inside $\widetilde{V}$, which are vertical-like;
 \item[c)] the horizontal boundary of $B^{v}$ is a subset of the horizontal boundary of $\widetilde{V}$ and consists of connected components of leaves of the foliation of $U^{-}$, which are horizontal-like by the construction of the set $\widetilde{V}$;
 \item[d)] $B^{v}\cap H^{-1}(\widetilde{V})=\emptyset$.
 \end{enumerate}
 \item Similarly, there exists a set $B^{h}$ such that
 \begin{enumerate}
 \item[a)] $\{(x,y)\in\widetilde{V},\ |y|\leq 1+|a|\}\subset B^{h}\subset \{(x,y)\in\widetilde{V},\ |y|\leq \gamma_{1}R\}$;
 \item[b)] the horizontal boundary of $B^{h}$ consists of leaves of the foliation of $U^{-}$ inside $\widetilde{V}$, which are horizontal-like;
 \item[c)] the vertical boundary of $B^{h}$ is a subset of the vertical boundary of $\widetilde{V}$ and consists of leaves of the foliation of $U^{+}$ inside $\widetilde{V}$, which are vertical-like by the construction of the set $\widetilde{V}$;
 \item[d)] $B^{h}\cap H(\widetilde{V})=\emptyset$.
 \end{enumerate}
 \end{enumerate}
 
 See Figure \ref{fig:VandTube} for a schematic drawing of the tube $B^{v}$.
 \end{prop}
 \proof
 Recall that the polydisk $\{1+|a|\leq|x|\leq\gamma_{2}R,\ |y|\leq R\}$ is contained $V\setminus H^{-1}(V)$, and consequently in $U^{+}$, and is disjoint from the union of all the backward iterates of $B^{v}_{\delta}$. Therefore the foliation of $U^{+}$ in this region consists of vertical-like leaves.

\begin{figure}[htb]
\begin{center}
\includegraphics[scale=1]{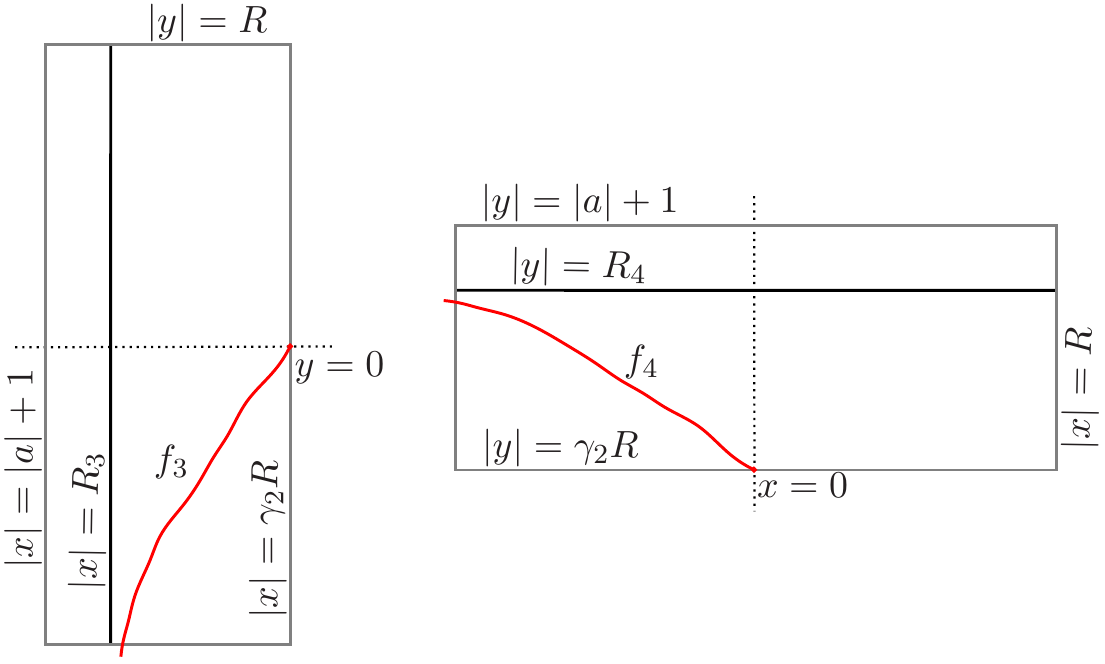}
\end{center}
\caption{Choosing $R_{3}$ and $R_{4}$.}
\label{fig:fibers3C}
\end{figure}

 There exists $R_{3}$ such that $1+|a|\leq R_{3}<\gamma_{2}R$ and the vertical-like leaves of the foliation of $U^{+}$ inside $V$ which pass through the circle $|x|=\gamma_{2}R$, $y=0$ do not cross $|x|=R_{3}$ (see Figure \ref{fig:fibers3C}). 
 
 If
 \begin{equation}\label{eq:cond1}
 \frac{R|a|}{2R_{3}-1}<\gamma_{2}R-R_{3},
 \end{equation}
 then this follows from Corollary \ref{cor:proj}. 
 Condition \eqref{eq:cond1} is equivalent to
 \[
 2R_{3}^{2}-R_{3}(2\gamma_{2}R+1)+R(|a|+\gamma_{2})<0.
 \]
 We first need to choose $R$ so that the discriminant $\Delta_{1}=(2\gamma_{2}R-1)^{2}-8R|a|\geq 0$, and then require that
 \[
 \frac{2\gamma_{2}R+1-\sqrt{\Delta_{1}}}{4}<R_{3}<\frac{2\gamma_{2}R+1+\sqrt{\Delta_{1}}}{4}.
 \]
 Notice that Equation \eqref{eq:Rlower} gives
 \begin{equation}\label{eq:aux5}
 \frac{2\gamma_{2}R+1+\sqrt{\Delta_{1}}}{4}>\frac{\gamma_{2}R}{2}> \frac{\gamma_{2}}{\gamma_{1}^{2}-\gamma_{2}^{2}}(1+|a|),
 \end{equation}
 which is greater than $1+|a|$ provided that 
 $\gamma_{2}^{2}+\gamma_{2}-\gamma_{1}^{2}> 0$. This last inequality is equivalent to $\gamma_{2}>\frac{-1+\sqrt{1+4\gamma_{1}^{2}}}{2}$. By assumption, we have $\gamma_{2}\geq \frac{1}{\sqrt{2}}\gamma_{1}$, and an easy computation shows that $\frac{1}{\sqrt{2}}\gamma_{1}>\frac{-1+\sqrt{1+4\gamma_{1}^{2}}}{2}$ for $0<\gamma_{1}<1$. Thus it is possible to find such $R_{3}\geq 1+|a|$.
 
 Mirroring the arguments above for $U^{-}$, we can show that there exists $R_{4}$ such that $1+|a|\leq R_{4}<\gamma_{2}R$ and the horizontal-like leaves of the foliation of $U^{-}$ inside $V$ which pass through the circle $|y|=\gamma_{2}R$,  $x=0$ do not cross $|y|=R_{4}$ (see Figure \ref{fig:fibers3C}).
 If
 \begin{equation}\label{eq:cond2}
 \frac{R}{2R_{4}-|a|}<\gamma_{2}R-R_{4},
 \end{equation}
 then this follows by Corollary \ref{cor:proj}. Condition \eqref{eq:cond2} is equivalent to
 \[
 2R_{4}^{2}-R_{4}(2\gamma_{2}R+|a|)+R(1+\gamma_{2}|a|)<0.
 \]
 As before, we first need to choose $R$ so that the discriminant $\Delta_{2}=(2\gamma_{2}R-|a|)^{2}-8R\geq 0$, and then ask that
 \[
 \frac{2\gamma_{2}R+|a|-\sqrt{\Delta_{2}}}{4}<R_{4}<\frac{2\gamma_{2}R+|a|+\sqrt{\Delta_{2}}}{4}
 \]
As in \eqref{eq:aux5} we have
 \begin{equation}\label{eq:aux6}
 \frac{2\gamma_{2}R+|a|+\sqrt{\Delta_{1}}}{4}>\frac{\gamma_{2}R}{2},
 \end{equation}
 which was shown above to be greater than $1+|a|$ when $\gamma_{2}\geq \frac{1}{\sqrt{2}}\gamma_{1}$. Note that \eqref{eq:aux5} and \eqref{eq:aux6} are rather rough inequalities, so other conditions on $\gamma_{1}$ and $\gamma_{2}$ might also work, as we will discuss later.
Nonetheless, we have shown that it is possible to find $R_{4}\geq 1+|a|$.

We now address the sign of the discriminants $\Delta_{1}$ and $\Delta_{2}$:
\[
\Delta_{1} \geq 0 \Leftrightarrow 4\gamma_{2}^{2}R^{2}-(4\gamma_{2}+8|a|)R+1\geq 0
\]
and
\[
\Delta_{2} \geq 0 \Leftrightarrow 4\gamma_{2}^{2}R^{2}-(4\gamma_{2}|a|+8)R+|a|^{2}\geq 0.
\]
From the two quadratic equations we get two lower bounds:
\[
R \geq l_{1}:=\frac{\gamma_{2}+2|a|+\sqrt{(\gamma_{2}+2|a|)^{2}-\gamma_{2}^{2}}}{2\gamma_{2}^{2}}
\]
and
\[
R \geq l_{2}:=\frac{\gamma_{2}|a|+2+\sqrt{(\gamma_{2}|a|+2)^{2}-|a|^{2}\gamma_{2}^{2}}}{2\gamma_{2}^{2}}.
\]
If $\ds \gamma_{2}\geq \frac{1}{\sqrt{2}}\gamma_{1}$ then $\ds \frac{\gamma_{2}}{\gamma_{1}^{2}-\gamma_{2}^{2}}\geq \frac{2}{\gamma_{2}^{2}}$. By Equation \eqref{eq:Rlower} we know that
 \begin{equation}\label{aux7}
R\geq\frac{2}{\gamma_{1}^{2}-\gamma_{2}^{2}}(1+|a|)\geq  \frac{2}{\gamma_{2}^{2}}(1+|a|).
 \end{equation}
We can now use \eqref{aux7} in two ways:
\begin{eqnarray*}
\frac{2}{\gamma_{2}^{2}}(1+|a|) &>&\frac{\gamma_{2}+2|a|}{\gamma_{2}^{2}}>l_{1}\\
\frac{2}{\gamma_{2}^{2}}(1+|a|) &>& \frac{\gamma_{2}|a|+2}{\gamma_{2}^{2}}>l_{2}\
\end{eqnarray*}
It follows that both discriminants are non-negative.

In order to finish the construction of $B^{v}$, it suffices to construct its vertical boundary $B^{v}$. For this we take the connected components of the vertical-like leaves of the foliation of $U^{+}$ in $\widetilde{V}$ which pass through the circle $|x|=\gamma_{2}R$, $y=0$. We know that they do not cross $|x|=R_{3}$ in $V$ by the discussion above. They cannot intersect $H^{-1}(\widetilde{V})$ either, because the vertical boundaries of $H^{-1}(\widetilde{V})\cap\widetilde{V}$ consist of vertical-like leaves of the same foliation of $U^{+}$ by the construction of the dynamical set $\widetilde{V}$. So each leaf from the vertical boundary of $B^{v}$ intersects the horizontal boundary of $\widetilde{V}$. The construction for $B^{h}$ is similar, we leave the details to the reader.
\qed

\noindent \textbf{$\HOV_{\beta}$ conditions:} So far we have imposed the following conditions:
 \begin{equation*}
|c|>\frac{2(\gamma_{1}^{2}+\gamma_{2}^{2})}{(\gamma_{1}^{2}-\gamma_{2}^{2})^{2}}(1+|a|)^{2},
\end{equation*}
 where $0<\gamma_{2}<\gamma_{1}<1$, $\ds \gamma_{1}\leq \frac{1+\sqrt{1+3\gamma_{2}^{2}}}{3}$ and $\ds \gamma_{1}\leq \sqrt{2}\gamma_{2}$.
 
 \begin{remark} The minimum of the 2-parameter function $(\gamma_{1},\gamma_{2})\mapsto\frac{2(\gamma_{1}^{2}+\gamma_{2}^{2})}{(\gamma_{1}^{2}-\gamma_{2}^{2})^{2}}$ taken after all $\gamma_{1}$ and $\gamma_{2}$ subject to the restrictions above is $\beta=18.75$, obtained when $\gamma_{2}=\frac{2\sqrt{2}}{5}$ and $\gamma_{1}=\frac{4}{5}$. In this case, 
$\frac{1+\sqrt{1+3\gamma_{2}^{2}}}{3}= \sqrt{2}\gamma_{2}$.

On the boundary of this HOV$_{\beta}$ region, that is when $|c|=\beta(1+|a|)^{2}$, $\beta=18.75$, we have $\eta_{1}=\gamma_{1}R$, $\eta_{2}=\gamma_{2}R$, and the lower and upper bounds for $R$ from \eqref{eq:Rlowup} coincide, giving only one choice for $R=6.25(1+|a|)$. This shows that our strategy for constructing a dynamical polydisk $\widetilde{V}$ and dynamical tubes $B^{v}$, $B^{h}$ which will later serve as trapping regions for the critical locus is somehow strong.
\end{remark}

\begin{remark} A better HOV$_{\beta}$ region with $\beta<18.75$ such that Propositions \ref{prop:technical} , \ref{lemma:proj}, \ref{prop:dynV}, \ref{prop:dyntubes} hold is possible, and the value $18.75$ is not the most optimal one. By analyzing carefully the proof of Proposition \ref{prop:dyntubes} and using better inequalities we can, for example, get HOV$_{\beta'}$ with $\beta'\simeq 17.82$, but of course even lower bounds may be possible. If we also put bounds on the size of $|a|$ we can produce larger HOV-like regions up to a certain Jacobian.
\end{remark}

\begin{remark} If we just want to build a dynamical polydisk $\widetilde{V}$, and not the inner tubes $B^{v}$ and $B^{h}$ with dynamical boundaries (keep $B^{h}_{\delta}$ and $B^{v}_{\delta}$ or replace them with different sets using different arguments, other than those from Proposition \ref{prop:dyntubes}), then we just take $\gamma_{2}=0$ in Propositions \ref{prop:technical} and \ref{prop:dynV}. We then find $|c|>\frac{2}{\gamma_{1}^{2}}(1+|a|)^{2}$ and $0<\gamma_{1}\leq \frac{2}{3}$, which is optimized for $\gamma_{1}=\frac{2}{3}$ and $\beta=\frac{9}{2}$. So we get the HOV$_{\beta}$ region defined by $|c|>\frac{9}{2}(1+|a|)^{2}$.
\end{remark}

\noindent\textbf{Dynamical properties of the set $\widetilde{V}$:}
\begin{enumerate}

\item  The Julia sets can be obtained by taking nested intersections:
\begin{eqnarray*}
J^{+} \cap \widetilde{V} = \bigcap_{n\geq 0} H^{-n}(\widetilde{V}), \quad
 J^{-}\cap \widetilde{V} = \bigcap_{n\geq 0} H^{n}(\widetilde{V}), \quad
J = \bigcap_{n\in \Z} H^{n}(\widetilde{V}). 
\end{eqnarray*}

\item If $z\in \widetilde{V}\setminus H^{-1}(\widetilde{V})$, then $z\in U^{+}$, and if we let $(x_{n},y_{n})=H^{n}(z)$, $n\geq 0$, then $|x_{1}|>\gamma_{1}R$ and $|x_{n}|>|x_{n-1}|$ for every $n>1$.
\item If $z\in \widetilde{V}\setminus H(\widetilde{V})$, then $z\in U^{-}$, and if we let $(x_{n},y_{n})=H^{-n}(z)$, $n\geq 0$, then $|y_{1}|>\gamma_{1}R$ and $|y_{n}|>|y_{n-1}|$ for every $n>1$.

\smallskip
\proof Recall that $V_{1}\subset \widetilde{V} \subset V$, where $V_{1}$ is the polydisk of radius $\gamma_{1} R$. Note that both $R$ and $\gamma_{1}R$ satisfy inequality \eqref{eq:R}, therefore both $V_{1}$ and $V$ are eligible polydisks for the dynamical filtration of $\C^{2}$ (see Figure \ref{fig:filtration}).
 The proofs for parts (2) and (3) follow elementary from the standard properties of the filtration. 
 
In addition, by construction, the vertical boundary of $\widetilde{V}$ is in $V\setminus H^{-1}(V)$, and the  horizontal boundary of $\widetilde{V}$ is included in $V\setminus H(V)$. Therefore 
\begin{eqnarray*}
\widetilde{V} \cap H^{-2}(V) \quad \subseteq  & \widetilde{V}\cap H^{-1}(\widetilde{V}) & \subseteq \quad \widetilde{V}\cap H^{-1}(V)\\
\widetilde{V} \cap H^{2}(V) \quad \subseteq & \widetilde{V}\cap H(\widetilde{V}) &\subseteq \quad \widetilde{V} \cap H(V).
\end{eqnarray*}
Hence $\bigcap_{n\geq 0} H^{-n}(\widetilde{V}) = \widetilde{V}\cap \bigcap_{n>0} H^{-n}(V)$, $\bigcap_{n\geq 0} H^{n}(\widetilde{V}) = \widetilde{V}\cap \bigcap_{n>0} H^{n}(V)$, and 
 $\bigcap_{n\in \Z} H^{n}(V) = \bigcap_{n\in \Z} H^{n}(\widetilde{V})$. The conclusion of part (1) then follows from the horseshoe construction in Equation \eqref{eq: horseshoe}.
\qed
\end{enumerate}

\noindent\textbf{Binary coding of the tubes.}  In this subsection, we will make use of the horseshoe labeling process, to explain how one can apply it to do a binary coding of the forward images of the tube $B^{h}$ and the backward images of the tube $B^{v}$ from Proposition \ref{prop:dyntubes}. We assume that the horseshoe construction is known (see e.g. \cite{MMKT}, \cite{MNTU}).

We first label the two connected components of the set $H(\widetilde{V})\cap \widetilde{V}$ as $T^{h}_{0}$ and  $T^{h}_{1}$.  This induces a matching labelling on the two connected components of the set $H^{-1}(\widetilde{V})\cap \widetilde{V}$, since $H^{-1}\left(H(\widetilde{V})\cap \widetilde{V}\right)=H^{-1}(\widetilde{V})\cap \widetilde{V}$. Hence we denote 
\[
T^{v}_{i}:=H^{-1}(T^{h}_{i}),\ \  \mbox{for}\ i=0,1.
\]

Following the standard horseshoe construction, when $k\geq 2$, we can label inductively the $2^{k }$ connected components of $H^{k}(\widetilde{V})\cap \widetilde{V}$ as $T^{h}_{w}$, using finite strings $w=w_{0}\ldots w_{k-1}$ with $k$ letters from the alphabet $\{0,1\}$, by applying the rule
\[T^{h}_{w}:=T^{h}_{w_{0}}\cap H(T^{h}_{w_{1}})\cap \ldots \cap H^{k-1}(T^{h}_{w_{k-1}}).\]
Equivalently, if we denote by $|w|$ the length of $w$, then the rule can be stated as $w_{j}=i$ iff $H^{-j}(T^{h}_{w})\subset T^{h}_{i}$ where $j=1,\ldots,|w|-1$ and $i\in \{0,1\}$.

We now follow the standard horseshoe coding process for the vertical tubes. When $k\geq 2$, we can label inductively the $2^{k}$ connected components of $H^{-k}(\widetilde{V})\cap \widetilde{V}$ as $T^{v}_{w}$, using finite strings $w=w_{-k}\ldots w_{-1}\in\{0,1\}^{k}$ by applying the rule
\[T^{v}_{w}:=T^{v}_{w_{-1}}\cap H^{-1}(T^{v}_{w_{-2}})\cap \ldots \cap H^{-k}(T^{v}_{w_{-k}}).\]
Equivalently, if we denote by $|w|$ the length of $w$, then the rule can be stated as $w_{j}=i$ iff the forward iterates satisfy $H^{|j|}(T^{v}_{w})\subset T^{v}_{i}$ where $j=-2,\ldots,-|w|$ and $i\in \{0,1\}$.

Note that $H^{|w|}(T^{v}_{w})=T^{h}_{\alpha}$, where $|w|=k\geq 1$, $w=w_{-k},\ldots, w_{-1}$, $\alpha=\alpha_{0},\ldots, \alpha_{k-1}$, and $\alpha_{0}=w_{-k},\ldots, \alpha_{k-1}=w_{-1}$.

The Julia set $J$ of any H\'enon map from the HOV region is a Cantor set, and the map
\begin{equation}\label{def:conjugacy}
\tau:\Sigma_{2}\rightarrow J,\ \  \tau(w)=\bigcap_{k=-\infty}^{\infty}H^{k}(T^{h}_{w_{k}})
\end{equation}
is a homeomorphism conjugating $\sigma$, the shift map on two symbols, to the H\'enon map restricted to its Julia set $J$.

The tubes $B^{h}$ and $B^{v}$ do not intersect the forward, respectively backward iterate of $\widetilde{V}$, so they do not have any induced binary labeling from the horseshoe construction. However, all forward iterates of $B^{h}$ naturally  inherit the labeling of the forward iterates of $\widetilde{V}$; analogously, all backward iterates of $B^{v}$ inherit the labeling of the backward iterates of $\widetilde{V}$. 

Note that 
for any $k\geq 1$, the set $H^{-k}(B^{v})\cap \widetilde{V}$ has $2^{k}$ connected components, which lie in the $k$-th, but not in the $(k+1)$-th preimage of $\widetilde{V}$ (see Figure \ref{fig:horseshoe} as a guideline).
Therefore there exists a unique coding assignment such that $B^{v}_{w}\subset T^{v}_{w}$, where $|w|=k$. 
Likewise, the forward iterates of $B^{h}$ inherit the coding of the forward iterates of $\widetilde{V}$, hence for $k\geq 1$, each of the $2^{k}$ connected components of $H^{k}(B^{h})\cap \widetilde{V}$ has a unique labeling $B^{h}_{\omega}$, where $\omega$ is the unique binary string on length $k$ such that $B^{h}_{\omega}\subset T^{h}_{\omega}$.

We denote the intersections  by 
\begin{equation}\label{eq: poly}\Delta_{\lambda, \omega}:=B^{v}_{\lambda}\cap B^{h}_{\omega}.
\end{equation} 
For the intersections between the labeled tubes $B^{v}_{\lambda}$ and  $B^{h}_{w}$ with the unlabeled tubes $B^{h}$ and $B^{v}$, or between the unlabeled tubes $B^{h}$ and $B^{v}$, we will use the complimentary notations:
\begin{eqnarray}\label{eq:notations}
X_{\lambda}:=\Delta_{(\lambda, \emptyset)}=B^{v}_{\lambda}\cap B^{h},\quad\quad
Y_{\omega}:=\Delta_{( \emptyset, \omega)}=B^{v}\cap B^{h}_{\omega} 
\end{eqnarray}
and 
\begin{equation}\label{eq:B}
B:=\Delta_{( \emptyset, \emptyset)}=X_{\emptyset}=Y_{\emptyset}=B^{v}\cap B^{h}.
\end{equation}
In order to make the notations more intuitive, we use the upper scripts h/v and the notations X/Y for sets within the horizontal/vertical tubes.

\section{The critical locus and trapping regions}\label{sec:CriticalLocus}

We proceed with the description of the foliations of $U^{+}$ and $U^{-}$ and their tangencies.
It is easy to see from Proposition \ref{prop:foliation1} that on the set 
\[(U^+\cap U^-) \setminus \left (\bigcup_{n\geq 0} H^{-n}(S_{\delta}^v) \cup \bigcup_{n\geq 0} H^{n}(S_{\delta}^h)\right),\] 
the tangent spaces to the leaves of foliation of  $U^+$, respectively $U^-$,  belong to vertical, and respectively horizontal cones; therefore the foliations of $U^+$ and $U^{-}$ are everywhere transverse to each other, and the critical locus in this region is the empty set.

In fact, since $H(S^{v}_{\delta})=S^{h}_{\delta}$, we see that the critical locus $\mathcal{C}$ is contained in the union of forward and backward images of the strip $S_{\delta}^h$, a fact which we formulate as Lemma \ref{cor:c} below, to use for future reference:

\begin{lemma}\label{cor:c} For each point $z\in \Crit$, there exists $n\in {\mathbb Z}$ such that $H^n(z)\in{S_{\delta}^h}$.
\end{lemma}

\begin{thm}[\textbf{Fundamental domain}]\label{thm:FundamentalDomain} Let $F:=\left(\widetilde{V}\cup S_{\delta}^h\right)\backslash H(\widetilde{V})$. Then $\Crit \cap F$ is a fundamental domain of the critical locus.
\end{thm}
\proof We proceed in two steps.
\begin{enumerate}
\item We first show that for each $z\in \Crit$, there exists $n\in {\mathbb Z}$ such that $H^n(z)\in F$.
 By Corollary \ref{cor:c}, for each $z\in \Crit$ there exists $n\in {\mathbb Z}$ such that $\omega=H^{n}(z)\in S_{\delta}^h$. If $\omega\not \in H(\widetilde{V})$, then $\omega\in F$, and we are done. If $\omega\in H(\widetilde{V})$, then $H^{-1}(\omega)\in \widetilde{V}$. Since $\omega$ belongs to the escaping set $U^{-}$, there exists a smallest positive integer $m\geq 1$ such that $H^{-k}(\omega) \in \widetilde{V}$, for $1\leq k<m$, and $H^{-m}(\omega) \notin \widetilde{V}$. This implies that $H^{-(m-1)}(\omega)\in \widetilde{V}\backslash H(\widetilde{V})$. Hence $H^{-(m-1)}(\omega)\in F$.

\item Given $z\in \Crit\cap F$, we show that there cannot exist another iterate $n\neq0$ such that $\omega=H^n(z)\in F$. By changing the roles of $z$ and $\omega$, it suffices to prove this for $n>0$. 
\begin{enumerate} 
\item If $z\in S_{\delta}^h\backslash \widetilde{V}$, then let $z=(x,y)$ and $H^n(z)=(x_n,y_n)$ for $n>0$. By the properties of the Hubbard filtration $V_{\gamma_{1}R}$ we know that 
\[y_1=x, \mbox{ and } |y_n|>|y_1|>\gamma_{1}R>\delta \mbox{ for all } n>1.\] 
Hence $w\not \in S_{\delta}^h$, so it is not in the set $F$, which is a contradiction.

\item Assume now that $z\in \widetilde{V}$. Since $z$ belongs to the escaping set $U^{+}$, there exists a smallest integer $m>0$ such that $H^k(z)\in \widetilde{V}$, for $0\leq k<m$ and $H^{m}(z)\notin\widetilde{V}$. As in part (a), this implies that for $k>m$, we have $H^k(z)\notin S_{\delta}^h$, hence $H^k(z)\not\in F$. For $0\leq k<m$, we have $H^k(z) \in \widetilde{V}$. 
This implies $H^k(z)\in H(\widetilde{V})$ for $0<k\leq m$, so these iterates cannot belong to the set $F$ either. Contradiction.  \qed
\end{enumerate}
\end{enumerate}

\begin{cor}\label{cor:fund} $F'\cap\mathcal{C}$ is another fundamental domain of the critical locus, where $F':=\left(\widetilde{V}\cup S_{\delta}^v\right)\backslash H^{-1}(\widetilde{V})$.
\end{cor}

\begin{remark} Note that in the Lyubich-Robertson case, $S_{\delta}^h\cap \Crit$ is a fundamental domain of the critical locus. On the other hand, in the HOV$_{\beta}$-region, $S_{\delta}^{h}\cap \Crit$ is not a fundamental domain, since there are points of the critical locus that belong to $S_{\delta}^h\cap H^n(S_{\delta}^h)$ for $n\neq 0$.
\end{remark}

\begin{thm} There exists $R'>0$ large enough such that the critical locus in the region $S_{\delta}^h \cap \{|x|>R'\}$ is a punctured disk with a hole at $x=\infty$, asymptotic to the $x$-axis. Similarly, the critical locus in the region $S_{\delta}^v \cap \{|y|>R'\}$ is a punctured disk with a hole at $x=\infty$, asymptotic to the $y$-axis.

\end{thm}
\proof 
The rigorous description of the critical locus at $\infty$ has been given in \cite{LR} and \cite{F}. The value of $R'$ such that $V\cup H(V)\subset \{|x|<R', |y|<R'\}$ works. We briefly explain why the critical locus in the strip $S^{v}_{\delta}$ is tangent at $\infty$ to the $y$-axis. The Botcher coordinates satisfy $\varphi^{-}(x,y)\simeq y$ when $(x,y)\rightarrow \infty$ in $V^{-}$, whereas $\varphi^{+}(x,y)\simeq x$ when $(x,y)\rightarrow \infty$ in $V^{+}$. Assume for simplicity, that for $|z_{0}|$ large enough, the leaf of the foliation of $U^{-}$ given by $(\varphi^{-})^{-1}(z_{0})\cap S^{v}_{\delta}$ is the horizontal disk $(z, z_{0})$, $|z|\leq \delta$, whereas the leaf of the foliation of $U^{+}$ given by $(\varphi^{+})^{-1}(z_{0})\cap S^{h}_{\delta}$ is the vertical disk $(z_{0}, z)$, $|z|\leq \delta$, whose preimage $H^{-1}((\varphi^{+})^{-1}(z_{0}))\cap S^{v}_{\delta}$ is therefore a vertical parabola parametrized by $(z, (z^{2}-z_{0})/a)$,  $|z|\leq \delta$. The two leaves develop a tangency inside $S^{v}_{\delta}$ exactly when $z=0$, which explains how the critical locus inside $S^{v}_{\delta}$ is tangent to the $y$-axis at $\infty$.
\qed

We will first study the critical locus inside $S^{v}_{\delta}$ and $S^{h}_{\delta}$, so we need to describe the foliations of $U^{+}$ and $U^{-}$ in a large part of these strips.

\begin{defn}\label{def:paraleaf} We say that a leaf of a foliation is {\it vertical parabolic-like} if its projection on the second coordinate is two-to-one except at one point. Similarly, we say that the leaf is {\it horizontal parabolic-like} if its projection on the first coordinate is two-to-one except at one point. 
\end{defn}

\begin{prop}[\textbf{Parabolic-like leaves}]\label{prop:foliation2} The foliation of $U^+$ in the set 
\[
S_{\delta}^v \cup B^{v}\setminus H^{-1}(\widetilde{V})
\] is (vertical) parabolic-like. 
The foliation of $U^-$ inside the region $S_{\delta}^h \cup B^{h} \setminus H(\widetilde{V})$ is (horizontal) parabolic-like.
\end{prop}
\proof
We know that $H(S_{\delta}^v) \subset S_{\delta}^h$, $H(B^{v})\subset \{|x|>\gamma_{1}R, |y|<\gamma_{1}R\}$, $H(B^{v})\cap \widetilde{V}=\emptyset$, and that the foliation of $U^{+}$ in the set $W:=\{|y|\leq \gamma_{1}R\}\setminus \widetilde{V}$ is vertical-like. Let us consider a part of a vertical-like fiber of the foliation of $U^{+}$ in $W$, parametrized by $(\varphi(z),z), |z|\leq \gamma_{1}R$. Its pull-back under $H$ is given by $L(z)=(z, (p(z)-\varphi(z))/a)$. The derivative of the second coordinate is equal to $0$ iff $2z=\varphi'(z)$. 

By Lemma \ref{prop:cones} we have the following estimate:
\[
|\varphi'(z)|<\frac{|a|}{2|\varphi(z)|-1}<\frac{|a|}{2\gamma_{1}R-1}
\]
In any case, $\ds R>\frac{2}{\gamma_{1}^{2}-\gamma_{2}^{2}}(1+|a|)>\frac{2}{\gamma_{1}^{2}}(1+|a|)$, so
\begin{equation}\label{eq:gamma1R}
\frac{|a|}{2\gamma_{1}R-1}<\frac{|a|\gamma_{1}}{4|a|+4-\gamma_{1}}<\frac{|a|}{4|a|+3}<1.
\end{equation}

However, an easy application of Rouch\'e's Theorem shows that this happens for exactly one point $|z_{0}|<\delta$, since on the boundary $|z|=\delta$ we have 
\[\ds |2z|=|a|+1\geq 1>\frac{|a|}{2\gamma_{1}R-1}> |\varphi'(z)|,\] therefore the functions $2z-\varphi'(z)$ and $2z$ have the same number of zeros inside the disk $|z|<\delta$.  We use Rouch\'e's Theorem once more to show that the same equation $2z=\varphi'(z)$ has no solutions in the annulus $\delta\leq |z|\leq \gamma_{1}R$. When $z\neq z_{0}$, the degree of the map $(p(z)-\varphi(z))/a$ is equal to two, so $L$ projects two-to-one over the $y$-axis. 

Using Rouch\'e's Theorem and the estimates for the horizontal cone from Lemma \ref{prop:cones}, we get that the foliation of $U^-$ inside $S_{\delta}^h\cup B^{h} \setminus H(\widetilde{V})$ is (horizontal) parabolic-like, and the ``tips of the parabolas'' lie in the strip $S_{\delta}^h$.
\qed

The following lemma will be of use in Lemmas \ref{lemma:NOsing} and \ref{lemma:NOsingN}, where we analyze the critical locus inside $B^{h}$ and $B^{v}$.

\begin{lemma}\label{lemma:NOhorver} The leaves of the foliation of $U^{-}$ in the horizontal tube $B^{h}$ have no horizontal tangent lines. The leaves of the foliation of $U^{+}$ in the vertical tube $B^{v}$  have no vertical tangent lines. 
\end{lemma}
\proof Suppose that there exists a connected component $\mathcal{L}$ of a leaf of the foliation of $U^{+}$ in $B^{v}$  with a vertical tangent line at some point $(x_{0},y_{0})$. Let $(x_{n},y_{n})=H^{n}(x_{0},y_{0})$ for $n\geq 1$. The derivative of the H\'enon map sends a vertical vector $v_{0}=(0,\eta)$ to a horizontal vector $v_{1}=(-a\eta,0)$. Then $H(\mathcal{L})$ has a horizontal tangent line at the point $(x_{1},y_{1})$. Notice that the tube $B^{v}$ is mapped under $H$ outside of the set $\widetilde{V}$, where by Proposition \ref{prop:cones} we have an invariant family of horizontal cones. Hence for every $n\geq 1$, the tangent line to the leaf $H^{n}(\mathcal{L})$ at the point $(x_{n},y_{n})$ belongs to the horizontal cone at $(x_{n},y_{n})$, However, $(x_{n},y_{n})\rightarrow \infty$ in $V^{+}$, and $\varphi^{+}(x,y)\simeq x$ as $(x,y)\rightarrow \infty$ in $V^{+}$, hence the tangent line to $H^{n}(\mathcal{L})$ at the point $(x_{n},y_{n})$ must belong to the vertical cone at $(x_{n},y_{n})$ as $n\rightarrow \infty$. Contradiction.

The second part of the lemma is proved identically, making use of the fact that the derivative of the inverse H\'enon map sends a horizontal vector to a vertical vector, which thereafter remains inside the vertical cones under all backward iterations of the H\'enon map, thus contradicting the fact that $\varphi^{-}(x,y)\simeq y$ as $(x,y)\rightarrow \infty $ in $V^{-}$.
\qed

\begin{prop}\label{prop:longleaves} The foliation of $U^{-}$ in the set 
\[
	\widetilde{V}_{1}:=\widetilde{V} \setminus \left(\bigcup\limits_{n\geq  0} H^{n}(B^h)\cup J^{-}\right)
\] consists of long horizontal-like holomorphic disks which can be parametrized by $(z,\psi(z))$, with $|z|\leq M$ , $M\geq \gamma_{1}R$ and contraction factor $ |\psi'(z)|<\frac{1}{2R_{4}-|a|}<\frac{1}{|a|+2}$.

The foliation of $U^{+}$ in the set \[
	\widetilde{V}_{2}:=\widetilde{V} \setminus \left(\bigcup\limits_{n\geq  0} H^{-n}(B^v)\cup J^{+}\right)
\] consists of long vertical-like holomorphic disks which can be parametrized by $(\varphi(z),z)$, with $|z|\leq M,\ M\geq \gamma_{1}R$ and contraction factor $ |\varphi'(z)|<\frac{|a|}{2R_{3}-1}<\frac{|a|}{2|a|+1}$.

The set $\widetilde{V}_{2}$ is backward invariant in the sense that $H^{-1}(\widetilde{V}_{2})\cap \widetilde{V}\subset \widetilde{V}_{2}$. Likewise, the set $\widetilde{V}_{1}$ is forward invariant in the sense that $H(\widetilde{V}_{1})\cap \widetilde{V}\subset \widetilde{V}_{1}$.
\end{prop}

\proof The proof is done by induction, making use of the horseshoe structure, of the dynamical construction of the sets $\widetilde{V}$, $B^{h}$ and $B^{v} $ from Propositions \ref{prop:dynV}, \ref{prop:dyntubes}, and of the estimates in the families of invariant cones from Proposition \ref{prop:cones}. It suffices to do the proof  for example for $U^{-}$, as it can easily be adapted step-by-step for $U^{+}$.

By Proposition \ref{prop:dyntubes}, the tube $B^{h}$ is a subset of $\widetilde{V}\setminus H(\tilde{V})$, and the horizontal boundary of $B^{h}$ does not intersect the inner tube $|y|\leq R_{4}$, where $R_{4}>1+|a|$. All the forward images $H^{n}(B^{h}_{\delta})\cap \widetilde{V}$, $n>0$, where we loose our invariant family of horizontal cones, are inside $H(\widetilde{V})$. On the set $W_{1}:=\widetilde{V}\setminus (B^{h}\cup H(\widetilde{V}))$ the foliation of $U^{-}$ is horizontal-like, with leaves of the form $(z,\psi(z))$ and 
\[
|\psi'(z)|<\frac{1}{2R_{4}-|a|}<\frac{1}{2(1+|a|)-|a|}=\frac{1}{2+|a|}.
\]
In addition, any leaf of the foliation of $U^{-}$ inside $W_{1}$ cannot exit $W_{1}$ through the horizontal boundary, as this consists of the horizontal boundaries of the sets $\widetilde{V}, H(\widetilde{V})$, and $B^{h}$, all laminated by leaves of the same foliation of $U^{-}$. Hence, such a leaf will exit $W_{1}$ through its vertical boundary, which is just a subset of the vertical boundary of $\widetilde{V}$. The set $\widetilde{V}$ contains the polydisk $V_{\gamma_{1}R}$ inside, so the parameter $z$ belongs to a disk of radius larger than $\gamma_{1}R$.

By induction on $n\geq1$, we can repeat the argument for each of the sets 
\[W_{n+1}:=H(W_{n})\cap \widetilde{V}=H^{n}(W_{1})\cap \widetilde{V}.
\]
The contraction factor is the same (or better), because each set $W_{n}$,$n>1$ is disjoint from the tube $|y|\leq R_{4}$. We can conclude the proof by passing to the limit, and noting that in the horseshoe region we have $\bigcap_{n\geq 0} H^{n}(\widetilde{V}) \cap \widetilde{V}=J^{-}\cap \widetilde{V}$, hence $\bigcup_{n\geq 1} W_{n}=\widetilde{V}_{1}$.
\qed

\begin{cor}\label{cor:fitcontinuously} The foliation of $U^{-}$ in the set $\widetilde{V}_{1}$ and the lamination of $J^{-}$ inside $\widetilde{V}$ fit together continuously to form a locally trivial lamination of the set $\widetilde{V} \setminus \bigcup_{n\geq  0} H^{n}(B^h)$.

The foliation of $U^{+}$ in the set $\widetilde{V}_{2}$ and the lamination of $J^{+}$ inside $\widetilde{V}$ fit together continuously, to form a locally trivial lamination of the set 
$\widetilde{V} \setminus \bigcup_{n\geq  0} H^{-n}(B^v)$.
\end{cor}

To pass from the description of the foliations of $U^{\pm}$ to the critical locus, the following standard result will become handy (we refer to \cite{LR} for the details of its proof):
\begin{thm}[\cite{LR}]\label{thm:LRcontact}
Consider a pair of holomorphic foliations $\mathcal{F}_{1}$ and $\mathcal{F}_{2}$
defined on some complex two dimensional manifold. Let $\mathcal{C}$ be the critical locus.
If the leaves of $\mathcal{F}_{1}$ and $\mathcal{F}_{2}$ have order of contact two at every point of some
component $X$ of $\mathcal{C}$ then $X$ is smooth, $X$ meets no other component of $\mathcal{C}$, $X$ is
a component of $\mathcal{C}$ with multiplicity one, and $X$ is everywhere transverse to both
foliations.
\end{thm}

In Theorem \ref{prop:critlocus-Sdelta} below, we will prove the fact that throughout the $\HOV_{\beta}$ region, the order of contact of the foliations of $U^{+}$ and $U^{-}$ in the sets $S_{1}$ and $S_{2}$ is two, therefore we can make use of Theorem \ref{thm:LRcontact} to describe the critical locus in these two regions.

\begin{thm}\label{prop:critlocus-Sdelta} Denote by $\mathcal{C}^{v}$ the critical locus in the region 
\begin{equation}\label{def:S1}
	S^{v}:=\left (S_{\delta}^v \cup B^{v} \right) \setminus \left(\bigcup\limits_{n\geq  0} H^{n}(B^h) \cup H^{-1}(\widetilde{V})\cup J^{-}\right).
\end{equation}
$\mathcal{C}^{v}$ is a punctured holomorphic disk tangent at infinity to the $y$-axis, with a Cantor set removed and with punctures at the horizontal boundaries of each of the sets $H^{n}(B^h)$, $n\geq 0$ and $H^{-1}(\widetilde{V})$. It projects one-to-one to the $y$-axis, and is everywhere transverse to the foliations of $U^{+}$ and $U^{-}$, and to the horizontal boundaries of $\widetilde{V}$, $H^{-1}(\widetilde{V})$ and $H^{n}(B^{h})$, $n\geq 0$.

Similarly, if we denote by $\mathcal{C}^{h}$ the critical locus in the region 
\begin{equation}\label{def:S2}
	S^{h}:=\left(S_{\delta}^h\cup B^{h}\right) \setminus \left(\bigcup\limits_{n\geq 0} H^{-n}(B^v)\cup H(\widetilde{V}) \cup J^{+}\right),
\end{equation}
then $\mathcal{C}^{h}$ is a punctured disk tangent at infinity to the $x$-axis, with a Cantor set removed and with punctures at the vertical boundaries of each of the sets $H^{-n}(B^v), n\geq 0$ and $H(\widetilde{V})$. It projects one-to-one to the $x$-axis, and is everywhere transverse to the foliations of $U^{+}$ and $U^{-}$, and to the vertical boundaries of $\widetilde{V}$, $H(\widetilde{V})$, and $H^{-n}(B^{v})$, $n\geq 0$.

In the region $S^{h}\cup S^{v}$, the order of contact of the foliations of $U^{+}$ and $U^{-}$ is two. 
\end{thm}

 \begin{figure}[htb]
\begin{center}
\includegraphics[scale=1]{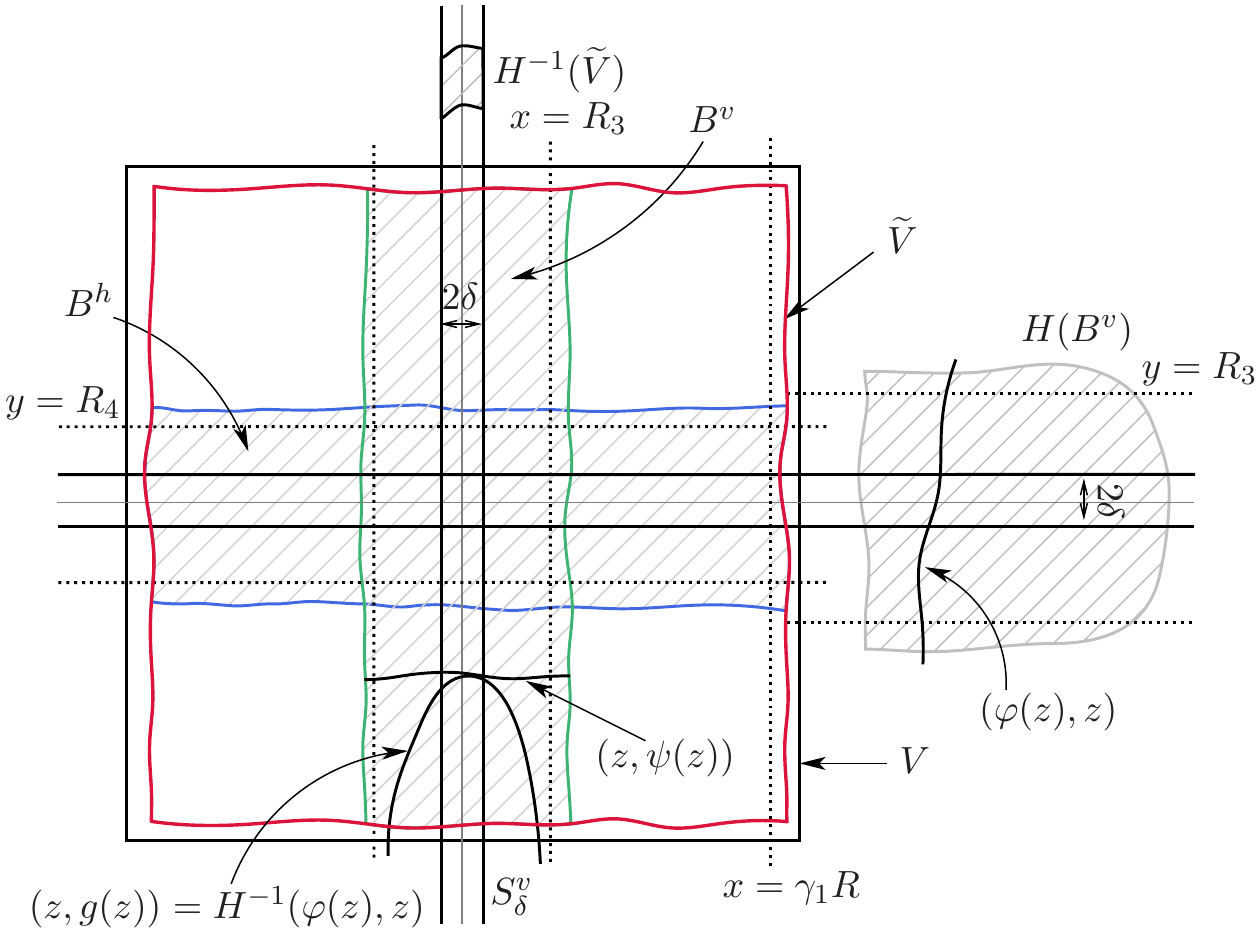}
\end{center}
\caption{The tangency point inside $S^{v}_{\delta}$  between a parabolic-like leaf of the foliation of $U^{+}$ and a horizontal-like leaf of the foliation of $U^{-}$ gives rise to a point on the critical locus $\mathcal{C}^{v}$.}
\label{fig:tangency}
\end{figure}

\proof We will work in the region $S^{v}$. A schematic drawing of the critical locus $\mathcal{C}^{v}$ from the region $S^{v}$ is done in Figure \ref{fig:components}. By Proposition \ref{prop:foliation2}, and equation \ref{eq:gamma1R}, the foliation of $U^{+}$ in the region $S_{\delta}^v \setminus H^{-1}(\widetilde{V})\subset U^{+}$ is given by (vertical) parabolic-like leaves of the form $(z, g(z))=(z,(p(z)-\varphi(z))/a))$, $|z|\leq \delta$, where $ |\varphi'(z)|<\frac{|a|}{2\gamma_{1}R-1}<\frac{|a|}{4|a|+3}<1$. 

The tube $B^{v}$ maps outside $\widetilde{V}$ under one iterate of the H\'enon map, as shown in Figure \ref{fig:tangency}. Since the vertical boundary of $\widetilde{V}$ does not intersect $|x|=\gamma_{1}R$, we can say that $B^{v}$ maps under one iterate in the region $\{|x|>\gamma_{1}R,\ |y|<\gamma_{2}R\}$, and escapes to infinity in forward time. Also, $B^{v}$ contains a tube of width $|x|\geq R_{3}$ inside, where $\gamma_{2}R>R_{3}>1+|a|$, so we can take $|z|\leq R_{3}$.

The foliation of $U^{-}$ in the region  $S^v\cup B^{v}$  
consists of horizontal-like holomorphic disks of the form $(z,\psi(z))$. Note that these are part of long horizontal-like leaves which exit $\widetilde{V}$ through its vertical boundary. Note that all of the removed sets are disjoint from $S_{\delta}^v\setminus B^{v}=S_{\delta}^{v}\setminus \widetilde{V}$. If we look only inside $B^{v}$, we have \[\ds |\psi'(z)|<\frac{1}{2R_{4}-|a|}<\frac{|a|}{|a|+2}<1,\, |z|\leq R_{3}.\] Of course, in $S^{v}_{\delta}\setminus B^{v}$ we have even stronger estimates \[\ds |\psi'(z)|<\frac{1}{2\gamma_{1}R-|a|}<\frac{1}{2R_{4}-|a|},\, |z|\leq \gamma_{1}R,\] so it suffices to work with the weaker ones in both cases.

A horizontal-like leaf and a vertical parabolic-like leaf have a tangency point if and only if $g(z)=\psi(z)$ and $g'(z)=\psi'(z)$. However, the last equation is equivalent to $2z=\varphi'(z)+a\psi'(z)$ and we can count its solutions inside the disk $|z|\leq \delta$ by making use of Rouch\'e's Theorem once again (see Figure \ref{fig:tangency}).

On the boundary $|z|=\delta=\frac{1+|a|}{2}$, we have 
\begin{eqnarray*}
	|2z|=1+|a| &>& |a|\frac{1}{2R_{4}-|a|}+\frac{|a|}{2\gamma_{1}R-1}\\
	&\geq& |a||\psi'(z)|+|\varphi'(z)|\geq |a\psi'(z)+\varphi'(z)|
\end{eqnarray*}
hence  $2z-\varphi'(z)-a\psi'(z)$ has exactly one zero inside $|z|<\delta$, and of course, no zeros on the boundary $|z|=\delta$. Here we do not really need good bounds on $\gamma_{1}, R, R_{4}$; it is enough to have $\varphi'(z)<1$ and $|\psi'(z)|\leq 1$. 
Of course, applying Rouch\'e on any other disk of radius $R_{3}>r>\delta$ will give the same information: there exists exactly one tangency inside the disk of radius $r$, which we had already located inside the disk of radius $\delta$.

It is also clear that a horizontal-like leaf cannot be tangent to more than one parabolic-like leaf of the foliation of $U^{+}$. Thus the critical locus in $S^{v}$ is trapped in the strip $S^{v}_{\delta}$. 

We now prove that the order of contact of these two foliations is two. Suppose that the order of contact between a horizontal-like leaf and vertical parabolic-like leaf is three at a tangency point $z$ with $|z|<\delta$. This implies that $g''(z)=\psi''(z)$, which is equivalent to $2=a\psi''(z)+\varphi''(z)$. In what follows, we show that this is not possible.

Let $\Omega=\{|z-w|=\delta\}$ be a simple closed curve, positively oriented around $z$. $\Omega$ bounds a domain which is contained in the disk of radius $R_{3}$. By Cauchy's integral formula we write
\[
\psi''(z) = \frac{1}{2\pi i}\int_{\Omega}\frac{\psi'(w)}{(w-z)^{2}}\, dw
\]
and use  Proposition \ref{prop:cones} to get 
\[
|\psi''(z)|\leq \frac{\sup_{\Omega}|\psi'(w)|}{\delta}<\frac{1}{2R_{4}-|a|}\cdot\frac{1}{\delta}. 
\]
Similarly, we find $\ds |\varphi''(z)|<\frac{|a|}{2\gamma_{1}R-1}\cdot\frac{1}{\delta}$, which leads to 
\[
|a\psi''(z)+\varphi''(z)|\leq |a||\psi''(z)|+|\varphi''(z)|<\frac{1+|a|}{\delta}=2.
\]
By Theorem \ref{thm:LRcontact}, we know that the critical locus in $S^{v}$
is connected, smooth and transverse on the foliations of $U^{+}$ and $U^{-}$. It implies that the critical locus $C^{v}$ projects one-to-one to the $y$-axis along the horizontal-like leaves of the foliation of $U^{+}$, using the holonomy map between the two transversals, $C^{v}$ and $x=0$.

From Rouch\'e's Theorem it followed that the critical locus does not intersect the vertical boundary of $S_{\delta}^{v}$ given by $|x|=\delta$, so it must intersect the horizontal boundary.  The fact that the horizontal boundaries of $H^{n}(B^{h})$, $n\geq 0$, $\widetilde{V}$, and $H^{-1}(\widetilde{V})$ are all dynamically defined (horizontal-like leaves of the foliation of $U^{-}$) ensures that the critical locus intersects them transversely.
\qed

Using Theorem \ref{prop:critlocus-Sdelta} and Proposition \ref{prop:longleaves} , together with the binary coding in Section \ref{sec:DynamicalBox} and the fact that the $n$-th forward or backward iterate of the H\'enon map is a polynomial mapping of degree $2^{n}$, we can formulate the following corollary, about a substantial part of the critical locus:

\begin{cor}[\textbf{Iterates of the critical locus components $\mathcal{C}^{v}/\mathcal{C}^{h}$}]\label{cor:critlocus-Sdelta} 
Denote by $\mathcal{C}^{v, -n}:=H^{-n}(\mathcal{C}^{v})$, and by $\mathcal{C}^{h,n}:=H^{n}(\mathcal{C}^{h})$, for $n\geq 1$. 

The critical locus $\mathcal{C}^{v,-n}$ projects $n$-to-$1$ on the $y$-axis, along the horizontal-like leaves of the foliation of $U^{-}$ in the region $\widetilde{V}_{1}$, and $\mathcal{C}^{v,-n} \cap \widetilde{V}$ consists of $2^{n}$ connected components, one inside each tube $B^{v}_{\lambda}$, which can be labeled accordingly as  $\mathcal{C}^{v}_{\lambda}$, where $\lambda=\lambda_{-n}\ldots \lambda_{-1}$ belongs to $\{0,1\}^{n}$.  

The critical locus $\mathcal{C}^{h,n}\cap \widetilde{V}$ projects $n$-to-$1$ on the $x$-axis, along the vertical-like leaves of the foliation of $U^{+}$ in the region $\widetilde{V}_{2}$, and $\mathcal{C}^{h,n} \cap \widetilde{V}$ consists of $2^{n}$ connected components, one inside each horizontal tube $B^{h}_{w}$, which can be labeled accordingly as  $\mathcal{C}^{h}_{w}$, where $w=w_{0}\ldots w_{n-1}\in\{0,1\}^{n}$.

 All components $\mathcal{C}^{v,-n}$ and $\mathcal{C}^{h,n}$ are transverse to the foliations of $U^{\pm}$.
\end{cor}

\section{\texorpdfstring{The critical locus in $\widetilde{V}$}{The critical locus in V}}\label{sec:CriticalLocusB}

Denote by $Sing(\mathcal{C})$ the singular set of the critical locus $\Crit$. The singular set is a codimension one complex analytic subvariety of $\Crit$, hence it is just a set of points.  At the end of this section, we will actually show that $Sing(\mathcal{C})=\emptyset$ in the HOV$_{\beta}$ region.

We have already shown the critical locus to be smooth in certain regions, for example in the vertical and horizontal tubular regions $S^{v}$ and $S^{h}$ of Theorem \ref{prop:critlocus-Sdelta}. 
By Theorem \ref{thm:FundamentalDomain}, a fundamental region for the critical locus is $\Crit \cap F$, where 
\[
F=(\widetilde{V}\cup S_{\delta}^h)\backslash H(\widetilde{V}).
\] 
Outside $\widetilde{V}$, the set $S^{h}$ is equal to $S_{\delta}^h \setminus H(\widetilde{V})$, where the critical locus is included in the smooth component $\mathcal{C}^{h}$, hence we are left to show that $\mathcal{C}$ is smooth in 
$\widetilde{V}\setminus H(\widetilde{V})$.

By Lemma \ref{prop:critlocus-Sdelta} and Corollary \ref{cor:critlocus-Sdelta} , we have a complete description of the critical locus in $\widetilde{V}$, outside the intersections of the vertical-like and horizontal-like tubes $H^{-k}(B^v)\cap H^m(B^h)$ where $k,m\geq 0$: these are just backward iterates of the irreducible component $\mathcal{C}^{v}$ from $S^{v}$ and forward iterates of $\mathcal{C}^{h}$ from $S^{h}$, hence smooth.

The description of the critical locus in the intersection sets $H^{-k}(B^v)\cap H^m(B^h)$, $k,m\geq 0$ will follow from Ehresmann's Theorem \ref{thm:Ehresmann}. 
Of course, by dynamics, it suffices to consider only intersections of the form 
\begin{equation}\label{def:int}
B^{h}\cap H^{-k}(B^v)\subset S^{h} \mbox{ where } k\geq 0.
\end{equation}

In the HOV region, the intersection in Equation \ref{def:int} consists of $2^{k}$ distinct sets, labeled $X_{\lambda}=\Delta_{(\lambda,\emptyset)}$, where $\lambda=\lambda_{-k}\ldots \lambda_{-1}\in\{0,1\}^{k}$ when $k\geq 1$, or $\lambda=\emptyset$ when $k=0$ (see the binary coding in Section \ref{sec:DynamicalBox}, Equation \eqref{eq:notations}).
The vertical boundary of each $X_{\lambda}$ is a subset of the vertical boundary of $H^{-k}(B^{v})$, hence laminated by vertical-like leaves of the foliation of $U^{+}$. The horizontal boundary of $X_{\lambda}$ is a subset of the horizontal boundary of $B^{h}$, hence laminated by the horizontal-like leaves of the foliation of $U^{-}$. Let us denote by $\Crit_{X_{\lambda}}$ the critical locus in $X_{\lambda}$, that is $\Crit_{X_{\lambda}}:=\, X_{\lambda}\cap\Crit$.

It is worth pointing out here that the analysis above shows that another fundamental region of the critical locus is given by
\begin{eqnarray}\label{def: thirdFundamentalRegion}
(\mathcal{C}^{v}\cap \widetilde{V}) \cup \mathcal{C}^{h}\cup \bigcup_{\lambda\in \{0,1\}^{n}, n\geq 0} \Crit_{X_{\lambda}},
\end{eqnarray}
a fact which we will also exploit in Section \ref{sec:criticalModel} when building topological models for the critical locus.
 
The prototype for the critical locus $\Crit_{X_{\lambda}}$ will be the critical locus in the set $B=S^{v}\cap S^{h}$ from Equation \eqref{eq:B}.

Recall that by Proposition \ref{prop:foliation2} we know that in $B$, the foliation of $U^{+}$ is vertical parabolic-like and the foliation of $U^{-}$ is horizontal parabolic-like. Therefore we could in principle describe the critical locus in $B$ by hand, as the set of tangencies of two families of parabolas. However, we would like to give a more general argument that works for all the regions $X_{\lambda}$.

Let $\chi_{c,a}$ be the Euler characteristic of $\Crit_{X_{\lambda}}$ and $\mu_{c,a}$ be the sum of the Milnor numbers of singularities of the critical locus $\Crit_{X_{\lambda}}$ for the H\'{e}non map $H_{c,a}.$ 
We will use Theorem \ref{prop:critlocus-Sdelta}, the Ehresmann Fibration Theorem and its Corollary \ref{prop:topology} to describe the critical locus $\Crit_{X_{\lambda}}$.

\begin{thm}[\textbf{Ehresmann Fibration Theorem} \cite{D}]\label{thm:Ehresmann} Suppose that $M$ and $N$ are smooth manifolds and that $f:M\to N$ is a proper smooth submersion. Then $f$ is a locally trivial smooth fibration.
If $M$ is a manifold with boundary $\partial M$ and if both $f:M\to N$ and $\restr{f}{\partial M}:\,\partial M\to N$ are submersions, then both $f$ and $\restr{f}{\partial M}$ are locally trivial fibrations.
\end{thm}

\begin{cor}\label{prop:topology} Let $M$ be a complex $3$-dimensional manifold with smooth boundary $\partial M$, and $f:M\to \C^2$ a holomorphic map. Assume that $\restr{f}{\partial M}:\partial M \to \C^2$ is a proper submersion. Assume further that each level set has at most finitely many critical points and not all level sets are critical. 
Let $L=f^{-1}(\epsilon)$ be a level set of $f$, and denote by $\chi_L$ its Euler characteristic, and by $\mu_L$ the sum of the Milnor numbers of the singularities on $L$. Then $\mu_L-\xi_L$ is an invariant which does not depend on $L$.
\end{cor}

\proof We apply the Ehresmann Fibration Theorem to the map $f:f^{-1}(U_{\epsilon})\to U_{\epsilon}$, where $\epsilon=(\epsilon_1, \epsilon_2)\in\C^{2}$ is a non-singular value, and $U_{\epsilon}$ is a small enough neighborhood of $\epsilon$ (such level sets exist by our assumption). By Theorem \ref{thm:Ehresmann}, all level sets $f^{-1}(\epsilon')$ where $\epsilon'\in U_{\epsilon}$ are diffeomorphic, and their Euler characteristics are equal. Since they are non-singular, the sums of their Milnor numbers are equal to $0$. Thus, $\xi_L$ is a local constant on the set of non-singular values. The space of singular values has codimension at least one, so it does not separate the space. Therefore, the set of non-singular values is connected. Hence, $\xi_L-\mu_L$ is an invariant on the set of non-singular values.

Let $\epsilon^0=(\epsilon^0_1,\epsilon^0_2)$ and assume that $L=f^{-1}(\epsilon^0)$ is a singular level set. Since the set of singular values has codimension one, there is a plane $\Lambda$ going through $\epsilon^0$ so that $\epsilon^0$ is an isolated critical value on this plane.  For a singular point $a$, let us take a sphere $S$ around it, such that nearby level sets $\epsilon\in \Lambda$ are transverse to $S$. The Euler characteristic of $L$ inside $S$ is $\chi(L\cap int(S))=1$. Denote the Milnor number of $a$ by $\mu(a)$. The nearby level sets $L'$ are homotopy
equivalent to a bouquet of $\mu(a)$ one-dimensional spheres. Thus, $\chi(L\cap int(S))-\mu(L\cap int(S))=1$. The Euler characteristic of the intersection of a level set with the sphere is equal to $0$. Outside of the spheres, the nearby level sets are diffeomorphic to each other by the Ehresmann Fibration Theorem. Hence, $\chi(L)-\mu(L)$ is a global invariant.
\qed

\begin{thm}\label{thm: zero} All non-singular critical loci $\Crit_{X_{\lambda}}$ for $(c,a)\in \HOV_{\beta}$ are diffeomorphic. Moreover, $\mu_{c,a}-\xi_{c,a}=0$ in $B_{w}$ for all $(c,a)\in \HOV_{\beta}$. 
\end{thm}

\proof Region $X_{\lambda}$ depends holomorphically on $c$ and $a$ in the $\HOV_{\beta}$-region, but we will not add an index to mark this dependency, in order to simplify notations. We fix $c$ and let $a$ vary, and we omit $c$ in the subscripts.  Since for every $a$, there exists $R$ such that $|a|<R$, it is enough to prove the statement on the sets $\D_{r_{0}}=\{a:|a|<r_{0}\}$, where $r_{0}$ is arbitrary. Consider the manifold $M'$ that is a fibration over the disk $\D_{r_{0}}$, with fiber the region $X_{\lambda}$.  By Equation \eqref{eq:definingfunction}, the critical locus $\mathcal{C}_{a}$ is the zero-set of a holomorphic function $w_{a}$. We set $f(x,y,a)=(w_a(x,y), a)$. Since the critical locus $\mathcal{C}_{a}$ is transverse to the boundary of its corresponding $X_{\lambda}$,  there exists $\epsilon_0$ such that the level set $f^{-1}(\epsilon, a)$ is transverse to the boundary of $X_{\lambda}$ for all $a\in {\mathbb D}_{r_{0}}$ and all $\epsilon \in {\mathbb D}_{\epsilon_0}$. Let $M=f^{-1}({\mathbb D}_{r_{0}}\times {\mathbb D}_{\epsilon_0})$. 
By considering the case $a=0$, we see that not all level sets are singular.  Let $\epsilon=(\epsilon_1,\epsilon_2)$ and assume that $L=f^{-1}(\epsilon_1,\epsilon_2)$ is a non-singular level set. Let $M=f^{-1}(U_{\epsilon})$, where $U_{\epsilon}$ is a neighborhood of $\epsilon$.  Applying Ehresmann Fibration Theorem to $M$, we see that non-singular level sets close to $L$ are diffeomorphic to $L$. The set of non-singular level sets is connected, hence all non-singular level sets are diffeomorphic. Let $L$ be a singular level set. Let $a$ be a critical point on $L$, with Milnor number $\mu(a)$. There exists a one dimensional curve such that on this curve $a$ is an isolated critical point. Then all nearby level sets are a collection of $\mu-1$ one-dimensional spheres.

We apply Corollary \ref{prop:topology} to $f:M\to \C^2$, and get that $\mu(a)-\xi(a)$ is a global invariant. 
We calculate the invariant $\mu_{c,a}-\xi_{c,a}=0$ for $a=0$ and arbitrary $c$ in the region $\HOV_{\beta}$. To simplify notations, we denote it by $\mu_0-\chi_0$ and show that it is equal to $0$. When $a=0$, the critical locus $\Crit_{X_{\lambda}}$ is the union of two intersecting disks \cite{F}. Hence $\mu_0=1$ and $\chi_0=1$, which implies that $\mu_0-\chi_0=0$. Hence $\mu_{c,a}-\xi_{c,a}=0$ for all $(c,a)\in \HOV_{\beta}$. 
\qed

\begin{thm}\label{thm:handle} Let $\Crit$ be a complex analytic set in a neighborhood of a polydisk $\Delta$ in $\C^{2}$, with $\mu(C)-\chi(C)=0$. Assume that $\Crit$ is smooth in a neighborhood of the boundary of $\Delta$, and that the boundary $\partial (\Crit\setminus \Delta)$ is the disjoint union of two real analytic curves $\mathcal{E}^{h}$ and $\mathcal{E}^{v}$, each homeomorphic to a circle, such that $\mathcal{E}^{h}$ is a subset of the horizontal boundary of $\Delta$ and $\mathcal{E}^{v}$ is a subset of the vertical boundary of $\Delta$.  

Suppose further that the set of singular points of $\Crit$ is nonempty.  Then $\Crit$ in the region $\Delta$ is the union of two holomorphic disks intersecting at one point.
\end{thm}

\proof
Each connected component of $\Crit\setminus Sing(\Crit)$ inside $\Delta$ has a boundary component on $\partial \Delta$. Therefore, there are at most two irreducible components of $\Crit$ inside $\Delta$ . So we consider these two cases:
\begin{enumerate} 
\item There is one irreducible component of $\mathcal C$ in $\Delta$ that has two boundary components $\mathcal{E}_{h}$ and $\mathcal{E}_{v}$. But then $\mathcal C \cap \Delta$ is smooth, a contradiction. 

The fact that $\Crit$ needs to be smooth can also be seen by computing the Milnor numbers. Since $\partial ({\mathcal C} \setminus \Delta)$ has two connected components, its Euler characteristic is less than or equal to $0$. On the other hand, the Milnor number is greater than or equal to 0. Since $\mu(C)-\chi(C)=0$, both the Euler characteristic and the sum of Milnor numbers of singularities have to be $0$. This is a contradiction with the fact that $Sing(\Crit)$ is assumed non-empty in $\Delta$.
\item There are two irreducible components of $\mathcal C$ in $\Delta$. Since each component has a boundary on $\partial\Delta$, the Euler characteristic of each component is at most one. In order for the component of the critical locus to be non-smooth, the sum of the Milnor numbers of the singularities has to be at least $1$. Hence, $\mu(C)$ must be equal to $1$, and the two connected components must have Euler characteristic $1$, and hence be homeomorphic to disks. Therefore $\Crit\cap \Delta$ is the union of two disks intersecting transversally in one point. The transversality property follows from the fact that a complex variety can never be a differentiable manifold (not even of class $C^{1}$, see \cite{M2}) throughout a neighborhood of a singular point.
\qed
\end{enumerate}

We can use the dynamics of the H\'enon map and the properties of the foliations of $U^{+}$ and $U^{-}$, to show that the critical locus inside each of the sets $X_{\lambda}$ from Equation \eqref{def:int} cannot be the union of two disks intersecting in one point. We will start with the polydisk $B=X_{\emptyset}$.

\begin{figure}[htb]
\begin{center}
\includegraphics[scale=1]{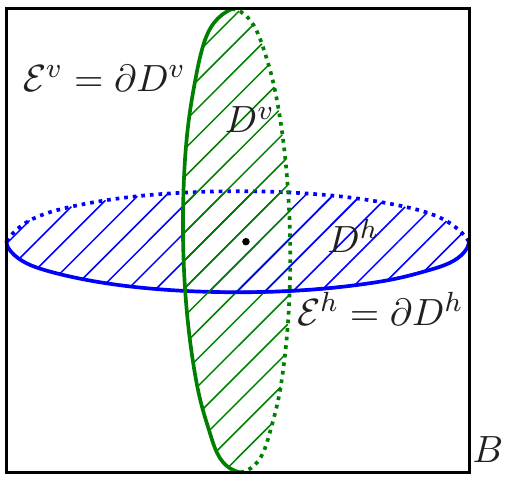}
\end{center}
\caption{Intersection of two disks of the critical locus inside $B$.}
\label{fig:intersectingdisks}
\end{figure}

\begin{lemma}\label{lemma:NOsing} The critical locus $\Crit$ in the polydisk $B$ cannot be  the union of two holomorphic disks intersecting at one point, $D^{v}$ with boundary $\mathcal{E}^{v}$ on the vertical boundary of $B$, and $D^{h}$ with boundary $\mathcal{E}^{h}$ on the horizontal boundary of $B$, as in Figure \ref{fig:intersectingdisks}.
\end{lemma}
\proof The set $B=B^{v}\cap B^{h}$ maps outside $\widetilde{V}$ under one iterate (forward or backward) of the \He map, therefore on $B$ the functions $(\varphi^{+})^{2}$ and $(\varphi^{-})^{2}$ are well defined holomorphic functions. 

On the vertical boundary of $B$, the leaves of the foliation of $U^{+}$ are vertical-like. Any such leaf is a level-set $(\varphi^{+})^{2}=\xi$ for some value of $\xi\in\C-\D$, so the gradient $\nabla (\varphi^{+})^{2}$ is equal to $0$ along the level-set.
A tangent vector to this leaf at some point $(x,y)$ on the circle $\mathcal{E}^{v}$
is  perpendicular to the gradient line, hence it is a scalar multiple of 
\[ 
\left(-\partial_{y}(\varphi^{+})^{2}(x,y), \partial_{x}(\varphi^{+})^{2}(x,y)\right).
\] 
The tangent vector belongs to the vertical cone at $(x,y)$, hence 
\[
\left|\partial_{y}(\varphi^{+})^{2}(x,y)\right| < k\left|\partial_{x}(\varphi^{+})^{2}(x,y)\right|,
\] 
where $k<1$ is a fixed constant strictly less than one  whose value depends only on the HOV$_{\beta}$ region. Inside the vertical tube $B^{v}$, the foliation of $U^{+}$ can have horizontal tangent lines, hence  $\partial_{x}(\varphi^{+})^{2}(x,y)$ will be equal to $0$ for some points $(x,y)$ inside $B$. However, we claim that  this does not happen on the disk $D^{v}$ of the critical locus. In the region $B^{h}$, the foliation of $U^{-}$ does not admit any horizontal tangent lines, by Lemma \ref{lemma:NOhorver}. Moreover, on $D^{v}$ the tangent lines to the foliations of $U^{+}$ and $U^{-}$ coincide, hence
\[
\partial_{x}(\varphi^{+})^{2}(x,y)\neq 0,\ \ \mbox{for all}\ \  (x,y)\in D^{v}.
\]  
Therefore the function 
\[
\ds g_{1}:=\frac{\partial_{y}(\varphi^{+})^{2}}{\partial_{x}(\varphi^{+})^{2}}\] 
is well defined and holomorphic on the closure of the disk $D^{v}$. By the Maximum Modulus principle, we have $|g_{1}(x,y)|<k<1$ for all $(x,y)\in D^{v}$, since the inequality is satisfied on the boundary $\mathcal{E}^{v}$. This is equivalent to saying that the foliation of $U^{+}$ at any point on the disk $D^{v}$ is vertical-like.

Take now any point $(x,y)$ on the circle $\mathcal{E}^{h}$. This circle is part of the horizontal boundary of $B$, which is laminated by horizontal-like leaves of the foliation of $U^{-}$. Therefore, a tangent vector at $(x,y)$ to the foliation of $U^{-}$ is of the form 
\[
\left(-\partial_{y}(\varphi^{-})^{2}, \partial_{x}(\varphi^{-})^{2}\right), \ \ \mbox{with}\ \ \left|\partial_{x}(\varphi^{-})^{2}\right|<k\left|\partial_{y}(\varphi^{-})^{2}\right|. 
\]  
Inside the vertical tube $B^{v}$, the foliation of $U^{-}$ admits horizontal tangent lines, hence  $\partial_{y}(\varphi^{+})^{2}(x,y)$ will be equal to $0$ for some points $(x,y)$ in $B$. However, by Lemma \ref{lemma:NOhorver}, in the vertical tube $B^{v}$, the foliation of $U^{+}$ admits no vertical tangent lines. Since the disk $D^{h}$ belongs to the critical locus, the foliation of $U^{-}$ will have no vertical tangent lines either on $D^{h}$, hence $\partial_{y}(\varphi^{-})^{2}\neq 0$
 on $D^{h}$. Therefore the function 
\[
\ds g_{2}:=\frac{\partial_{x}(\varphi^{-})^{2}}{\partial_{y}(\varphi^{-})^{2}}\] 
is well defined and holomorphic on the closure of the holomorphic disk $D^{h}$. By the Maximum Modulus principle, we have $|g_{2}(x,y)|<k<1$ for all $(x,y)\in D^{h}$, since the inequality is satisfied on the boundary $\mathcal{E}^{h}$. This is equivalent to saying that the foliation of $U^{-}$ at any point on the disk $D^{h}$ is horizontal-like.

Assume that the critical locus inside $B$ is the union of the intersecting disks $D^{v}$ and $D^{h}$. The foliation of $U^{+}$ is vertical-like on $D^{v}$, whereas the foliation of $U^{-}$ is horizontal-like on $D^{h}$, hence the two foliations cannot have a common tangent line at the intersection point of $D^{v}$ and $D^{h}$. Hence the intersection of the two disks cannot belong to the critical locus, which is a contradiction. 
\qed

\begin{remark}
In the set $B$ we encounter a symmetry between $\varphi^{+}$ and $\varphi^{-}$ which cannot be reproduced in the other sets $X_{\lambda}$ from Equation \ref{def:int} when $|\lambda|=n>0$. Each $X_{\lambda}$ is a subset of $B^{h}$, so the foliation of $U^{-}$ does not have any horizontal tangent lines, by Lemma \ref{lemma:NOhorver}. However, the foliation of $U^{+}$ can and will admit vertical tangent lines (for example in the case of $X_{0}$ and $X_{1}$, the two connected components of $H^{-1}(B^{v})\cap B^{h}$, these vertical tangents are the backward images of the horizontal tangents to the foliation of $U^{+}$ inside the vertical tube $B^{v}$). So we need to adjust the argument of Lemma \ref{lemma:NOsing} and combine it with a Hartogs Extension Lemma in order to claim the more general statement formulated in Lemma \ref{lemma:NOsingN} about the critical locus $\mathcal{C}_{X_{\lambda}}$.
\end{remark}

The following short topological digression on linking numbers will be useful to us. Let $C_1$, $C_2$ be two smooth oriented disjoint circles in $\mathbb S^3$ without self-intersections. 

Let $D_1\subset {\mathbb S}^3$ be a smooth disk with boundary $C_1$, oriented so that the positive orientation on $D_1$ induces an orientation on $C_1$.  Assume that $D_1$ intersects $C_2$ transversally. At each point $p$ of intersection we consider the basis $\{v_1,v_2,v_3\}$, where $\{v_1, v_2\}$ form the positive basis of $TD_1$ and $v_3$ defines the positive orientation on $C_2$. If $v_1,v_2$ and $v_3$ define the positive orientation of $S^3$, we say that $sign(p)=+1$, otherwise $sign(p)=-1$. The linking number of $C_1$ and $C_2$ is equal to 
\[
\mbox{Link} (C_1,C_2)=\sum_{p\in D_{1}\cap C_{2}} sign(p).
\]

The sphere ${\mathbb S}^3$ is the boundary of the disk ${\mathbb D}^3$. Let $D_1$ and $D_2$ be two disks in ${\mathbb D}^3$ with boundaries $C_1$ and $C_2$ such that the positive orientation on $D_1$ and $D_2$ induces a positive orientation on $C_1$ and $C_2$. Assume that $D_1$ intersects $D_2$ transversally. At each point $p$ of this intersection we consider the following basis $\{v_1, v_2, v_3, v_4\}$,  where $\{v_1, v_2\}$ form a positive basis of $TD_1$ at the point $p$, while $\{v_3, v_4\}$ form a positive basis of $TD_2$ at the point $p$. If $\{v_1, v_2, v_3, v_4\}$ is a positive basis of $\D^3$, then we say that $sign(p)=+1$, otherwise $sign(p)=-1$.

 \begin{lemma} 
 \[ 
 \mbox{Link} (C_1,C_2)=\sum_{p\in D_{1}\cap D_{2}} sign(p)
 \]
\end{lemma}

\begin{cor} Let $C_1,$ $C_2$ be non self-intersecting circles at horizontal and vertical boundaries of a polydisk ${\mathbb D}\times \mathbb D$. Then their linking number is equal to $1$.
\end{cor}

For a complex manifold $M$, the complex structure induces a positive orientation: if $\{e_1, \ldots, e_n\}$ is a basis in the complex tangent space of $M$ at a point $p$, then naturally the set $\{e_1, ie_1,\ldots, e_n, ie_n\}$ is a positively oriented basis in the real tangent space. If $V,W\subset M$ are complex submanifolds of complimentary dimension of a manifold $M$, then the complex structures on $V$ and $W$ are induced by the complex structure on $M$. If $\{e_1, \ldots, e_k\}$ is a basis in the complex tangent space to $V$,  $\{e_{k+1}, \ldots, e_n\}$ 
is a basis of $W$, then $\{e_1, ie_1, \ldots, e_k, ie_k\}$ is a positive oriented basis for $TV$, $\{e_{k+1},\ldots, e_n\}$ is a positive basis for $TW$ and consequently $\{e_1, ie_1,\ldots, e_n, ie_n\}$ is a positive basis for $TM$.

\begin{lemma}\label{lemma:linking} Let $P$ be a complex manifold which is a smooth image of a polydisk ${\mathbb D}\times \mathbb D$. Let $C_1$ and $C_2$ be non self-intersecting circles at the images of horizontal and vertical boundaries. Let $D_1$ and $D_2$ be complex disks with boundaries $C_1$, respectively $C_2$. Then $D_1$ and $D_2$ intersect at exactly one point. 
\end{lemma}

\begin{lemma}[\textbf{Critical Locus in $X_{\lambda}$}]\label{lemma:NOsingN} Let $\lambda\in\{0,1\}^{n}, n>0$. The critical locus $\Crit$ intersects the boundary of the $X_{\lambda}$ transversely in two disjoint circles, $\mathcal{E}^{v}$ on the vertical boundary of $X_{\lambda}$ and $\mathcal{E}^{h}$ on the horizontal boundary of $X_{\lambda}$. 
The critical locus $\Crit$ inside $X_{\lambda}$ cannot be the union of two intersecting disks, $D^{v}$ with boundary $\mathcal{E}^{v}$, and $D^{h}$ with boundary $\mathcal{E}^{h}$.
\end{lemma}
\proof The set $X_{\lambda}=B^{v}_{\lambda}\cap B^{h}$ maps outside $\widetilde{V}$ under one backward iterate, respectively under $n+1$ forward iterates of the \He map, therefore on $X_{\lambda}$ the functions $(\varphi^{-})^{2}$ and $(\varphi^{+})^{2^{n+1}}$ are well defined holomorphic functions. Note that $(\varphi^{+})^{2^{n+1}}=(\varphi^{+})^{2}\circ H^{n}$ is well defined and holomorphic in $H^{-n}(B^{v})$, in particular in the entire tube $B^{v}_{\lambda}$. 
Likewise, $(\varphi^{-})^{2}$ is well defined and holomorphic in the entire tube $B^{h}$.

Any leaf of the foliation of $U^{+}$ inside $B^{v}_{\lambda}$ is part of a level-set $(\varphi^{+})^{2^{n+1}}=\xi$ for some value of $\xi\in\C-\D$. 
A tangent vector to this leaf at some point $(x,y)$
is  perpendicular to the gradient line, hence it is of the form 
\[
\left(-\partial_{y}(\varphi^{+})^{2^{n+1}}(x,y), \partial_{x}(\varphi^{+})^{2^{n+1}}(x,y)\right).
\]

The vertical boundary $B^{v}_{\lambda}$ is laminated by vertical-like leaves of the foliation of $U^{+}$ and we denote it by $\partial_{v}(B^{v}_{\lambda})$. If we take any point $(x,y)$ on $\partial_{v}(B^{v}_{\lambda})$, 
then a tangent vector at $(x,y)$ to the foliation of $U^{+}$ will belong to the vertical cone at $(x,y)$, hence it will satisfy the inequality
\begin{equation}\label{ineq:fol}
\left|\partial_{y}(\varphi^{+})^{2^{n+1}}\right|<k\left|\partial_{x}(\varphi^{+})^{2^{n+1}}\right| ,
\end{equation}
where $k<1$ is a fixed constant strictly less than one  whose value depends only on the HOV$_{\beta}$ region. In particular, we must have 
\[
\partial_{x}(\varphi^{+})^{2^{n+1}}\neq 0 \mbox{ on } \partial_{v}(B^{v}_{\lambda}).
\]
In fact, inequality  \eqref{ineq:fol} holds not only on $\partial_{v}(B^{v}_{\lambda})$, but in a larger neighborhood of the vertical boundary of $X_{\lambda}$, since the foliation of $U^{+}$ is vertical-like in the part of the tube $B^{v}_{\lambda}$ corresponding to the $n^{th}$ preimage of the entire region $\mathcal{R}$ enclosed between the purely straight tube $S^{v}_{\delta}\subset B^{v}$ and $\partial_{v}B^{v}$ inside $H^{n}(B^{v}_{\lambda})\subset B^{v}$ (recall Propositions  \ref{prop:dyntubes} part (a) and  \ref{prop:foliation1}).  

By Lemma \ref{lemma:NOhorver}, in the entire horizontal tube $\overline{B^{h}}$, the foliation of $U^{-}$ admits no horizontal tangent lines. Since the closed disks $\overline{D^{h}}$ and $\overline{D^{v}}$ belong to the critical locus, the foliation of $U^{+}$ will have no horizontal tangent lines here either, hence 
\[
\partial_{x}(\varphi^{+})^{2^{n+1}}\neq 0 
 \mbox{ on } \overline{D^{h}}\cup \overline{D^{v}}. 
 \]
 Therefore, by continuity, the function 
\[\ds g(x,y):=\frac{\partial_{y}(\varphi^{+})^{2^{n+1}}(x,y)}{\partial_{x}(\varphi^{+})^{2^{n+1}}(x,y)}\] 
is well defined and holomorphic in a neighborhood of $\overline{D^{h}}\cup \overline{D^{v}}\cup \partial_{v}(B^{v}_{\lambda})$. 

Just as in the proof of Lemma \ref{lemma:NOsing}, we can apply the Maximum Modulus Principle on the disk $\overline{D^{v}}$ to infer that $|g(x,y)|<k<1$ for all $(x,y)\in D^{v}$, since the inequality is satisfied on the boundary $\mathcal{E}^{v}$. This is equivalent to saying that the foliation of $U^{+}$ at any point on the disk $D^{v}$ is vertical-like. As a consequence, the disk $\overline{D^{v}}$ cannot intersect any part where the foliation of $U^{-}$ is horizontal-like, in particular, it must be bounded away from the horizontal boundary of $X_{\lambda}$, so 
$\overline{D^{v}}$ is finitely sheeted over the first coordinate.

The key observation now is that $\partial_{v}(B^{v}_{\lambda}) \cup \overline{D^{v}}$ is a {\it distorted Hartogs figure} in $\C^{2}$. Therefore, by Hartogs Lemma (see e.g. \cite{Chi1}, \cite{G}), $g$ has a unique holomorphic extension to the set 
\begin{equation}\label{eq:hartogs}
B^{v}_{\lambda}\cap \left\{|y|\leq \gamma_{1}R\right\}.
\end{equation}
Note here that by Proposition \ref{prop:dynV} (see also Figures \ref{fig:horseshoe} and \ref{fig:VandTube}), the set $X_{\lambda}$ is compactly contained in the open straight polydisk $V_{1}$ of radius $\gamma_{1}R$, therefore the disk $\overline{D^{v}}$ does not intersect the set $|y|=\gamma_{1}R$.

By the Maximum Modulus principle in  horizontal disks in $B^{v}_{\lambda}$, we get $|g(x,y)|<k<1$ for all $(x,y)\in B^{v}_{\lambda}\cap \{|y|\leq \gamma_{1}R\}$, since the inequality is satisfied on the vertical boundary $\partial_{v}(B_{\lambda})$. This implies that the foliation of $U^{+}$ at any point in $X_{\lambda}$ is vertical-like.

However, on the horizontal boundary of $X_{\lambda}$, the leaves of the foliation of $U^{-}$ are horizontal-like. Since the circle $\mathcal{E}^{h}$ belongs to critical locus, and to the horizontal boundary of $X_{\lambda}$, it follows that the foliation of $U^{+}$ is also horizontal-like at any point on $\mathcal{E}^{h}$. This implies that the tangent vector to the foliation of $U^{+}$ at any point on $\mathcal{E}^{h}$ belongs to the corresponding horizontal cone, hence $|g(x,y)|>1$ on $\mathcal{E}^{h}$. 
We have reached a contradiction with inequality \eqref{ineq:fol}. This contradiction shows that the critical locus inside $X_{\lambda}$ cannot be the union of two disks, which concludes our proof.

In what follows we emphasize why $\partial_{v}(B^{v}_{\lambda}) \cup \overline{D^{v}}$ is a distorted Hartogs figure. This is easier to see if we map it forward under $H^{n}$ in the vertical tube $B^{v}$.

\begin{figure}[htb]
\begin{center}
\includegraphics[scale=1]{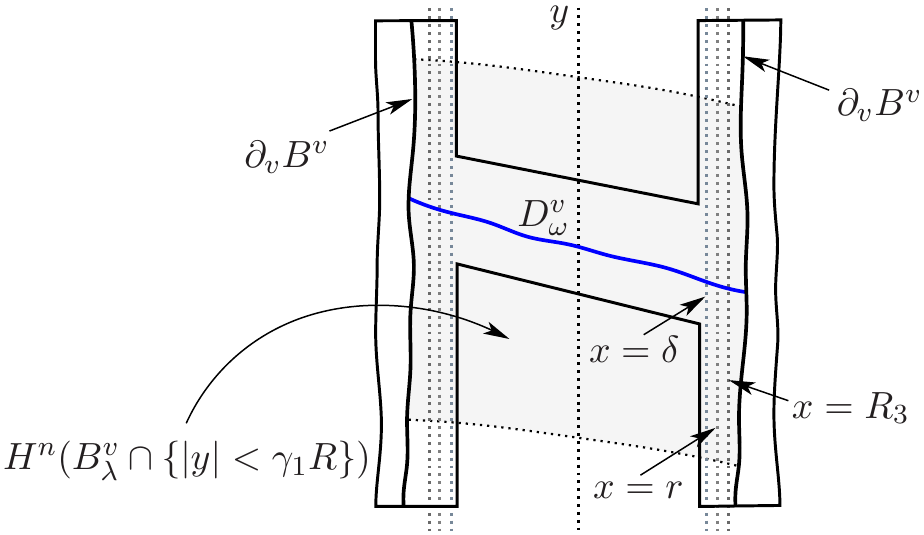}
\end{center}
\caption{A distorted Hartogs Figure formed by the disk $D^{v}_{w}$ of the critical locus together with the vertical boundary of $B^{v}$ foliated  by vertical-like leaves of the foliation of $U^{+}$.}
\label{fig:Hartogs}
\end{figure}

Consider the string $w=w_{0}\ldots w_{n-1}\in\{0,1\}^{n}$, where  $w_{0}=\lambda_{-n},\ldots, w_{n-1}=\lambda_{-1}$. With the binary coding of Section  \ref{sec:DynamicalBox}, we have
\begin{equation*}
Y_{w} = H^{n}(X_{\lambda})=B^{v}\cap B^{h}_{w}\ \subset\ B^{v}\cap T^{h}_{w} = H^{n}(B^{v}_{\lambda})
\end{equation*}
Keeping the same coding, we denote $D^{v}_{w}:=H^{n}(D^{v})$ and $\mathcal{E}^{v}_{w}:=H^{n}(\mathcal{E}^{v})$. Note that $D^{v}_{w}$ is a subset of $Y_{w}$, disjoint from the horizontal boundary of $Y_{w}$. $D^{v}_{w}$ intersects transversely the vertical-like foliation of $U^{+}$ in the region $\mathcal{R}$ (see Figure \ref{fig:Hartogs}), and $\mathcal{E}^{v}_{w}$ projects one-to-one to the $x$-axis along the leaves of the foliation of $U^{+}$ via the holonomy map.  

There exists $\delta\leq r\leq R_{3}$ such that the critical $D^{v}_{w}$ intersects transversely the cylinder $|x|=r$. Let $i:\D\rightarrow D^{v}$ be an injective holomorphic map parametrizing the disk $D^{v}$ and let $(\alpha_{1}(z),\alpha_{2}(z))=H^{n}\circ i(z)$.  Assume by contradiction that for every $r$ there exists $z_{r}$ with $|\alpha_{1}(z_{r})|=r$ such that $\alpha_{1}'(z_{r})=0$; then $\alpha'_{1}(z)=0$ on a set which contains an accumulation point, therefore $\alpha_{1}'$ is constantly equal to $0$ on $\D$, hence $\alpha_{1}$ is a constant and the disk $D^{v}_{w}$ is a vertical disk. However, this is impossible, since $D^{v}_{w}$ does not intersect the horizontal boundary of $Y_{w}$.

The fact that $D^{v}_{w}$ is a holomorphic disk in $Y_{w}$ whose closure does not intersect the horizontal boundary of $Y_{w}$ also implies that the restriction of the projection $\pi(x,y)=x$ to $D^{v}_{w}$ is a proper holomorphic mapping. In particular,  $\pi: D^{v}_{w}\cap \{|x|\leq r\} \rightarrow \D_{r}$ is a branched covering map over the disk $\D_{r}$ of finite degree $m>0$, which is a covering map over the circle of radius $r$.
By Corollary \ref{cor:critlocus-Sdelta} $\mathcal{E}^{v}_{w}$ projects one-to-one to the $x$-axis along the leaves of the vertical-like foliation of $U^{+}$. By Lemma \ref{lemma:linking} it follows that $m=1$, hence $D^{v}_{w}$ is the graph of a holomorphic function $\xi$, of the form $\{(x,\xi(x)), x\in\D_{r}\}$.

The fact that $D^{v}_{w}$ is a holomorphic disk in $Y_{w}$ whose closure does not intersect the horizontal boundary of $Y_{w}$ also implies that the restriction of the projection $\pi(x,y)=x$ to $D^{v}_{w}$ is a proper holomorphic mapping. In particular,  $\pi: D^{v}_{w}\cap \{|x|\leq r\} \rightarrow \D_{r}$ is a branched covering map over the disk $\D_{r}$ of finite degree $m>0$, which is a covering map over the circle of radius $r$.
By Corollary \ref{cor:critlocus-Sdelta} $\mathcal{E}^{v}_{w}$ projects one-to-one to the $x$-axis along the leaves of the vertical-like foliation of $U^{+}$. By Lemma \ref{lemma:linking} it follows that $m=1$, hence $D^{v}_{w}$ is the graph of a holomorphic function $\xi$, of the form $\{(x,\xi(x)), x\in\D_{r}\}$. 
 
The  complex-valued function 
\[ 
T(x,y):=g\circ H^{n}(x,y)
\] 
is well defined and holomorphic on $\overline{D^{v}_{w}}\cup \mathcal{R}$. It is worth mentioning that even if we have a straightforward composition rule for the rate of escape function $\varphi^{+}\circ H=(\varphi^{+})^{2}$, this does not extend to the ratio of partial derivatives of $\varphi^{+}$, so we do not have any nice reduction formula for $g\circ H$ to work with.

Consider now a subset of the vertical tube $B^{v}$ represented by the set 
\[\Lambda:=\left\{|x|\leq r\right\} \cap H^{n}\left(B^{v}_{\lambda}\cap \left\{|y|\leq \gamma_{1}R\right\}\right).\]  
By Lemma \ref{lemma:polydisks}, $\Lambda$ is biholomorphic to the standard polydisk $\D_{r}\times \D_{\gamma_{1}R}$ via the map $G(x,y)=(x, h^{2}_{-n}(x,y))$. Note that the map $G$ preserves vertical lines, and maps the disk $D^{v}_{w}$ of the critical locus to the graph 
\[
\mathcal{O}:=G(D^{v}_{w})=\{(x, g_{-n}(x,\xi(x))), x\in\D_{r}\}.
\]

Inside the polydisk $\D_{r}\times \D_{\gamma_{1}R}$ we have a still distorted, but more classical Hartogs figure $\mathcal{H}$ formed by taking a neighborhood of the vertical boundary of the polydisk $\{|x|=r, |y|\leq \gamma_{1}R\}$ and of the holomorphic disk $\mathcal{O}$. The holomorphic function $f:=T\circ G:\mathcal{H}\rightarrow \C$. 
extends as a holomorphic function to the entire polydisk $\D_{r}\times \D_{\gamma_{1}R}$, and this extension is unique. 
Tracing back our steps, it means that $g$ has a unique holomorphic extension to the set in equation \eqref{eq:hartogs}, which completes our proof.
\qed

\noindent\textbf{Proof of Theorem \ref{thm:CL}.} 
Consider the polydisk regions $X_{\lambda}$ from \eqref{def:int}. By Theorem \ref{thm: zero}, the sum of the Milnor numbers of the singularities of the critical locus inside each $X_{\lambda}$ satisfies the relation $\mu_{c,a}-\xi_{c,a}=0$ for all $(c,a)\in \HOV_{\beta}$. 
In Theorem \ref{thm:handle}, using the general theory of analytic sets and their singularities, we showed that if the critical locus in $X_{\lambda}$ is not smooth and $\mu_{c,a}-\xi_{c,a}=0$, then it has a very rigid description: necessarily it is the union of two holomorphic disks intersecting at one point. However, in Lemma \ref{lemma:NOsingN}, we use the dynamics of the H\'enon map to show that the critical locus inside $X_{\lambda}$  cannot be the union of two holomorphic disks intersecting at one point. Therefore the critical locus in $X_{\lambda}$ is smooth. By Theorem \ref{thm:FundamentalDomain}, it follows that the critical locus is smooth in the entire fundamental domain, therefore it is everywhere smooth.

We can in fact describe more accurately the critical locus $\mathcal{C}_{X_{\lambda}}$. By the analysis above, it is a Riemann surface with zero Euler characteristic.

By \cite{F}, in the perturbative setting (i.e. when the Jacobian $a$ is very small), the critical locus inside $X_{\lambda}$ is a connected sum of two disks, hence a cylinder. In \cite{F} the case $a=0$ is first analyzed separately (when $a=0$ the H\'enon map $H_{c,a}$ degenerates to the quadratic polynomial $p_{c}$,  and its critical locus inside $X_{\lambda}$ is non-smooth, and is the union of two holomorphic disks intersecting transversely at one point). Then, bifurcation theory is used to claim that  the holomorphic perturbation of two disks intersecting at one point is either smooth (in which case it is a cylinder) or it has singularities (in which case it is again the intersection of two disks). 

By the Ehresmann Fibration Theorem, adapted to the H\'enon map in Theorem  \ref{thm: zero}, all smooth critical loci in the HOV$_{\beta}$ region are diffeomorphic. 
Therefore, for each $(c,a)$ in HOV$_{\beta}$, the critical locus of $H_{c,a}$ in each set $X_{\lambda}$ is a connected sum of two disks, hence a cylinder. It is easy to see now that the entire critical locus $\Crit$ is connected, since the critical loci $\mathcal{C}^{v}$ and $\mathcal{C}^{h}$ from the tubes $S^{v}$ and $S^{h}$ are connected by a cylinder inside the polydisk $B$.

We have shown that $\mathcal{C}$ has no singularities, so all points in $\mathcal{C}$ are regular points. Moreover $\mathcal{C}$ is connected. An analytic variety is irreducible if and only if the set of regular points is connected. Hence $\mathcal{C}$ is irreducible.
\qed

\section{Model of the critical locus}\label{sec:criticalModel}

In this section we construct the truncated sphere model of the critical locus and prove Theorem \ref{thm:SM}.

We will consider the fundamental region of the critical locus given in Equation \ref{def: thirdFundamentalRegion}.
 In Theorem \ref{prop:critlocus-Sdelta} we understood the critical loci $\mathcal{C}^{v}$ and $\mathcal{C}^{h}$ from the two large tubular regions $S^{v}$ around the $y$-axis and $S^{h}$ around the $x$-axis. 
These parts of the critical locus are naturally glued together by dynamics, since the H\'enon function maps $S^{v}\cap\widetilde{V}$ outside $\widetilde{V}$ and onto $S^{h}\cap H(\widetilde{V})$. Hence $H(\mathcal{C}^{v})$ is the unique irreducible component of the critical locus which extends $\mathcal{C}^{h}$ inside $H(\widetilde{V})$.

\begin{figure}[htb]
\begin{center}
\includegraphics[scale=.9]{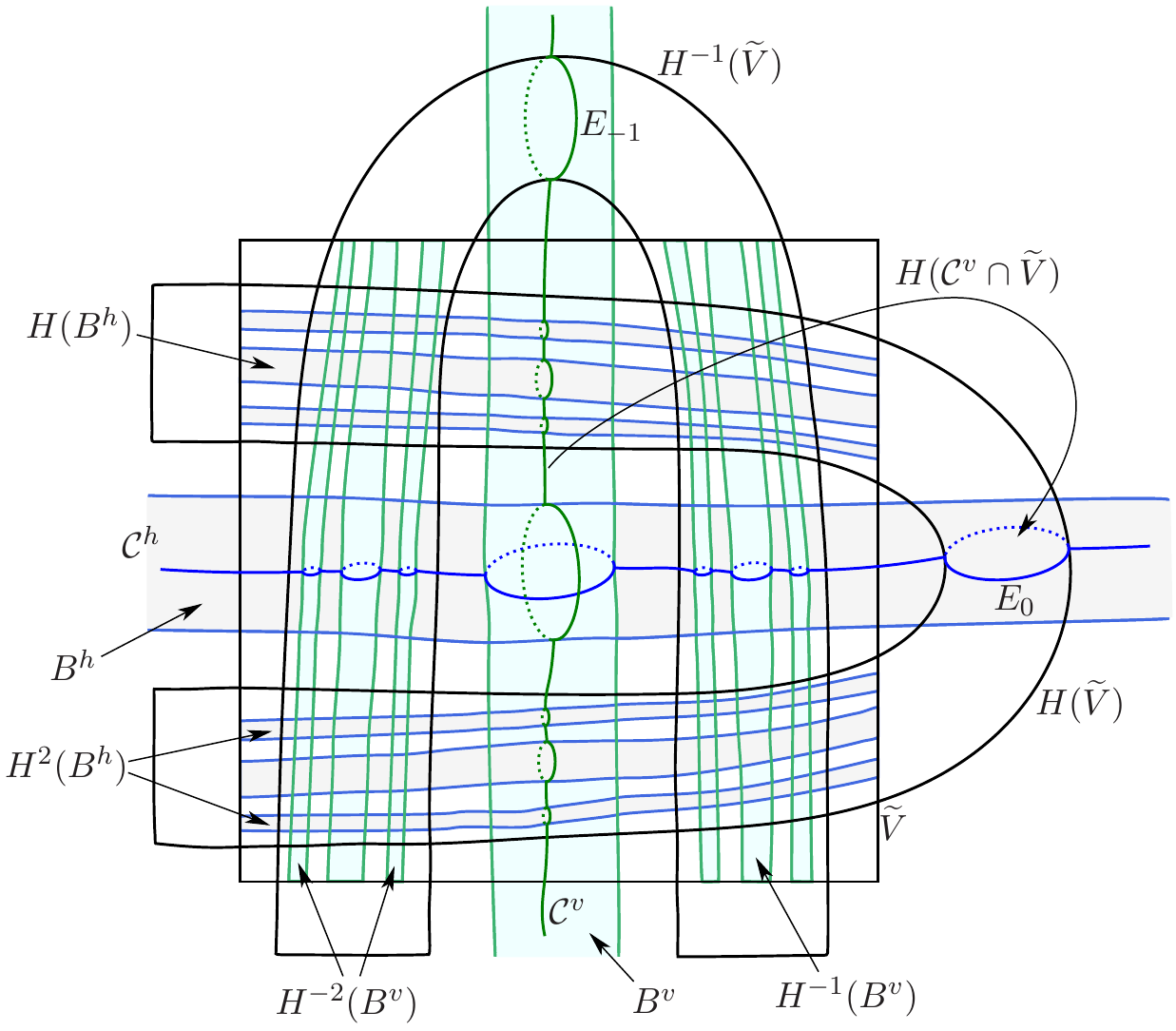}
\end{center}
\caption{The components $\mathcal{C}^{v}$, $\mathcal{C}^{h}$ of the critical locus, and the glueing of $H(\mathcal{C}^{v}\cap\widetilde{V})$ and $\mathcal{C}^{h} $ on the boundary of $H(\widetilde{V})$.}
\label{fig:components}
\end{figure}

Recall from Theorem \ref{prop:critlocus-Sdelta} that $\mathcal{C}^{h}$ is a punctured disk tangent at infinity to the $x$-axis, with a Cantor set removed (corresponding to the intersection of the closure of $C^{h}$ with $J^{+}$) and with punctures at the vertical boundaries of each of the sets $H^{-n}(B^v), n\geq 0$ and at the boundary of $H(\widetilde{V})$, as depicted in Figure \ref{fig:components}.  By compactifying $\mathcal{C}^{h}$ with a point at $\infty$ on the $x$-axis, we can view it as a subset of the lower hemisphere of the unit sphere in $\C^{2}$. The equator corresponds to the boundary circle $E_{0}$ given by the intersection of $\overline{\mathcal{C}^{h}}$ with the boundary of $H(\widetilde{V})$, 
\[E_{0}=H(\partial\mathcal{C}^{v}\cap \partial\widetilde{V})= \partial\mathcal{C}^{h}\cap \partial H(\widetilde{V}).\]
To model the critical locus $\mathcal{C}^{h}$, from the hemisphere we remove a Cantor set corresponding to the intersection of $\overline{\mathcal{C}^{h}}$ with $J^{+}$, a set of disks corresponding to $H^{-n}(B^v), n\geq 0$, and an extra point representing $\infty$. Since the resulting model is merely topological, we will later depict the removed $\infty$ point on the equator.

Also by Theorem \ref{prop:critlocus-Sdelta}, we know that $\mathcal{C}^{v}\cap \widetilde{V}$ is a  punctured holomorphic disk, which intersects the horizontal boundary of $\widetilde{V}$ transversely in a set homeomorphic to a circle, with a Cantor set removed (corresponding to the Julia set $J^{-}$) and with punctures at the horizontal boundaries of each of the sets $H^{n}(B^h)$, $n\geq 0$. Hence we can view $H(\mathcal{C}^{v}\cap \widetilde{V})$ as a subset of the upper hemisphere of the unit sphere in $\C^{2}$. 

To model the critical locus $H(\mathcal{C}^{v}\cap\widetilde{V})$, from the lower hemisphere we shall remove a Cantor set corresponding to the intersection of the closure of $\mathcal{C}^{v}$ with $J^{-}$ and a disk for each set $H^{-n}(B^v), n\geq 0$.

\begin{figure}[htb]
\begin{center}
\vglue -1cm
\begin{tikzpicture}
\node[] (image) at (0,2) {\includegraphics[scale=0.35, angle=30]{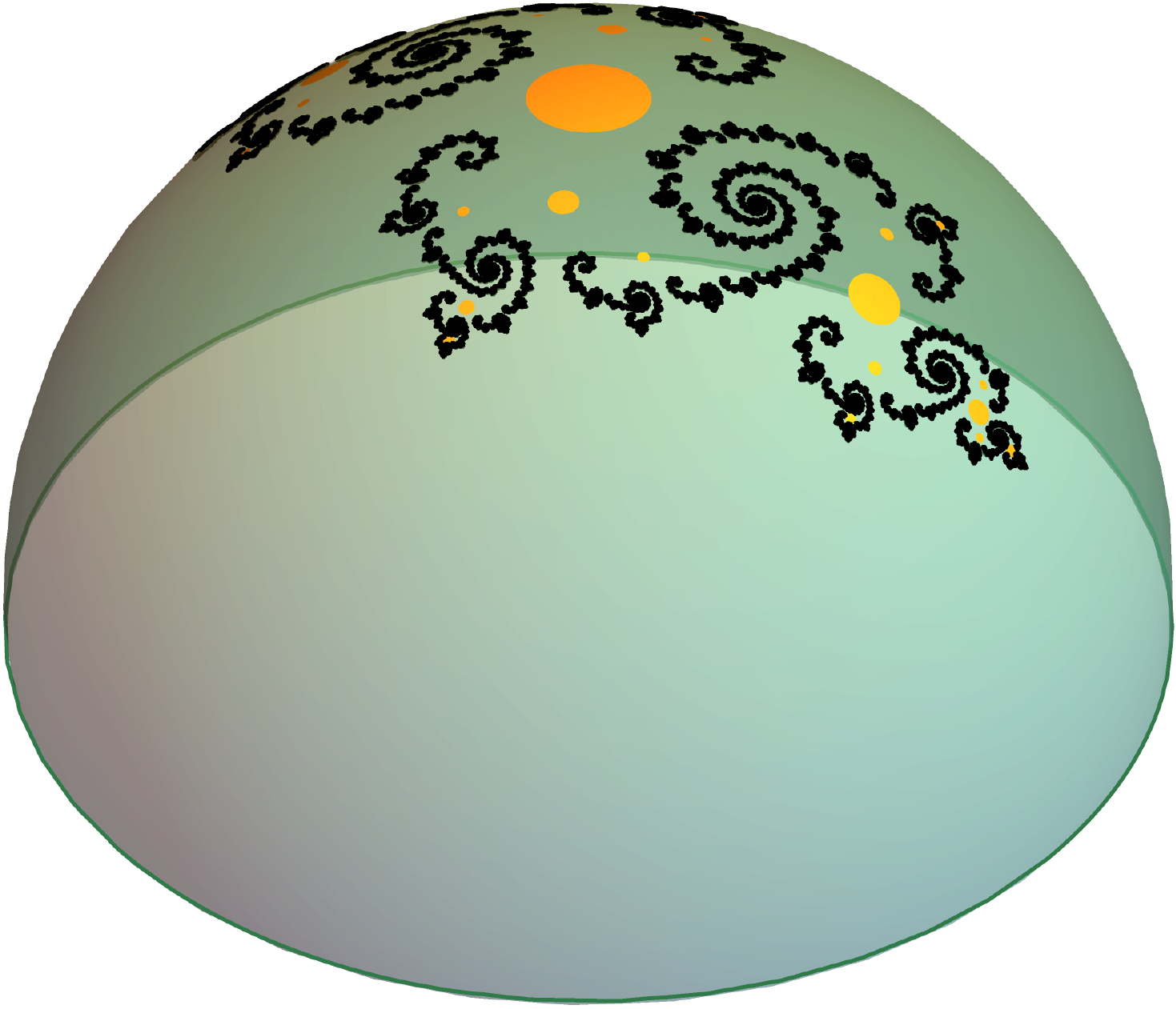}};
\node[] (image) at (5.25,0) {\includegraphics[scale=0.35, angle=30]{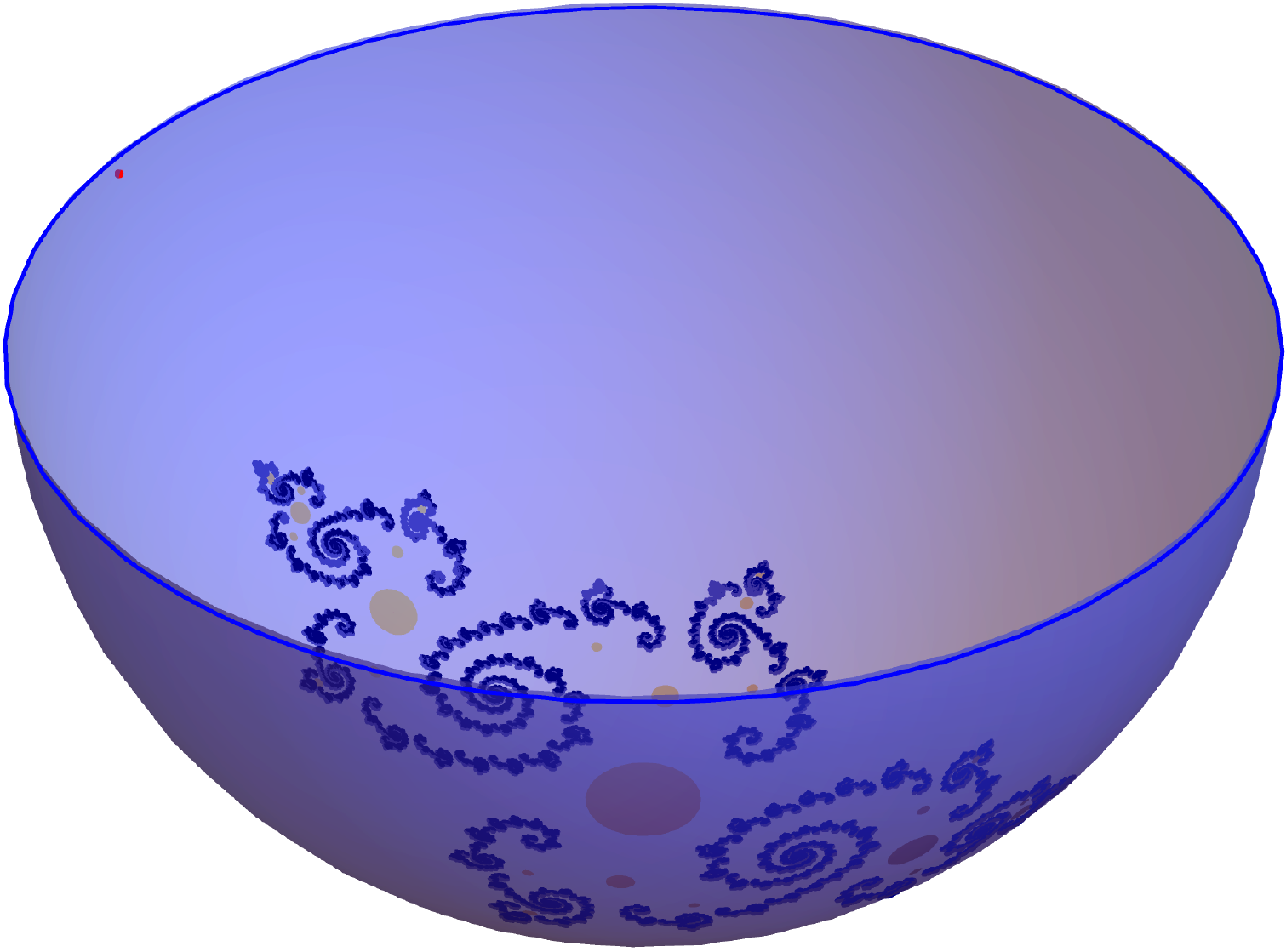}};
\end{tikzpicture}
\end{center}
\vglue -1cm
\caption{ {\sc right:} The critical locus $\mathcal{C}^{h}$ is on the lower hemisphere. {\sc left:} The critical locus $H(\mathcal{C}^{v})\cap H(\widetilde{V})$ is on the upper hemishpere.}
\label{fig:hemishperes}
\end{figure}

We can glue the lower and the upper hemisphere depicted in Figure \ref{fig:hemishperes} along the equator. The glued object is a truncated sphere, homeomorphic to the critical locus $\mathcal{C}^{h}\cup H(\mathcal{C}^{v}\cap \widetilde{V})$. This truncated sphere represents the building block $\mathcal{S}_{0}$ of the truncated sphere model described in Theorem  \ref{thm:SM}.

It is easy to see on the truncated sphere $\mathcal{S}_{0}$ in Figure \ref{fig:hemishperes} that the Cantor set removed from the lower hemisphere is the accessible boundary of the critical locus $\mathcal{C}^{h}$ which lies in $J^{+}$. However note that $\mathcal{C}^{h}$ union this Cantor set is not a holomorphic punctured disk, in fact it is not even a topological manifold, because the removed yellow disks accumulate on this Cantor set, therefore at no point in the Cantor set can we find a neighborhood $U$ such that the intersection of $U$ with the sphere is homeomorphic to a disk in the plane.  
Similarly, the Cantor set in the upper hemisphere represents the accessible boundary of $H(\mathcal{C}^{v})$ which lies in $J^{-}$, but it cannot be added to the sphere to generate a larger Riemann surface.

For the sake of symmetry, we can also apply Theorem \ref{prop:critlocus-Sdelta} to glue $\mathcal{C}^{v}$ and  $H^{-1}(\mathcal{C}^{h}\cap \widetilde{V})$. This time however, we will put $\mathcal{C}^{v}$ on the upper hemisphere and $H^{-1}(\mathcal{C}^{h}\cap \widetilde{V})$ on the lower hemisphere. The advantage of making this choice is that we will have on the upper hemispheres on the truncated spheres only the Cantor sets corresponding to the intersection of the critical locus with the Julia set $J^{-}$. We obtain another truncated sphere, which corresponds to $\mathcal{S}_{-1}$ in the model of Theorem  \ref{thm:SM}. 

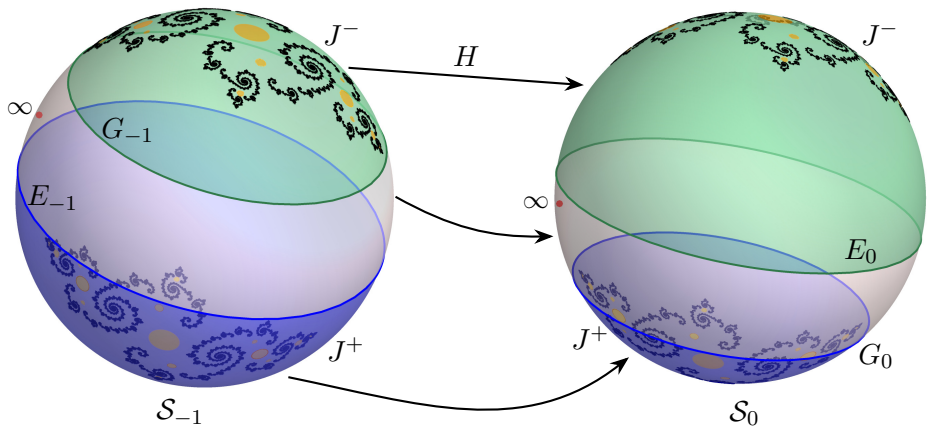
\begin{figure}[htb]
\begin{center}
\begin{tikzpicture}
\node[] (image) at (0,0) {\includegraphics[scale=.53]{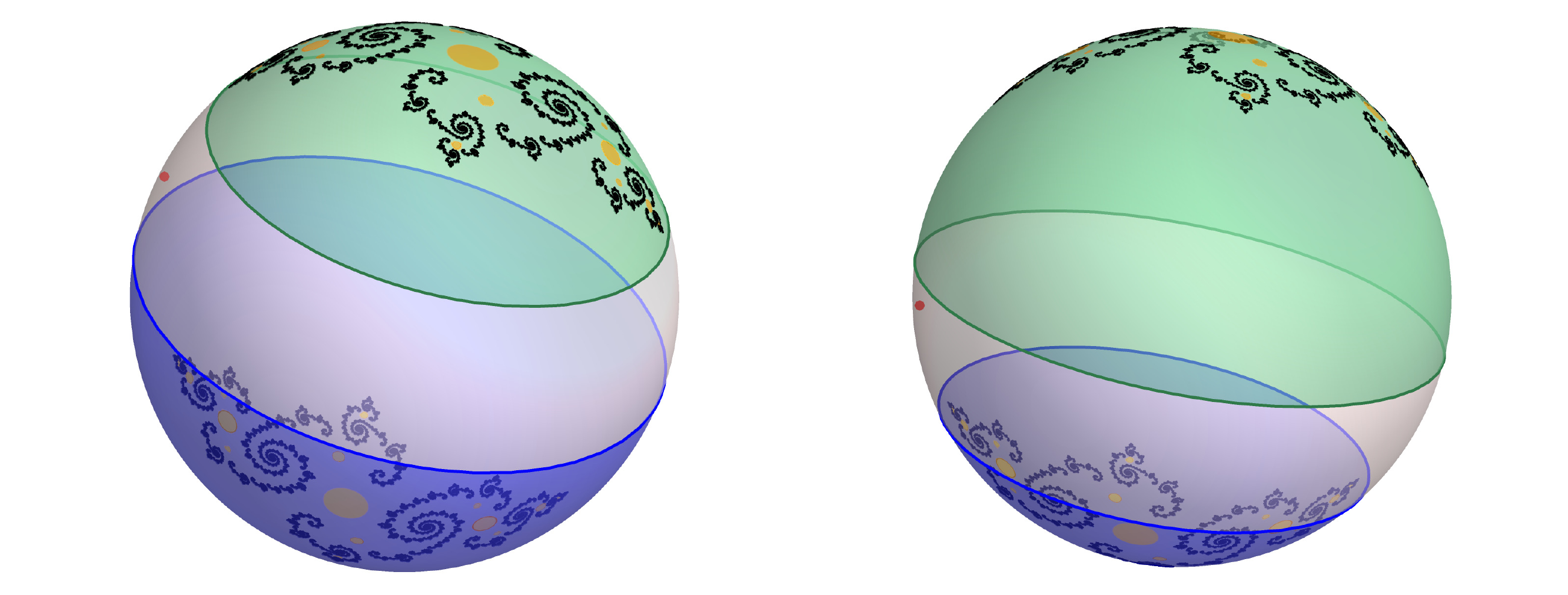}};
\node[] (a) at (-1.65, 2.2) {$J^{-}$};
\node[] (aa) at (5.5, 2.2) {$J^{-}$};
\node[] (b) at (-1.7 , 1.75) {};
\node[] (c) at (1.7, 1.5) {};
\node[] (d) at (-4.5 , 0.9) {$G_{-1}$};
\node[] (dd) at (5.25 , -0.7) {$E_{0}$};
\node[] (e) at (-5.5 , 0) {$E_{-1}$};
\node[] (ee) at (5.45 , -2.1) {$G_{0}$};
\node[] (f) at (-1.6, -2) {$J^{+}$};
\node[] (ff) at (1.65, -1.85) {$J^{+}$};
\node[] (g) at (-3.8 , -2.85) {$\mathcal{S}_{-1}$};
\node[] (gg) at (3.7 , -2.85) {$\mathcal{S}_{0}$};
\node[] (h) at (-2.5, -2.35) {};
\node[] (i) at (2.3, -2) {};
\node[] (j) at (-1.075, 0.1) {};
\node[] (k) at (1.295, -0.5) {};
\node[] (l) at (-5.93, 1.2) {\small$\infty$};
\node[] (m) at (0.93, -0.05) {\small$\infty$};
\draw[-{Stealth[scale=1.1]}, thick] (b)--(c) node[midway, above]{$H$};
\path[-{Stealth[scale=1.1]}, thick] (h) edge[bend right=0,  out=345, in=215] node [left] {} (i);
\path[-{Stealth[scale=1.1]}, thick] (j) edge[bend right=15] node [left] {} (k);
\end{tikzpicture}
\end{center}
\caption{The dynamics of the model map $H$ acting between the spheres $\mathcal{S}_{-1}$ and $\mathcal{S}_{0}$ of the model of the critical locus.}
\label{fig:spheredynamics}
\end{figure}

It is easy to see that the two truncated spheres $\mathcal{S}_{-1}$ and $\mathcal{S}_{0}$ are homeomorphic. The model H\'enon map takes the truncated sphere $\mathcal{S}_{-1}$ onto the truncated sphere $\mathcal{S}_{0}$.  The induced model map $H$ between $\mathcal{S}_{-1}$ and $\mathcal{S}_{0}$ 
extends continuously to the corresponding Cantor sets on the two truncated spheres (see Figure \ref{fig:spheredynamics}).
First recall that $\mathcal{C}^{v}$ and $\mathcal{C}^{h}$ intersect the boundary of $\widetilde{V}$ transversely, since the boundary is foliated by leaves of the foliation of $U^{\pm}$ and the critical locus is transverse to both foliations. Denote by  $G_{-1}$ the circle given by the intersection of $\mathcal{C}^{v}$ with the horizontal boundary of the set $\widetilde{V}$, and by $G_{0}$ the intersection of $\mathcal{C}^{h}$ with the vertical boundary of $\widetilde{V}$.  $G_{-1}$ is a subset of the critical locus $\mathcal{C}^{v}$, therefore it corresponds to a simple closed curve  inside the upper hemisphere of $\mathcal{S}_{-1}$, surrounding the Cantor set $J^{-}$. Similarly, but in the opposite direction, the curve $G_{0}$ is a subset of the critical locus $\mathcal{C}^{h}$, therefore it corresponds in the truncated sphere model to a simple closed curve inside the lower hemisphere of the sphere $\mathcal{S}_{0}$, surrounding the Cantor set $J^{+}$. 
 By construction, $G_{-1}$ is mapped by the model map $H$ to the equator $E_{0}$ of the sphere $\mathcal{S}_{0}$, whereas the equator  $E_{-1}$ of the sphere $\mathcal{S}_{-1}$ maps to the curve $G_{0}$ on the sphere $\mathcal{S}_{0}$. 
Hence, $H$ maps the lower blue hemisphere of the sphere $\mathcal{S}_{-1}$ strictly inside the lower hemisphere of sphere $\mathcal{S}_{0}$. 
Likewise, the inverse model map $H^{-1}$ maps the upper green hemisphere of sphere $\mathcal{S}_{0}$ strictly inside the upper hemisphere of $\mathcal{S}_{-1}$.

\begin{figure}[htb]
\begin{center}
\includegraphics[scale=.63]{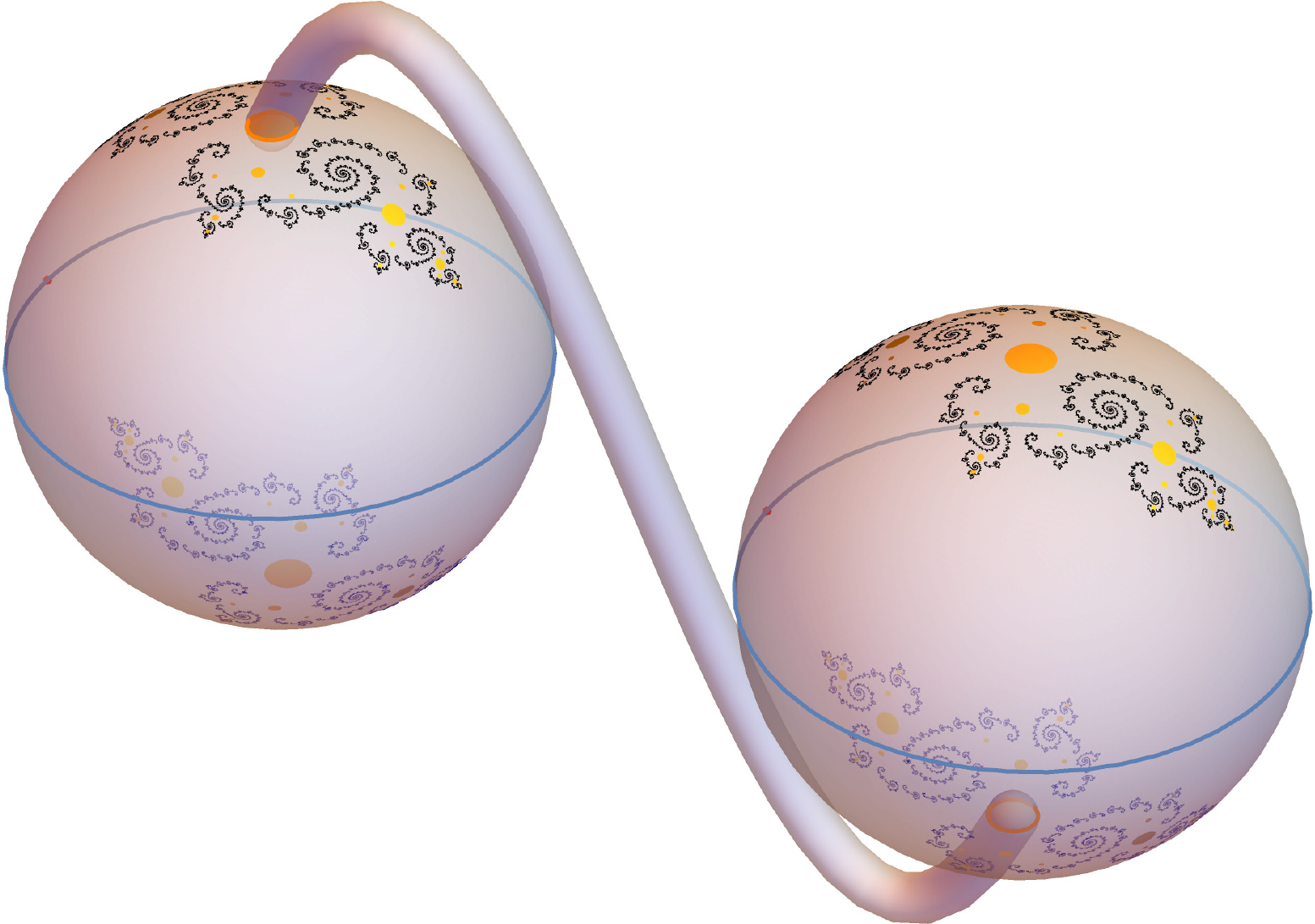}
\end{center}
\caption{The truncated spheres $\mathcal{S}_{-1}$ and $S_{0}$ are connected by one handle.}
\label{fig:onehandle}
\end{figure}

By the proof of Theorem \ref{thm:CL}, we know that there exists a unique irreducible component of the critical locus in the polydisk $B=B^{h}\cap B^{v}$ which extends $\mathcal{C}^{v}$ and $\mathcal{C}^{h}$ inside $B$, and moreover this component is homeomorphic to a cylinder. It has two boundary circles, one on the horizontal boundary of the polydisk (and on $\mathcal{C}^{v}$) and one on the vertical boundary of the polydisk (and on $\mathcal{C}^{h}$). We can model it by drawing a handle between the truncated spheres $\mathcal{S}_{-1}$ and $\mathcal{S}_{0}$, as in Figure \ref{fig:onehandle}; this handle is a cylinder connecting the biggest removed disk (shown in yellow) on the upper hemisphere of the truncated sphere $\mathcal{S}_{-1}$ to the biggest removed disk (shown in yellow) on the lower hemisphere of the truncated sphere $\mathcal{S}_{0}$. 

\begin{figure}[htb]
\begin{center}
\includegraphics[scale=0.95]{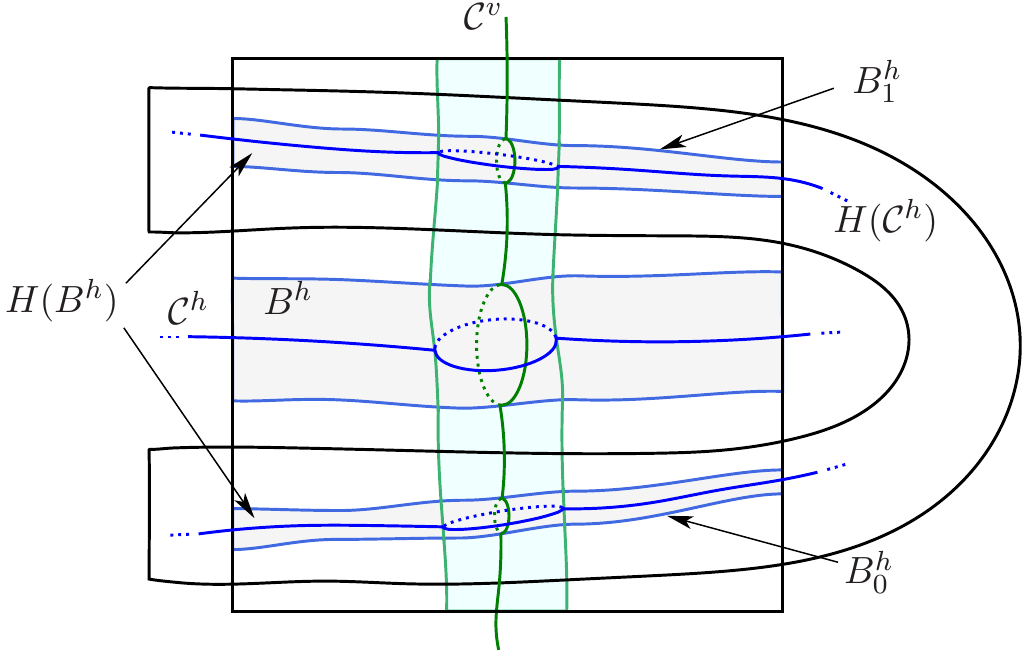}
\end{center}
\caption{The component $\mathcal{C}^{h,1}=H(\mathcal{C}^{h})$ of the critical locus, and its intersection with the vertical tube $B^{v}$.}
\label{fig:HofCritical}
\end{figure}

The proof of Theorem \ref{thm:CL} gives that all critical loci $\mathcal{C}_{X_{\lambda}}$ (or $\mathcal{C}_{Y_{w}}$ if working in the vertical tube) are homeomorphic to cylinders. It follows that the truncated sphere $\mathcal{S}_{-1}$ is connected by $2^{n}$ handles with the truncated sphere $\mathcal{S}_{n}$, $n\geq 0$. These handles connect the yellow disks on level $n$ in the upper hemisphere of sphere $\mathcal{S}_{-1}$ to the matching yellow disks in the lower hemisphere of sphere $\mathcal{S}_{n}$. It is perhaps instructive to illustrate where the two handles which connect  $\mathcal{S}_{-1}$ to $\mathcal{S}_{1}$ come from. 
We know that there exists a unique irreducible component of the critical locus in the polydisks $Y_{1}=B^{h}_{1}\cap B^{v}$ and respectively in $Y_{0}=B^{h}_{0}\cap B^{v}$ which extends $\mathcal{C}^{v}$ and $C^{h,1}:=H(\mathcal{C}^{h})$ inside $Y_{1}$ and $Y_{0}$ (see Figure \ref{fig:HofCritical}).

\begin{figure}[htb]
\begin{center}
\includegraphics[scale=.55]{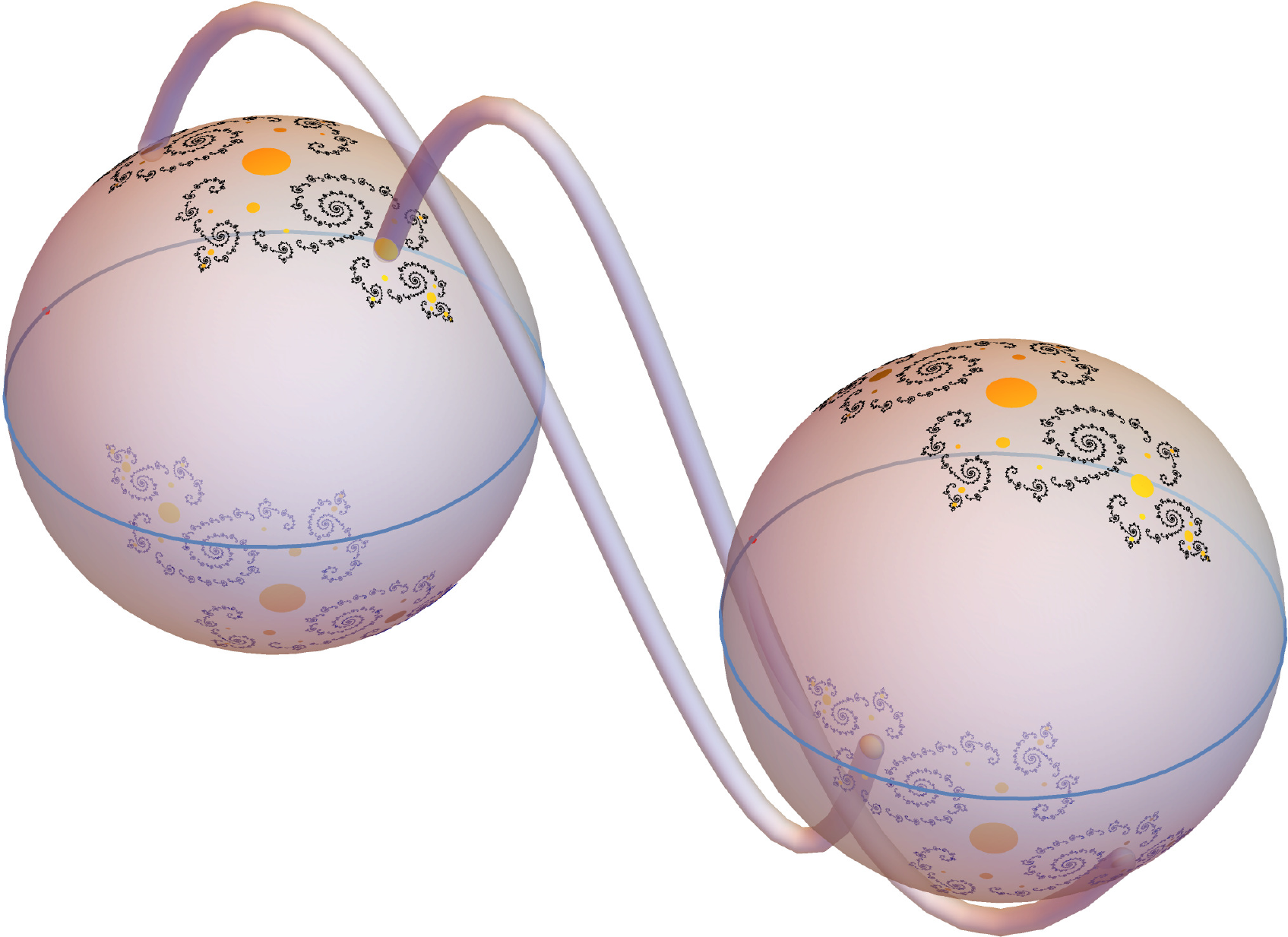}
\end{center}
\caption{$\mathcal{S}_{-1}$ and $\mathcal{S}_{1}$ are connected by two handles.}
\label{fig:twohandles}
\end{figure}

We can model it by drawing two handles between the truncated spheres $\mathcal{S}_{-1}$ and $\mathcal{S}_{1}$, as in Figure \ref{fig:twohandles}; these handles are cylinders connecting the two second biggest yellow disks on the upper hemisphere of the truncated sphere $\mathcal{S}_{-1}$ to the two second biggest yellow disks on the lower hemisphere of the truncated sphere $\mathcal{S}_{1}$.

Therefore, the critical locus in the fundamental region
\begin{eqnarray*}
\mathcal{C}^{v} \cup H^{-1}(\mathcal{C}^{h}\cap \widetilde{V}))\cup \bigcup_{w\in \{0,1\}^{n},\, n\geq 0} \Crit_{Y_{w}}
\end{eqnarray*}
is represented in the truncated sphere model by the  sphere $\mathcal{S}_{-1}$ together with all the handles that come out of the upper hemisphere of sphere $\mathcal{S}_{-1}$ (these handles join $\mathcal{S}_{-1}$ to the lower hemispheres of the  spheres $\mathcal{S}_{n}$, $n\geq 0$).

Similarly, the critical locus in the fundamental region
\begin{eqnarray*}
H(\mathcal{C}^{v}\cap \widetilde{V}) \cup \mathcal{C}^{h}\cup \bigcup_{\lambda\in \{0,1\}^{n},\, n\geq 0} \Crit_{X_{\lambda}}
\end{eqnarray*}
corresponds in the truncated sphere model to the  sphere $\mathcal{S}_{0}$ together with all the handles that come out of the lower hemisphere of sphere $\mathcal{S}_{0}$ (these handles join $\mathcal{S}_{0}$ to the upper hemispheres of the  spheres $\mathcal{S}_{n}$, $n<0$).

\section{Relation between the critical loci $\mathcal{C}$, $\mathcal{C}^s$ and $\mathcal{C}^u$}\label{sec:criticalLoci}

Recall that $\mathcal{C}^s$ represents the set of tangencies between the foliation of $U^-$ and the lamination of $J^+$, and respectively  $\mathcal{C}^u$ is the set of tangencies between the foliation of $U^+$ and the lamination of $J^-$. In this section we will discuss the relation between the various critical loci and prove Theorem \ref{thm:Acc}.

As we know, the critical locus $\mathcal{C}$ is 
an analytic subvariety of
$U^{+}\cap U^{-}$, so it is closed in $U^{+}\cap U^{-}$. Therefore, when we regard $\mathcal{C}$ as a subset of $\mathbb{C}^{2}$, its boundary must be contained in the boundary of $U^{+}\cap U^{-}$, that is in $J^{+}\cup J^{-}$. In fact, by \cite{BS5}, the closure of the critical locus $\mathcal{C}$ always intersects both $J^+$ and
$J^-$, that is 
\[\overline{\mathcal{C}}\cap J^+\cap U^-\neq
\emptyset\ \ \mbox{and}\ \ \overline{\mathcal{C}}\cap J^-\cap U^+\neq
\emptyset.
\] 
However, it is not true that $\mathcal{C}^s=\partial{\mathcal{C}}\cap (J^+\cap U^-)$ and $\mathcal{C}^u=\partial{\mathcal{C}}\cap (J^-\cap U^+)$. The relation between $\mathcal{C}$, $\mathcal{C}^s$ and $\mathcal{C}^u$ is  rather mysterious is general.

In this section, we will describe the relation between the three critical loci in the cases where we understand the model for the critical locus $\mathcal{C}$. Note that in all these cases, the H\'enon map is hyperbolic, so the laminations of $J^{+}$ and $J^{-}$ are everywhere transverse to each other.

In the case studied in \cite{LR}, \cite{T}, the Julia set is connected, and the foliation of $U^+$ and the lamination of $J^+$ fit together continuously. The boundary of $\mathcal{C}_0$, the primary component of the critical locus described in Theorem \ref{thm: LyubichRobertson}, belongs to $J^+$ and the stable critical locus $\mathcal{C}^s$ is the union of the boundaries of the forward and backward iterates of $\mathcal{C}_0$, that is,
\[
 \mathcal{C}^s =\bigcup_{n\in\Z} \partial{H^n(\mathcal{C}_0)}.
\]
Since $J$ is connected, it follows from \cite{BS5} that $\mathcal{C}^{u}=\emptyset$. However, this does not mean that the closure of $\mathcal{C}$ is disjoint from $J^{-}$, in fact the forward iterates of $\overline{\mathcal{C}_0}$ accumulate on $J^-\cap (U^+\cup J^+)$. 

In the case studied in \cite{F}, as well as in the general horseshoe region HOV$_{\beta}$ studied in this paper, the accessible boundary of the critical locus is $\Cs\cup\Cu$. 

\begin{lemma}\label{lemma:csu} Consider the regions $S^{v}$ and $S^{h}$ defined in Equations \eqref{def:S1} and \eqref{def:S2}. Denote by $\mathcal{C}^{hs}$ and $\mathcal{C}^{hu}$ the stable and unstable critical loci in the closure $\overline{S^{h}}$.  Denote by $\mathcal{C}^{vs}$ and $\mathcal{C}^{vu}$ the stable and unstable critical loci in $\overline{S^{v}}$. Then $\mathcal{C}^{vs}=\emptyset$ and  $\mathcal{C}^{vu}$ is a Cantor set, whereas $\mathcal{C}^{hu}=\emptyset$ and $\mathcal{C}^{hs}$ is a Cantor set.
\end{lemma}
\proof This is an easy Corollary of Theorem \ref{prop:critlocus-Sdelta}, using the horseshoe construction.
Note first that the closure of the set $S^{v}$ is included in the escaping set $U^{+}$, and therefore disjoint from the Julia set $J^{+}$. Hence $\mathcal{C}^{vs}=\emptyset$. Similarly, the closure of $S^{h}$ is included in the escaping set $U^{-}$, and therefore disjoint from the Julia set $J^{-}$. Hence $\mathcal{C}^{hu}=\emptyset$. 

The critical locus $\overline{C^{v}}$ is transverse to the foliation of $U^{-}$ and implicitly to the horizontal boundaries of the sets $H^{n}(\widetilde{V})$, for each $n\geq 0$. By the dynamical properties of the set $\widetilde{V}$ we have $J^{-}\cap \widetilde{V} = \bigcap_{n\geq 0} H^{n}(\widetilde{V})$. Furthermore, in the horseshoe region, we know that $J^{-}\cap S^{v}_{\delta}$ is isomorphic to the direct product of the disk $\D_{\delta}$ with a Cantor set $K_{1}$. By Theorem \ref{prop:critlocus-Sdelta}, the critical locus $\overline{\mathcal{C}^{v}}$ is trapped inside the straight tube $S^{v}_{\delta}$. Since the foliation of $U^{-}$ and the lamination of $J^{-}$ in $S^{v}$ consist of long horizontal-like leaves which fit together continuously to give a locally trivial lamination of the set $S^{v}$ (see Corollary \ref{cor:fitcontinuously}), by transversality, it follows that $\mathcal{C}^{vu}=\overline{C^{v}}\cap J^{-}$. Moreover, we obtain that $\mathcal{C}^{vu}$ is isomorphic to the Cantor set $K_{1}$, and is contained in the accessible boundary of $C^{v}$.
Note that the other parts of the accessible boundary of $C^{v}$ consist of the real analytic circles represented by the intersection of $C^{v}$ with the horizontal boundaries of the sets $H^{-1}(\widetilde{V})$ and $H^{n}(B^{h})$, $n\geq 0$. These circles are neither in $J^{-}$, nor part of the boundary of the big critical locus $\mathcal{C}$.

An identical analysis in the set $S^{h}$ leads to the description of $\mathcal{C}^{hs}$.
\qed
\begin{remark} Compare the asymmetry in Lemma \ref{lemma:csu} to Theorem \ref{thm:FundamentalDomain} and Corollary \ref{cor:fund}: both $F\cap \mathcal{C}$ and $F'\cap\mathcal{C}$ are fundamental domains of the critical locus $\mathcal{C}$, but $F\cap \mathcal{C}^{s}$ is a fundamental domain for $\mathcal{C}^{s}$, and $F' \cap \mathcal{C}^{u}$ is a fundamental domain for $\mathcal{C}^{u}$.
\end{remark}

\begin{lemma}\label{lemma:J} The Julia set $J$ is included in the closure of the critical locus $\mathcal{C}$.
\end{lemma}
\proof Let $\tau$ be the conjugacy map in Equation \eqref{def:conjugacy} between the shift map on $\Sigma_{2}$ and the H\'enon map on its Julia set $J$. Fix any point $z\in J$, and let $\omega=\tau^{-1}(z)\in \Sigma_{2}$. For any $m,n>0$, let $\lambda=\omega_{-n}\ldots \omega_{-1}$ and $w=\omega_{0}\ldots \omega_{m-1}$ be two finite substrings of $\omega$, and recall the binary coding in Section \ref{sec:DynamicalBox}, Equation \eqref{eq: poly}.  The critical locus $C$ has non-empty intersection with the set $\Delta_{(\lambda,w)}$. However, as $n,m\rightarrow\infty$, the size of the sets $\Delta_{(\lambda,w)}$ converges to $0$, and
\[
\mathcal{C}\cap \Delta_{(\lambda,w)} \subset \bigcap_{k=-n}^{m-1}H^{k}(T^{h}_{w_{k}})\rightarrow \bigcap_{k=-\infty}^{\infty}H^{k}(T^{h}_{w_{k}})= \tau(\omega)=z\in J.
\]
Thus any point in the Julia set belongs to the closure of the critical locus. 
 \qed
 
The H\'enon map in the HOV region is hyperbolic on its Julia set, therefore all intersections between the stable and unstable manifolds of points in $J$ are transverse intersections. However, the Julia set is accumulated by tangencies between the foliations of the escaping sets. This goes to show that while there are regions in $\widetilde{V}$ where the foliations of $U^{\pm}$ fit together continuously with the laminations of $J^{\pm}$, this does not happen in the neighborhood of any point in the Julia set $J$.

\medskip
\noindent{\textbf{Proof of Theorem \ref{thm:Acc}.}} First observe that all regions $\Delta_{(\lambda,w)}$ from $\widetilde{V}$ (which contain the handles in the truncated sphere model) are disjoint from $J^{+}\cup J^{-}$, hence \[\mathcal{C}^{u}\cap \Delta_{(\lambda,w)} = \mathcal{C}^{s}\cap\Delta_{(\lambda,w)}=\emptyset. \]
In Lemma \ref{lemma:csu} we showed that the accessible boundary of the critical locus $\mathcal{C}^{h}$ is the stable critical locus $\mathcal{C}^{hs}$, while the accessible boundary of the critical locus $\mathcal{C}^{v}$ is the unstable critical locus $\mathcal{C}^{vu}$. It follows from the truncated sphere model of the critical locus from Theorem \ref{thm:SM}   that
\begin{equation*}
\mathcal{C}^{s} = \bigcup_{n=-\infty}^{\infty}H^{n}(\mathcal{C}^{hs}) \ \  \mbox{and} \ \
\mathcal{C}^{u} = \bigcup_{n=-\infty}^{\infty}H^{n}(\mathcal{C}^{vu}), 
\end{equation*}
which implies that the accessible boundary of the critical locus $\mathcal{C}$ is $\mathcal{C}^{s}\cap \mathcal{C}^{u}$. 

To see that $\partial\mathcal{C}=J^{+}\cup J^{-}$, it suffices to prove this statement inside the set $\widetilde{V}$. Indeed, by construction $J\subset \widetilde{V}$ and by hyperbolicity $J^{+}=W^{s}(J)$ and $J^{-}=W^{u}(J)$, so  
\[
J^{+}=\bigcup_{n\geq 0} H^{-n}(J^{+}\cap \widetilde{V})\ \ \mbox{and} \ \ J^{-}=\bigcup_{n\geq 0} H^{n}(J^{-}\cap \widetilde{V}).
\]
Moreover, due to the symmetry between $J^{+}$ and $J^{-}$ in the HOV$_{\beta}$ region, it suffices to prove this statement for $J^{+}$. Pick a point $x\in J^{+}\cap \widetilde{V}$. Denote by $J^{+}_{\lambda}$, $\lambda\in \Sigma_{2}^{-}$, the leaf of the lamination of $J^{+}$ inside $\widetilde{V}$ passing through $x$. For each integer $n>0$, denote by $\bar{\lambda}^{n}$ the finite substrings of length $n$ of $\lambda$ of the form $\lambda_{-n}\ldots \lambda_{-1}$. By the horseshoe condition, the leaf $J^{+}_{\lambda}$ can be obtained as the nested intersection of the thin long vertical-like tubes $T^{v}_{\bar{\lambda}^{n}}$ whose horizontal width decrease to $0$ when $n\rightarrow\infty$. By combining this with Corollary \ref{cor:critlocus-Sdelta}, we get
\begin{equation}\label{eq:limit}
\mathcal{C}^{v}_{\bar{\lambda}^{n}}\subset B^{v}_{\bar{\lambda}^{n}}\subset T^{v}_{\bar{\lambda}^{n}}\rightarrow J^{+}_{\lambda}\ \ \mbox{as}\ \ n\rightarrow\infty. 
\end{equation}
Following the notations of Proposition \ref{prop:longleaves}, $x$ is either in the region $\widetilde{V}_{1}$, in $J^{-}$, or in $\bigcup_{n\geq 0} H^{n}(B^{h})$. In the first case, $x$ belongs to the closure of the union of the backward iterates $\mathcal{C}^{v}_{\bar{\lambda}^{n}}$ of the critical locus $\mathcal{C}^{v}$, as explained in Equation \eqref{eq:limit}. In the second case, $x$ belongs to $J$, so we can apply Lemma \ref{lemma:J}.
In the third case, let $B^{h}_{w}$ be the horizontal tube containing $x$, where $w=w_{0}\ldots w_{m-1}$ is a finite binary word of length $m>0$ or $w=
\emptyset$. Note that in $B^{h}_{w}$ we have a component of the critical locus, namely $\mathcal{C}^{h}_{w}$ whose closure intersects $J^{+}_{\lambda}$ transversely in a single point $y$ (Recall that $y$ is in the accessible boundary of $\mathcal{C}^{h}_{w}$). The polydisks $\Delta_{\bar{\lambda}^{n},w}$ are contained in the tubes $B^{v}_{\bar{\lambda}^{n}}$ and converge uniformly to $J^{+}_{\lambda}\cap B^{h}_{w}$ when $n\rightarrow\infty$. Denote by $\mathcal{E}^{h}_{\bar{\lambda}^{n}, w}$ and $\mathcal{E}^{v}_{\bar{\lambda}^{n}, w}$ the intersections of the critical locus $\mathcal{C}$ with the horizontal and respectively vertical boundary of $\Delta_{\bar{\lambda}^{n},w}$, similar to Lemma \ref{lemma:NOsingN}. As $n\rightarrow\infty$, the circles $\mathcal{E}^{h}_{\bar{\lambda}^{n}, w}$ converge to a real analytic circle $\mathcal{E}^{h}_{\lambda}$ given by the intersection of $J^{+}_{\lambda}$ with the horizontal boundary of $B^{h}_{w}$. The circles 
$\mathcal{E}^{v}_{\bar{\lambda}^{n},w}$ are all contained in the critical locus component $\mathcal{C}^{h}_{w}$ and converge to a single point, namely to the point $y$ in $J^{+}_{\lambda}$. The critical locus in the polydisk-like region $\Delta_{(\bar{\lambda}^{n},w)}$ converges to the interior of the region bounded by $\mathcal{E}^{h}_{\lambda}$ in $J^{+}_{\lambda}$. Hence $J^{+}_{\lambda}$ is included in the closure of the critical locus $\mathcal{C}$. 

This is not surprising, given that the intersections of the critical locus $\mathcal{C}$ with the vertical tubes $T^{v}_{\bar{\lambda}^{n}}$ are purely one-dimensional analytic sets with uniformly bounded one-dimensional volume, which converge geometrically to a a non-empty set $S\subseteq \widetilde{V}$ therefore by Bishop's Compactness Theorem \cite{B} (improved by Watanabe \cite{W}) $S$ is an analytic set, of pure dimension one (see \cite{Chi2}). In this case we have shown that $S=J^{+}_{\lambda}$.

In conclusion $J^{+}\cup J^{-}\subset \overline{\mathcal{C}}$. Since $\partial \mathcal{C}\subset J^{+}\cup J^{-}$, it follows that $\partial \mathcal{C}=J^{+}\cup J^{-}$.
\qed

\section{The monodromy action on the critical locus}

We say that a H\'{e}non map $H_{a,c}$ belongs to the horseshoe locus $\mathcal{HL}$ if the Julia set $J_{a,c}$ is a Cantor set and dynamics on the Julia set is conjugate to the two-sided shift on two-symbols $(\sigma, \Sigma_{2})$. 
The group $Aut(\Sigma_{2}, \sigma)$ of automorphisms of $\Sigma_{2}$ which commute with the shift map is enormous, and it was conjectured by Hubbard that $Aut(\Sigma_{2}, \sigma)$ embeds into the fundamental group of the horseshoe locus (see \cite{BS}). 

Any loop $\gamma:[0,1]\rightarrow \C^{2}$ with base point $(a_{0},c_{0})$  and whose image is contained in the horseshoe locus induces a natural map from the fundamental group $\pi_{1}\left(\mathcal{HL}, (a_{0},c_{0})\right)$ into the group $Aut(\Sigma_{2}, \sigma)$, called the monodromy action. This can be loosely described as a change in the binary labeling of the Cantor Julia set $J_{a_{0},c_{0}}$ when following continuously the points of $J_{a,c}$ along $\gamma$, back to the initial base point $(a_{0},c_{0})$. The motion gives rise to a homeomorphism $f$ of the Julia set $J_{a_{0},c_{0}}$ which commutes with the H\'enon map $H_{a_{0},c_{0}}$ and depends only on the homotopy type of $\gamma$. Following the conjugacy diagram between $(H_{a_{0},c_{0}}, J_{a_{0},c_{0}})$ and $(\sigma, \Sigma_{2})$, one obtains from $f$ a homeomorphism $g$ of $\Sigma_{2}$ which commutes with $\sigma$. The induced map $\gamma\rightarrow g$ is the monodromy action.

Similarly, we can analyze the monodromy action on the critical locus, and make the following remark:

\begin{remark} A loop in the HOV$_{\beta}$ locus induces an isotopy class of maps from the critical locus to itself, which commute with the H\'enon map. 
\end{remark}

The coding that one has on the Cantor Julia set $J_{a,c}$ is a result of the coding of the forward and backward iterates of the set $V$ from the dynamical filtration of $\C^{2}$ using one-sided sequences of $0'$s and $1$'s. It is the same coding that gives the labeling of the bidisks $\Delta_{(\lambda,\omega)}$, which we use to describe the $2^{|\lambda|+|\omega|-1}$ handles between the truncated spheres $\mathcal{S}_{-|\lambda|}$ and $\mathcal{S}_{|\omega|}$ and glue the corresponding disks of the two spheres.

Our parametric region HOV$_{\beta}\subset \rm{HOV}$ sits inside the horseshoe locus and it is easy to see that as subsets of $\C^{2}$, HOV$_{\beta}$ is isomorphic to HOV, therefore isomorphic to $\C^{*}\times (\C-\D)$. As a consequence, the fundamental group of the HOV$_{\beta}$ region is $\Z\times\Z$, with two simple loops as generators: a loop going around the Mandelbrot set, and a loop going around the axis $a=0$. Lipa \cite{L} works with explicit loops in the $(a,c)$-plane, $\gamma_{c}, \gamma_{a}:[0,1]\rightarrow$ HOV, $\gamma_{c}(t)=(1/100, -3e^{2\pi i t})$ and $\gamma_{a}(t)=(1/100 e^{2\pi i t}, -3)$,
which can be enlarged without changing their homotopy type, in order to fit in the HOV$_{\beta}$ region.

\begin{remark} Let $\gamma$ be a simple loop around $a=0$, with base point $(a_{0},c_{0})$ and image contained in the HOV$_{\beta}$ region. $\gamma$ induces a trivial monodromy action on the Julia set $J_{a_{0},c_{0}}\simeq \Sigma_{2}$, therefore it induces just the identity map $Id$ from the critical locus to itself.
\end{remark}

\begin{remark} Let $\gamma$ be a simple loop around the Mandelbrot set, with base point $(a_{0},c_{0})$ and image contained in the HOV$_{\beta}$ region. $\gamma$ induces the monodromy action $\delta$ on the Julia set $J_{a_{0},c_{0}}\simeq \Sigma_{2}$, where $\delta:\Sigma_{2}\rightarrow\Sigma_{2}$ is the homeomorphism which exchanges the symbols $0$ and $1$, that is $\delta(\omega)=\omega'$, where for every $n\in\Z$, $\omega'_{n}=1$ iff $\omega_{n}=0$, and $\omega'_{n}=0$ iff $\omega_{n}=1$. The loop $\gamma$ induces an involution map $\eta$ from the critical locus to itself, which commutes with the H\'enon map. 
\end{remark}
The involution $\eta$ can be described by keeping track of the changes in the labeling of the Julia set $J$ by the monodromy action $\delta$, and seeing how they affect the binary codings relevant for the critical locus. Firstly, changing the coding of the Julia set $J$ by $\delta$ can occur if and only if we change the codings of $J^{+}\cap\widetilde{V}\simeq \Sigma_{2}^{-}$ and $J^{-}\cap\widetilde{V}\simeq \Sigma_{2}^{+}$ by $\delta$. This implies that the finite binary string notations $\lambda, w \in \{0,1\}^{n}$ of all nested neighborhoods $T^{h}_{w}$ of $J^{-}\cap\widetilde{V}$ and  $T^{v}_{\lambda}$ of $J^{+}\cap\widetilde{V}$ have been consistently modified by interchanging the symbols $0$ and $1$. However, this further implies that all objects contained in the tubes $T^{v}_{\lambda}$ and $T^{h}_{w}$ undergo the same change in labeling. Applying the involution $\delta$, we change the labeling of the vertical tubes $B^{v}_{\lambda}$, horizontal tubes $B^{h}_{w}$, polydisks $\Delta_{(\lambda,w)}$, and critical locus components $\mathcal{C}^{v}_{\lambda}$ and $\mathcal{C}^{h}_{w}$. Consequently, the map $\eta$ acts as an involution on each truncated sphere $\mathcal{S}_{n}$. Following the notations of Theorem \ref{thm:SM}, the handle which connects the boundary of the disk $W_{\alpha_{n}}$ on the truncated sphere $\mathcal{S}_{k}$ and the boundary of $U_{\alpha_{n}}$ on $\mathcal{}S_{k+n+1}$ has been interchanged with the handle which connects the boundary of the disk $W_{\delta(\alpha_{n})}$ on the sphere $\mathcal{S}_{k}$ and the boundary of $U_{\delta(\alpha_{n})}$ on $\mathcal{S}_{k+n+1}$.

In fact, it is easy to see the only two automorphisms of $Aut(\Sigma_{2}, \sigma)$ compatible with the labelings of all auxiliary sets in our construction are just $Id$ and $\eta$, since the group $Aut(\Sigma_{2}^{+}, \sigma)$ of automorphisms of the space of one-sided sequences $\Sigma_{2}^{+}$ which commute with the shift $\sigma$ is generated by the homeomorphism $\delta$ which exchanges $0$ and $1$.

\end{document}